\documentclass[10pt,a4paper]{article}
\usepackage[utf8]{inputenc}
\usepackage[spanish]{babel}
\usepackage{amsmath}
\usepackage{amsfonts}
\usepackage{amssymb}
\usepackage{makeidx}
\usepackage{graphicx}
\usepackage{lmodern}
\usepackage{kpfonts}
\usepackage[left=2cm,right=2cm,top=2cm,bottom=2cm]{geometry}
\usepackage{dsfont}
\usepackage{hyperref}
\usepackage{verbatim}
 
\usepackage{fancyhdr}
\usepackage{pgf,tikz}
\decimalpoint

\usepackage{mathrsfs}
\usetikzlibrary{arrows}
\usepackage{subfigure}

\providecommand{\formula }[1]{
\begin{eqnarray*}
#1
\end{eqnarray*}
}

\newtheorem{prop}{Proposición}

\providecommand{\form}[1]{$#1$}

\providecommand{\der}[3]{\left(\dfrac{#1}{#1 #2} \right)^{#3}}

\providecommand{\ifr}[4]{ {}_{#1}^{}#2_{#3}^{#4} }

\providecommand{\ifrR}[5]{ {}^{#1}_{#2}#3_{#4}^{#5} }

\providecommand{\gam}[1]{\Gamma\left( #1\right)}
\providecommand{\si}[1]{\sin \left( #1 \right)}
\providecommand{\co}[1]{\cos \left( #1 \right)}
\providecommand{\ex}[1]{\exp \left( #1 \right)}

\providecommand{\com}[1]{$``$#1$"$}
\providecommand{\conj}[1]{\left\lbrace #1 \right\rbrace}
\providecommand{\re}[1]{Re\left(#1 \right)}

\providecommand{\formula }[1]{
\begin{eqnarray*}
#1
\end{eqnarray*}
}

\providecommand{\formul}[2]{
\begin{eqnarray}{#2}
#1
\end{eqnarray}
}

\title{Introducción al Cálculo Fraccional}
\author{
\begin{tabular}{cc}
A. Torres-Hernandez&
F. Brambila-Paz\\
 \href{mailto:anthony.torres@ciencias.unam.mx}{anthony.torres@ciencias.unam.mx}&
\href{mailto:fernandobrambila@gmail.com}{fernandobrambila@gmail.com}
\end{tabular}
}
\date{14-10-2017}

\begin{document}

El siguiente material fue creado con la idea de ser usado para un curso introductorio de cálculo fraccional.

\begin{titlepage}
\begin{center}
\begin{Huge}
\textsc{Introducción al Cálculo Fraccional}
\end{Huge}
\end{center}
\vspace{4cm}

\begin{figure}[!ht]
\centering
\subfigure{\includegraphics[height=5cm,width=8cm]{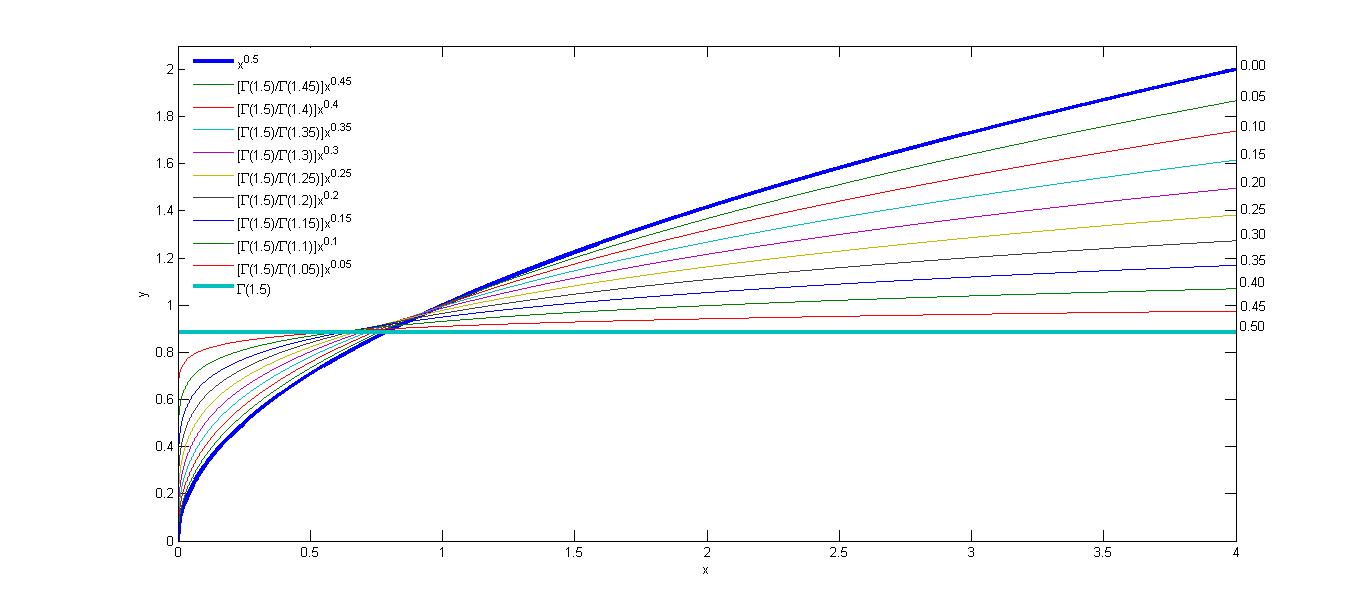}}
\subfigure{\includegraphics[height=5cm,width=8cm]{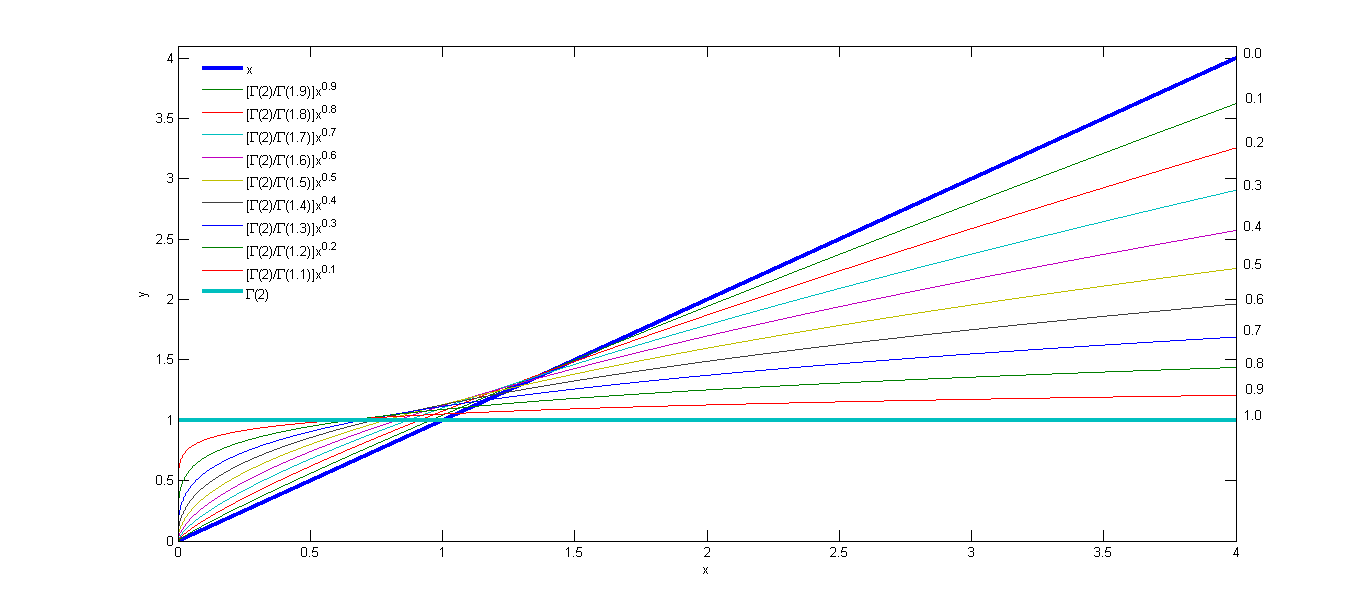}}
\subfigure{\includegraphics[height=5cm,width=8cm]{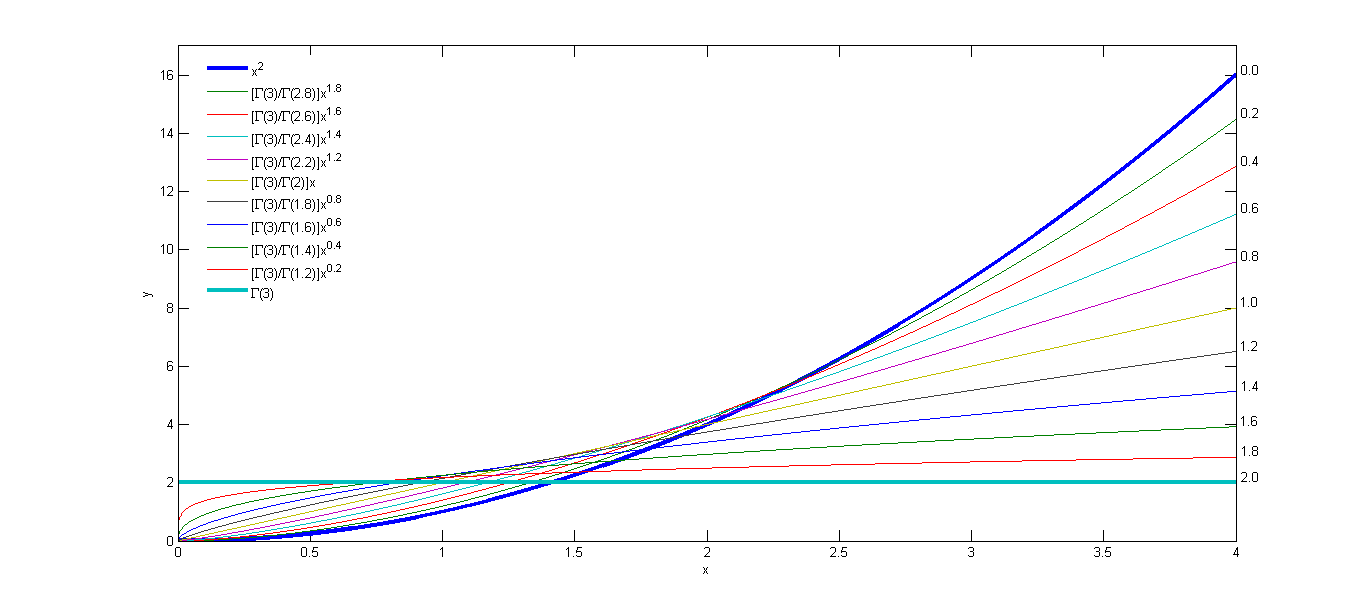}}
\subfigure{\includegraphics[height=5cm,width=8cm]{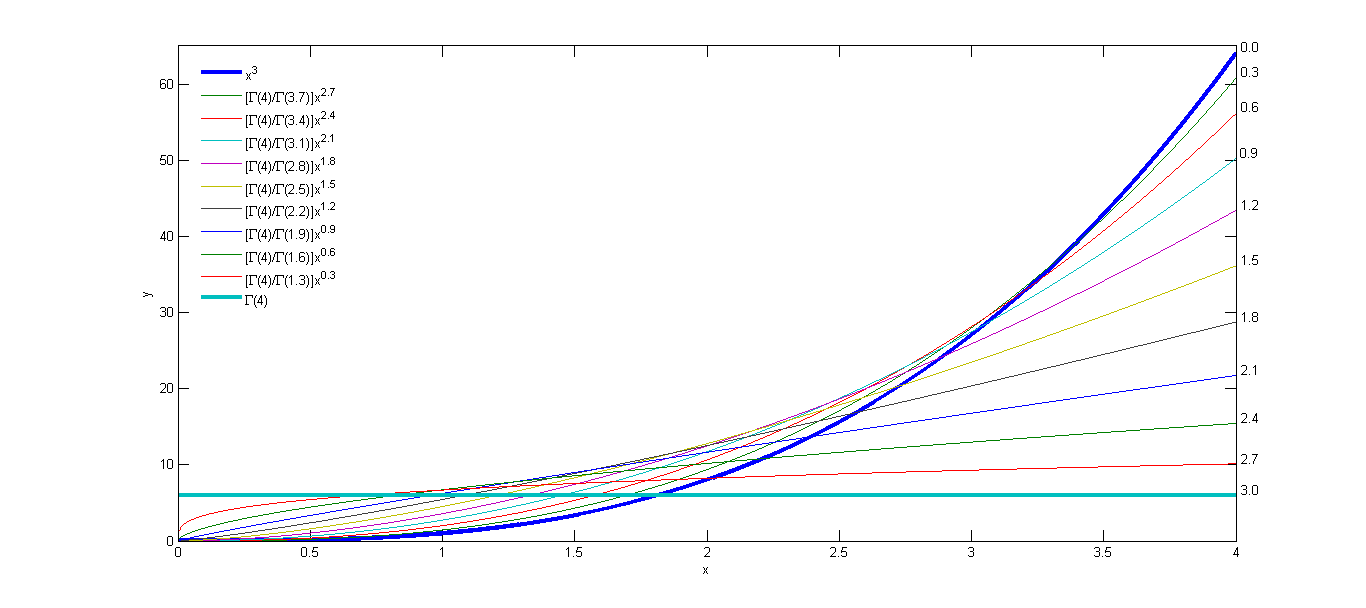}}
\subfigure{\includegraphics[height=5cm,width=8cm]{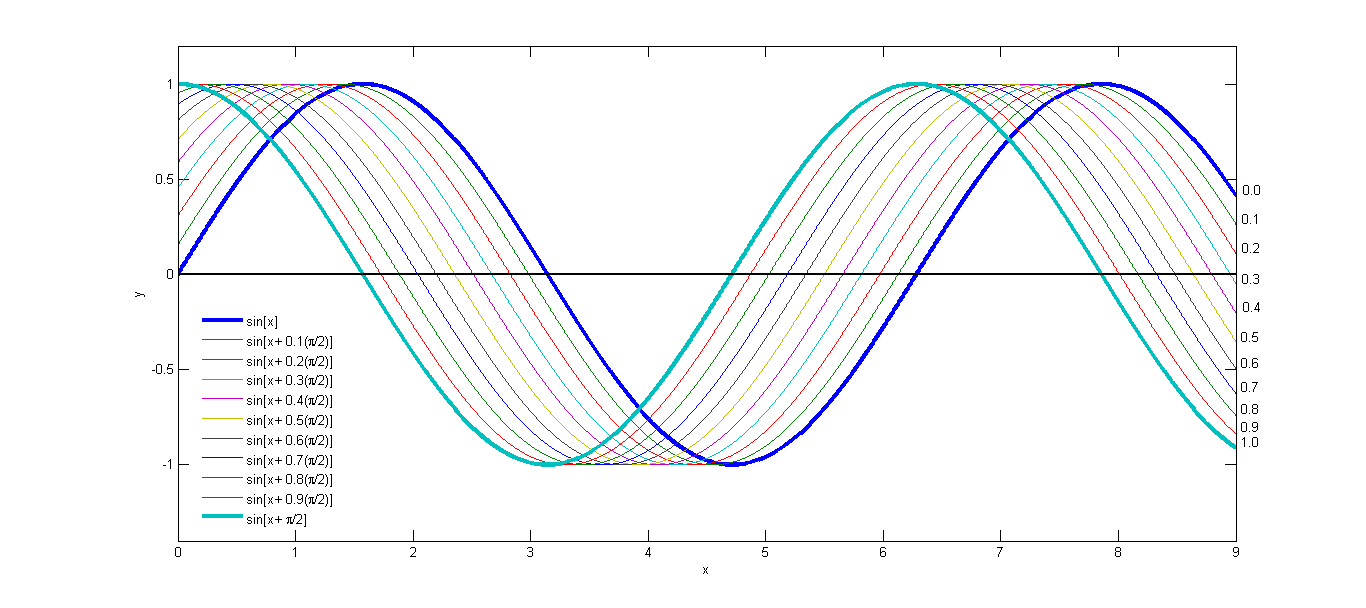}}
\subfigure{\includegraphics[height=5cm,width=8cm]{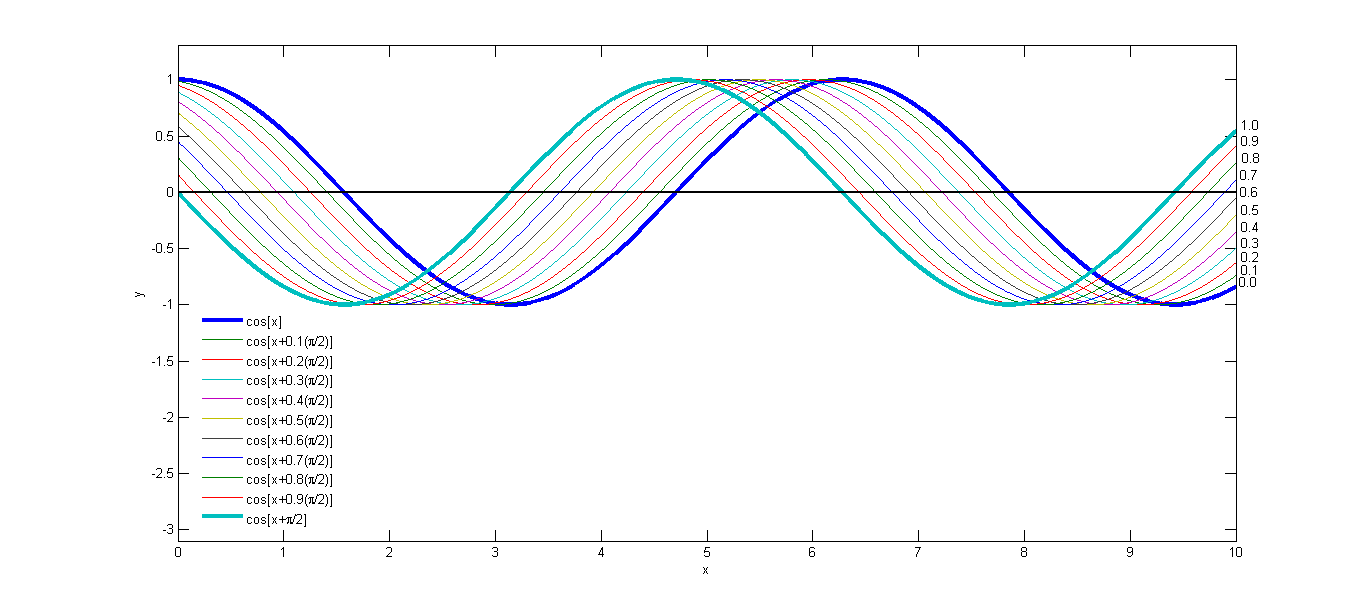}}
\end{figure}
\end{titlepage}

\addtocontents{toc}{\hfill \textbf{Página} \par}
\tableofcontents

\maketitle

\section{Cálculo Fraccional}

\subsection{Historia}

La pregunta que llevó al surgimiento de una nueva rama del análisis matemático conocida como cálculo fraccional fue: ¿Puede la definición de una derivada de orden entero \form{\left( \frac{d}{dx}\right)^n f(x)} ser extendida para el caso en que \form{n} sea una fracción y tener algún significado? Posteriormente la pregunta se volvió: ¿Puede  el orden de la derivada \form{n} ser cualquier número: racional, irracional o complejo? Debido a que la pregunta fue contestada de forma afirmativa, el nombre de cálculo fraccional se ha convertido en un nombre incorrecto y seria más apropiado llamarlo integración y diferenciación de orden arbitrario.

Leibniz inventó la notación \form{\left( \frac{d}{dx}\right)^nf(x)}. Quizás fue un juego por medio de símbolos lo que llevó a L'Hopital en \form{1695} a preguntarle a Leibniz: \com{¿Qué pasa si \form{n} es \form{\frac{1}{2}}?} Leibniz respondió \com{Usted puede ver con eso, señor, que uno puede expresar por una serie infinita una cantidad tal como \form{d^\frac{1}{2}\overline{xy}} o \form{d^{1:2}\overline{xy}}. Aunque la serie infinita y la geometríca son relaciones distantes, la serie infinita admite solamente el uso de exponentes que son enteros positivos y negativos, y aún no conoce el uso de exponentes fraccionarios.} \cite{miller93}

Posteriormente, en dicha carta, Leibniz continúa de forma profética: \com{Por lo tanto, \form{d^\frac{1}{2}x} será igual a \form{x\sqrt{dx:x}}. Esta es una aparente paradoja de la cual, un día, útiles consecuencias serán extraídas.}\cite{miller93}

El tema del cálculo fraccional no paso desapercibido de la atención de Euler. En \form{1730} él escribio \com{Cuando \form{n} es un número entero positivo, y  \form{p}  es una función de \form{x}, la razón \form{d^n p} a \form{dx^n} siempre se puede expresar algebraicamente, de modo que si \form{n=2} y \form{p=x^3}, entonces la razón \form{d^2x^3} a \form{dx^2} es \form{6x} a 1 .  Ahora me pregunto qué tipo de relación se puede hacer si \form{n} es una fracción. La dificultad en este caso puede ser fácilmente entendida. Si \form{n} es un entero positivo \form{d^n} se puede encontrar por diferenciación continua. Tal manera, sin embargo, no es evidente si \form{n} es una fracción. Pero aún con la ayuda de la interpolación que ya he explicado en esta disertación, uno puede ser capaz de agilizar el asunto.}\cite{miller93}

Lagrange contribuyó indirectamente al cálculo fraccional ya que en \form{1772} desarrolló la ley de exponentes para operadores diferenciales de orden entero:

\formula{
\der{d}{x}{m}\der{d}{x}{n}f(x) = \der{d}{x}{m+n}f(x).
}

Más tarde, cuando se desarrolló la teoría del cálculo fraccional, los matemáticos se interesaron en saber qué restricciones debían imponerse sobre funciones \form{f(x)} para que una regla análoga a la de ordenes enteros fuera válida para ordenes arbitrarios.

En \form{1812} Laplace definió una derivada fraccionaria por medio de una derivada  de orden entero, la primera mención de una derivada de orden arbitrario aparece en uno de los textos de Lacroix [\form{1819, pp. 409-410}], en el cual dedicó menos de dos páginas de las \form{700} en el texto a este tema. Desarrolló un único ejercicio matemático generalizando a partir de un caso de orden entero. Empezando con \form{f(x)=x^m} con \form{m} un entero positivo, Lacroix desarrolló fácilmente la derivada \form{n}-ésima de f(x) 

\formula{
\begin{array}{cccc}
\der{d}{x}{}f(x)&=& mx^{m-1},&\\
\der{d}{x}{2}f(x)&=& m(m-1)x^{m-2},&\\
\der{d}{x}{3}f(x)&=& m(m-1)(m-2)x^{m-3},&\\
&\vdots&\\
\der{d}{x}{n}f(x)&=&m(m-1)(m-2)\cdots(m-(n-1))x^{m-n}&\\
 &=&\dfrac{m!}{(m -n)!}x^{m-n},& \  m\geq n.
\end{array}
}

Haciendo uso del símbolo de Legendre para la generalización del factorial (la función Gamma), él obtiene

\begin{eqnarray}\label{eq:5}
\der{d}{x}{n}f(x)=\frac{\gam{m+1}}{\gam{m-n+1}}x^{m-n}.
\end{eqnarray}

Después, da el ejemplo para el caso \form{f(x)=x} con \form{n=\frac{1}{2}}, obteniendo

\formula{
\der{d}{x}{\frac{1}{2}}f(x)&=&\frac{\gam{2}}{\gam{\frac{3}{2}}}x^{\frac{1}{2}}\\
&=&\frac{2}{\sqrt{\pi}}x^\frac{1}{2}.
}

Cabe resaltar que el resultado obtenido por Lacroix es el mismo que se obtiene al utilizar la definición actual de una derivada fraccionaria de Riemann-Liouville. 

Fourier (\form{1822}) fue el siguiente en hacer mención sobre las derivadas de orden arbitrario. Su definición de un operador fraccionario se obtuvo de su representación de orden entero para \form{f(x):}

\formula{
f(x)=\frac{1}{2\pi}\int_{-\infty}^\infty f(\alpha)d\alpha \int_{-\infty}^\infty \co{p(x-\alpha)}dp.
}

ya que

\formula{
\der{d}{x}{n}\co{p(x-\alpha)}=p^n \co{p(x-\alpha)+n\frac{\pi}{2}},
}

para \form{n} un entero. Sustituyendo  \form{n} con \form{\nu} (\form{\nu} arbitrario), se obtiene la generalización

\formula{
\der{d}{x}{\nu}f(x)=\frac{1}{2\pi}\int_{-\infty}^\infty f(\alpha)d\alpha \int_{-\infty}^\infty p^\nu \co{p(x-\alpha)+\nu \frac{\pi}{2}}dp.
}

 El número \form{\nu} que aparece en la ecuación anterior se considerará como cualquier cantidad, ya sea positiva o negativa.

\subsubsection{Principales Aportaciones de Abel y Liouville al Cálculo Fraccional}

Leibniz, Euler, Laplace, Lacroix y Fourier hicieron mención sobre derivadas de orden arbitrario, pero el primer uso para operaciones fraccionarias fue hecho por Abel en \form{1823}, Abel aplicó el cálculo fraccional en la solución de una ecuación integral que surge en la formulación del problema de la tautócrona, i.e., el problema de determinar la forma de una curva de manera que el tiempo de descenso de una masa puntual sin fricción que se deslice por la curva bajo la acción de la gravedad sea independiente del punto de partida. Si el tiempo de deslizamiento es una constante conocida, entonces la ecuación integral de Abel es:

\formul{
k=\int_0^x (x-t)^{-\frac{1}{2}}f(t)dt.
}{\label{eq:1}}

La integral en la ecuación \eqref{eq:1} es, excepto por el factor multiplicativo \form{\gam{\frac{1}{2}}^{-1}}, un caso particular de una integral fraccionaria de orden \form{\frac{1}{2}} \cite{miller93}. Generalmente en ecuaciones integrales como en el caso anterior, la función \form{f} en el integrando es desconocida y debe ser determinada. Abel escribió el lado derecho de la ecuación como \form{\sqrt{\pi} \left( \frac{d}{dx}\right)^\frac{1}{2} f(x)}. Entonces él aplico a ambos lados de la ecuación  \form{\left( \frac{d}{dx}\right)^\frac{1}{2}} para obtener

\formul{
\der{d}{x}{\frac{1}{2}}k=\sqrt{\pi}f(x),
}{\label{eq:2}}

estos operadores fraccionarios (con condiciones adecuadas en las funciones \form{f}) tienen la propiedad de que \form{\left(\frac{d}{dx} \right)^\frac{1}{2} \left(\frac{d}{dx} \right)^{-\frac{1}{2}}f =\left(\frac{d}{dx} \right)^0f=f}. Así, cuando la derivada fraccionaria de orden \form{\frac{1}{2}} de la constante \form{k} en \eqref{eq:2} es obtenida, la función \form{f(x)} está determinada. Este es un logro notable de Abel en el cálculo fraccional. Es importante señalar que la derivada fraccionaria de una constante no siempre es igual a cero. 

El tema del cálculo fraccional permaneció inactivo durante un periodo de casi una década hasta que aparecieron las obras de Liouville.  Pero fue en \form{1974} que el primer texto [Oldham and Spanier]\cite{oldham74} dedicado exclusivamente a este tema fue publicado, durante estos años se impartió la primera conferencia relacionada con el cálculo fraccional.

Tal vez la fórmula integral de Fourier y la solución de Abel fueron unos de los principales trabajos que llamaron la atención de Liouville hacia el cálculo fraccional, quien realizo el primer estudio importante de este tema. Liouville tuvo éxito en la aplicación de sus definiciones a los problemas de la teoría potencial.

El punto de partida para su desarrollo teórico fue el resultado conocido para las derivados de orden entero de la función exponencial

\formula{
\der{d}{x}{m}\ex{ax}=a^m \ex{ax},
}

que extendió de una forma natural a derivadas de orden arbitrario

\formula{
\der{d}{x}{\nu} \ex{ax}=a^\nu \ex{ax}.
}

Él asumió que la derivada de orden arbitrario de una función \form{f(x)}, la cual se puede extender en una serie de la forma

\begin{eqnarray}\label{eq:3}
f(x)=\sum_{k=0}^\infty c_k \ex{a_k x}, \ \re{a_k}>0,
\end{eqnarray}

debería ser

\begin{eqnarray}\label{eq:4}
\der{d}{x}{\nu} f(x)= \sum_{k=0}^\infty c_k a_k^\nu \ex{a_k x}.
\end{eqnarray}

La ecuación \eqref{eq:4} es conocida como la primera definición de Liouville para una derivada fraccionaria. Esta formula generaliza de manera natural una derivada de orden arbitrario \form{\nu}, donde \form{\nu} es cualquier número: racional, irracional o complejo. Pero tiene la obvia desventaja de ser aplicable sólo a las funciones que pueden tomar la forma \eqref{eq:3}. Tal vez Liouville era consciente de esta desventaja, pues formuló una segunda definición.

Para obtener su segunda definición comenzó con una integral definida relacionada con la función Gamma:

\formula{
I=\int_0^\infty u^{a-1}\ex{-xu}du, \ a,x>0.
}

Tomando el cambio de variable \form{t=xu } se obtiene

\formula{
I&=& x^{-a}\int_0^\infty t^{a-1}\ex{-t}dt\\
&=& x^{-a}\gam{a},
}

o

\formula{
x^{-a}=\frac{1}{\gam{a}}I.
}

Después él aplico \form{\left(\frac{d}{dx} \right)^\nu} a ambos lados de la ecuación anterior, para obtener, según el supuesto básico de Liouville

\formula{
\der{d}{x}{\nu} x^{-a}=\frac{(-1)^\nu}{\gam{a}}\int_0^\infty u^{a+\nu-1}\ex{-xu}du.
}

De esta forma Liouville obtiene su segunda definición de una derivada fraccionaria:

\formula{
\der{d}{x}{\nu} x^{-a}=\frac{(-1)^\nu \gam{a+\nu}}{\gam{a}}x^{-a-\nu}, \ a>0.
}

Pero estas definiciones de Liouville eran demasiado restrictivas para permanecer un largo tiempo. La primera definición se restringe solo a funciones de la forma \eqref{eq:3}, y la segunda definición es útil sólo para la función del tipo \form{x^{-a}} (con \form{a>0}). Ninguno de los dos tipos es adecuado para ser aplicado a una clase más amplia de funciones.

\subsubsection{Aportación de Riemann al Cálculo Fraccional}

Riemann desarrolló principalmente su teoría sobre la integración fraccionaria durante su época como estudiante, pero rechazó la publicación de sus trabajos. Estos fueron publicados póstumamente en \form{1892}. Él buscó una generalización de una serie de Taylor y sus derivadas

\begin{eqnarray}\label{eq:6}
D^{-\nu}f(x)=\frac{1}{\gam{\nu}}\int_c^x(x-t)^{\nu-1}f(t)dt + \Psi(x).
\end{eqnarray}

Debido a la ambigüedad en el límite inferior \form{c} de la integral, Riemann consideró conveniente añadir una función complementaria \form{\Psi(x)}. Esta función complementaria es en esencia un intento de proporcionar una medida de la desviación de la ley de los exponentes

\formula{
\ifr{c}{D}{x}{-\mu}\ifr{c}{D}{x}{-\nu}f(x)=\ifr{c}{D}{x}{-\mu-\nu}f(x),
}

(donde los subíndices \form{c} y \form{x} en \form{D} se refieren a los límites de integración en \eqref{eq:6}) la cual es válida cuando los limites inferiores \form{c} son iguales. Riemann se preocupó por una medida de desviación para el caso \form{\ifr{c}{D}{x}{-\mu}\ifr{c'}{D}{x}{-\nu}f(x)} cuando \form{c\neq c'}.

A. Cayley (\form{1880}) comentó: \com{La mayor dificultad de la teoría de Riemann, me parece, es la cuestión del significado de una función complementaria que contiene una infinidad de constantes arbitrarias \cite{miller93}.} Cualquier definición satisfactoria de una operación fraccionaria exigirá que se elimine esta dificultad. De hecho, la definición actual de integración fraccionaria se presenta sin la función complementaria.

\subsubsection{Surgimiento de la Definición de Riemann-Liouville para Derivadas Fraccionarias}

El primer trabajo que llevó a lo que ahora se conoce como la definición de Riemann-Liouville parece ser el trabajo de Sonin (\form{1869}) \com{Sobre diferenciación con índice arbitrario.} Su punto de partida fue la fórmula integral de Cauchy. Letnikov escribió cuatro artículos sobre este tema de \form{1868} a \form{1872}. Su trabajo \com{Una explicación de los principales conceptos de la teoría de la diferenciación de índice arbitrario} (\form{1872}) es una extensión del trabajo de Sonin \cite{miller93}. La \form{n}-ésima derivada de la fórmula integral de Cauchy viene dada por

\begin{eqnarray}\label{eq:7}
D^nf(z)=\frac{n!}{2\pi i}\int_C \frac{f(\zeta)}{(\zeta-z)^{n+1}}d\zeta.
\end{eqnarray}

No presenta ningún problema la generalización de \form{n!} a valores arbitrarios debido a que \form{\nu!=\gam{\nu +1}}. Sin embargo, cuando \form{n} no es un entero, el integrando en \eqref{eq:7} ya no presenta un polo, sino un punto de ramificación. El contorno apropiado requeriría entonces un corte de rama que no se incluyó en el trabajo de Sonin y Letnikov, aunque se discutió.

No fue hasta que Laurent (\form{1884}) publicó su trabajo, que la teoría de operadores generalizados alcanzó un nivel adecuado como punto de partida para el matemático moderno. La teoría del cálculo fraccional está íntimamente conectada con la teoría de los operadores. El operador \form{D=\frac{d}{dx}} y \form{D^2=\left(\frac{d}{dx} \right)^2} denotan una regla de transformación que es familiar para cualquiera que haya estudiado el cálculo convencional. El punto de partida de Laurent fue también la fórmula integral de Cauchy. Su contorno era un circuito abierto en una superficie de Riemann, en contraste con el circuito cerrado de Sonin y Letnikov. El método de integración de contornos produjo la definición

\begin{eqnarray}\label{eq:8}
\ifr{c}{D}{x}{-\nu}f(x)=\frac{1}{\gam{\nu}}\int_c^x (x-t)^{\nu-1}f(t)dt, \ \re{\nu}>0,
\end{eqnarray}

para la integración de orden arbitrario.

Cuando \form{x>c} en la ecuación \eqref{eq:8}, tenemos la definición de Riemann pero sin una función complementaria. La versión más usual se presenta cuando \form{c=0},

\formul{
\ifr{0}{D}{x}{-\nu}f(x)&=&\ifrR{RL}{0}{D}{x}{-\nu}f(x)\nonumber\\
&=&\frac{1}{\gam{\nu}}\int_0^x (x-t)^{\nu-1}f(t)dt, \ \re{\nu}>0.
}{\label{eq:9}}

Esta forma de la integral fraccionaria a menudo se conoce como la integral fraccionaria de Riemann-Liouville. Una condición suficiente para que la ecuación \eqref{eq:9} converja es  \cite{miller93} 

\formula{
f\left(x^{-1}\right)\sim \mathcal{O}(x^{1-\epsilon}), \ \epsilon >0.
}

Las funciones integrables con la propiedad anterior se conocen a veces como funciones de la clase de Riemann. Por ejemplo, las constantes son de la clase de Riemann, al igual que

\formula{
x^a, \ a>-1.
}

Cuando \form{c} es el infinito negativo, la ecuación \eqref{eq:8} se convierte en

\formul{
\ifr{-\infty}{D}{x}{-\nu}f(x)=\frac{1}{\gam{\nu}}\int_{-\infty}^x (x-t)^{\nu-1}f(t)dt, \ \re{\nu}>0.
}{\label{eq:10}} 

Una condición suficiente para que \eqref{eq:10} converja es \cite{miller93}

\formula{
f(-x)\sim \mathcal{O}(x^{-\nu-\epsilon}), \ \epsilon >0, \  x\rightarrow \infty.
}

Las funciones integrables con la propiedad anterior se conocen a veces como funciones de la clase de Liouville. Por ejemplo

\formula{
x^{-a}, \ a>\nu>0,
}

es de clase Liouville. Una constante no lo es. Sin embargo, si \form{a} está entre \form{-1} y \form{0}, dependiendo del valor de \form{\nu}, las dos clases pueden coincidir. Tomando \form{f(t)=\ex{at}, \re{a}>0}, en la ecuación \eqref{eq:10}, se obtiene

\formul{
\ifr{-\infty}{D}{x}{-\nu}\ex{ax}=a^{-\nu}\ex{ax}.
}{\label{eq:11}}

Asumiendo que la ley de exponentes \form{D[D^{-\nu}f(x)]=D^{1-\nu}f(x)} se cumple, entonces si \form{0<\nu<1}, se obtiene \form{\mu=1-\nu >0} y la ecuación \eqref{eq:11} se convierte en

\formula{
\ifr{-\infty}{D}{x}{\mu}\ex{ax}=a^\mu \ex{ax}, \ \re{a}>0.
}

Cabe resaltar que la primera definición de Liouville se encuentra incluida en la ecuación \eqref{eq:10}.

Sin embargo, si \form{f(x)=x^{-a}, a>\nu >0}, entonces

\formul{
\ifr{-\infty}{D}{x}{-\nu}x^{-a}=(-1)^\nu \frac{\gam{a-\nu}}{\gam{a}}x^{-a+\nu},
}{\label{eq:12}}

para \form{x<0}, y si \form{0<\nu<1}, entonces \form{\mu=1-\nu >0} y

\formula{
\ifr{-\infty}{D}{x}{\mu}x^{-a}=(-1)^\mu \frac{\gam{a+\mu}}{\gam{a}}x^{-a-\mu}.
}

Este es el mismo resultado que se obtiene para la segunda definición de Liouville, excepto que él asumió que \form{x>0}. Si \form{x>0}, \eqref{eq:12}, es cierto sólo para el rango \form{0<\nu<a<1}.

Para \form{f(x)=x^a} y \form{\nu>0}, tenemos de \eqref{eq:9} que

\formula{
\ifrR{RL}{0}{D}{x}{-\nu}x^a = \frac{\gam{a+1}}{\gam{a+\nu+1}}x^{a+\nu}, \ a>-1,
}

y asumiendo nuevamente que \form{D[D^{-\nu}f(x)]=D^{1-\nu}f(x)}, se observa que si \form{0<\nu <1},

\formul{
\ifrR{RL}{0}{D}{x}{\nu}x^a = \frac{\gam{a+1}}{\gam{a-\nu+1}}x^{a-\nu}, \ a>-1.
}{\label{eq:13}}

Cabe señalar que para \form{f(x)=x} con \form{\nu=\frac{1}{2}}, la ecuación \eqref{eq:13} produce el mismo resultado obtenido por Lacroix. También se puede considerar la observación de Center \cite{miller93} sobre la derivada de orden arbitrario de una constante. Pero si \form{f(x)=1} con \form{\nu=\frac{1}{2}}, entonces \eqref{eq:13} da como resultado

\formula{
\ifrR{RL}{0}{D}{x}{\frac{1}{2}}(1)=\frac{1}{\sqrt{\pi}}x^{-\frac{1}{2}}.
}

Pero Center estaba equivocado cuando dijo que la definición de Liouville da cero para la derivada arbitraria de una constante. Porque usaba  \form{D^\nu x^{-a}=\frac{(-1)^\nu \gam{a+\nu}}{\gam{a}}x^{-a-\nu}, \ a>0}. Pero \form{1=x^0=(-x)^0} no está en la clase de Liouville.

\subsubsection{Finales del Siglo XIX e Inicios del Siglo XX}

Heaviside (\form{1892}) publicó una serie de artículos en los que mostró cómo ciertas ecuaciones diferenciales lineales pueden resolverse mediante el uso de operadores generalizados. Sus métodos han demostrado ser útiles en la teoría sobre la transmisión de corrientes eléctricas en cables, y han sido adoptados bajo el nombre de cálculo operacional de Heaviside.

El cálculo operacional de Heaviside se refiere a operadores funcionales lineales. Él denotaba al operador de diferenciación por la letra \form{p} y hacia uso de el como si se tratase de una constante en la solución de ecuaciones diferenciales. Por ejemplo, para la ecuación de calor en una dimensión

\formul{
\der{d}{x}{2}u=a^2 \frac{d}{dt}u,
}{\label{eq:15}}

donde \form{a^2} es una constante y \form{u} es la temperatura. Tomando

\formula{
\frac{d}{dt}=p,
}

entonces la ecuación \eqref{eq:15} toma la forma

\formul{
D^2 u =a^2 p u.
}{\label{eq:16}}

Gregory (\form{1841}), dijo ser el fundador del llamado cálculo de operaciones, había puesto la solución de de la ecuación \eqref{eq:15} en forma de un operador:

\formula{
u(x,t)=A\ex{axp^\frac{1}{2}}+ B\ex{-axp^\frac{1}{2}}.
}

Esto es lo que  se obtendría si se resolviera la ecuación \eqref{eq:16} asumiendo que \form{p} es una constante.

Pero fueron las brillantes aplicaciones de Heaviside las que aceleraron el desarrollo de la teoría de estos operadores generalizados. Él obtuvo resultados correctos al expandir la exponencial en serie de potencias de \form{p^\frac{1}{2}}, donde \form{p^\frac{1}{2}=\left( \frac{d}{dx}\right)^\frac{1}{2}=D^\frac{1}{2}}. En la teoría de circuitos eléctricos, Heaviside encontró un uso frecuente para los operadores \form{p^\frac{1}{2}}. Él interpretó que \form{p^\frac{1}{2}\rightarrow 1}, es decir  \form{D^\frac{1}{2}(1)}, obteniendo \form{(\pi t)^{-\frac{1}{2}}}. Ya que \form{f(t)=1} es una función de la clase de Riemann, está claro que el operador de Heaviside debe interpretarse en el contexto del operador de Riemann \form{\ifr{0}{D}{x}{\nu}}.

Su resultado era correcto, pero no pudo justificar sus procedimientos, comentó Kelland, en el intervalo de diez años entre la publicación de Fourier y las aplicaciones de Liouville  una situación similar siguió a las publicaciones de Heaviside, salvo que en este caso, un tiempo mucho más largo transcurrió antes de que sus procedimientos estuvieran justificados por Bromwich (\form{1919}).

Harold T. Davis (\form{1939}) dijo: \com{El período de desarrollo formal de los métodos operativos puede considerarse como terminado por \form{1900}. La teoría de las ecuaciones integrales estaba comenzando a despertar la imaginación de los matemáticos y a revelar las posibilidades de los métodos operacionales.}\cite{miller93}

En el periodo de \form{1900} a \form{1970} una modesta cantidad de trabajos publicados aparecieron sobre el tema del cálculo fraccional. Algunos de los que contribuyeron fueron Al-Bassam, Davis,  Erdélyi, Hardy,  Kober, Littlewood,  Love,  Osler,  Riesz,  Samko,  Sneddon,  Weyl y Zygmund.

En el año de \form{1974} se dio la primera conferencia internacional sobre el cálculo fraccional, celebrada en la Universidad de New Haven, Connecticut, y fue patrocinado por la Fundación Nacional de Ciencias. Se presento una cantidad moderada de trabajos relacionadas con el calculo fraccionario, incluyendo desigualdades obtenidas mediante el uso del cálculo fraccional, trabajos sobre el cálculo fraccional y las funciones generalizadas y aplicaciones del cálculo fraccional a la teoría de probabilidades.

En \form{1984} la segunda conferencia internacional sobre cálculo fraccional fue patrocinada por la Universidad de Strathclyde, Glasgow, Scotland. Entre los que contribuyeron se incluyen P. Heywood, S. Kalla, W. Lamb, J. S. Lowndes, K. Nishimoto, P. G. Rooney, y H. M. Strivastada, así como algunos de los matemáticos que participaron en la conferencia de New Haven. 

Una considerable actividad matemática relacionada con el cálculo fraccional se desarrollo en Japon en los años \form{80} con publicaciones de S. Owa (\form{1990}), M. Saigo (\form{1980}), y K. Nishimoto. El último autor publicó un trabajo de cuatro volúmenes (\form{1984,1987, 1989, 1991}) dedicados principalmente a las aplicaciones del cálculo fraccional a ecuaciones diferenciales ordinarias y parciales. En la Unión Soviética tres matemáticos, S. Samko, O. Marichev, y A. Kilbas, escribieron un texto enciclopédico sobre el cálculo fraccional y algunas de sus aplicaciones (\form{1987}).

La tercera conferencia internacional se celebró en la Universidad de Nihon  en Tokyo en \form{1989}. Algunas de los muchos que  contribuyeron fueron M. Al-Bassam, R. Bagley, Y. A. Brychkov, L. M. B. C. Campos, R. Gorenflo, J. M. C. Joshi, S. Kalla, E. R. Lobe, M. Mikolás, K. Nishimoto, S. Owa, A. P. Prudnikov, B. Ross, S. Samko, H. M. Srivastada.

El cálculo fraccional se utiliza en muchos campos de la ciencia y la ingeniería, incluyendo el flujo de fluidos, reología, transporte difusivo similar a la difusión, las redes eléctricas, la teoría electromagnética y la probabilidad. Algunos artículos de P. C. Phillips (\form{1989,1990}) han utilizado el cálculo fraccional en la estadística. R. L. Bagley (\form{1990}); Bagley y Torkiv (\form{1986}) han encontrado uso para el cálculo fraccional en los temas de viscoelasticidad y la electroquímica de corrosión.

\newpage

\subsection{Propiedades Básicas}

\subsubsection{Propiedades de la Función Gamma}

La función Gamma \form{\gam{z}} juega un papel importante en la teoría de la diferenciación. En consecuencia, es conveniente recopilar ciertas fórmulas y propiedades relacionadas con esta función. Una definición de \form{\gam{z}} con \form{z\in \mathds{C}}, es proporcionada por el límite infinito de Euler \cite{oldham74,arfken85}

\formul{
\gam{z}=\lim_{n\to \infty}\frac{n!}{z(z+1)(z+2)\cdots (z+n)}n^z,
}{\label{eq:19}}

\begin{figure}[!ht]
\centering
\includegraphics[height=0.35\textheight, width=\textwidth]{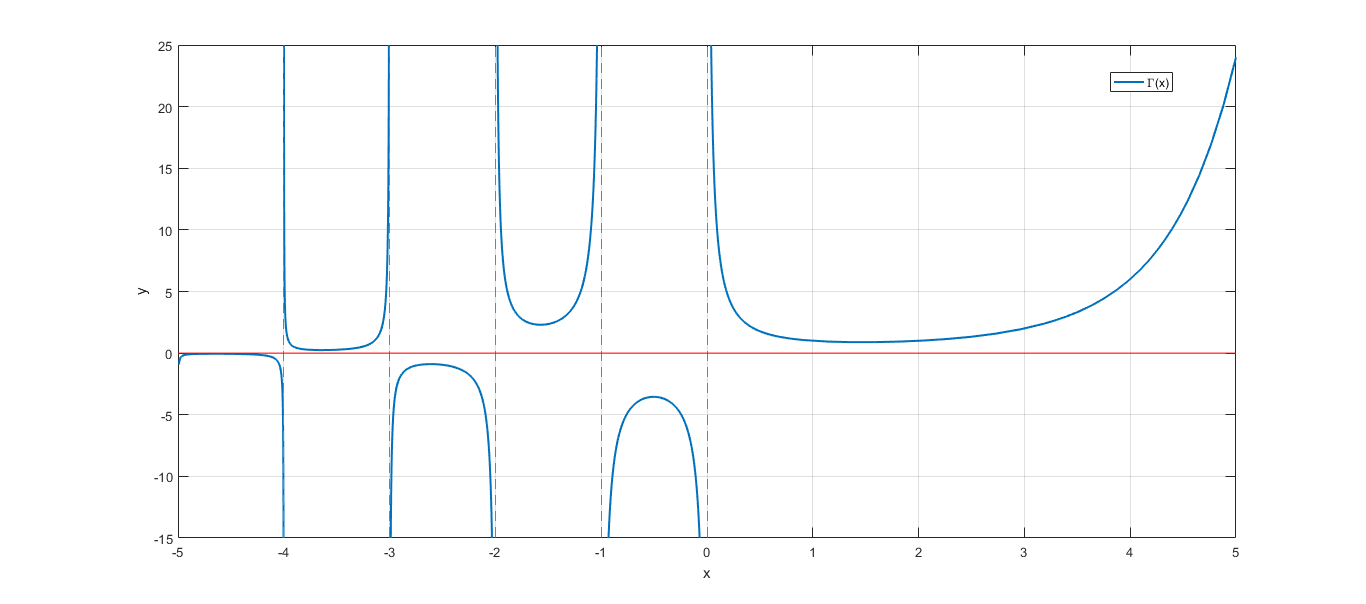}
\caption{Grafica de la función Gamma \form{\gam{x}} en el intervalo \form{[-5,5]}.}
\end{figure}

la definición anterior también puede ser escrita como

\formula{
\gam{z}=\lim_{n\to \infty}\frac{n!(z-1)!}{(z+n)!}n^z,
}

haciendo \form{z\to z+1} en la ecuación anterior

\formula{
\gam{z+1}&=&\lim_{n\to \infty}\frac{n!\cdot z!}{(z+n+1)!}n^{z+1}\\
&=&\lim_{n\to\infty}\frac{zn}{z+n+1}\cdot \frac{n!(z-1)!}{(z+n)!}n^z\\
&=&\lim_{n\to\infty}\frac{zn}{z+n+1}\cdot \lim_{n\to\infty} \frac{n!(z-1)!}{(z+n)!}n^z\\
&=&z \lim_{n\to\infty}\frac{n}{n+(z+1)}\cdot \lim_{n\to\infty} \frac{n!(z-1)!}{(z+n)!}n^z\\
&=& z \lim_{n\to\infty} \frac{n!(z-1)!}{(z+n)!}n^z\\
&=& z\gam{z},
}

se demuestra una de las propiedades fundamentales de la función Gamma

\formula{
\gam{z+1}=z\gam{z},
}

otra propiedad importante es la siguiente

\formula{
\gam{1}&=&\lim_{n\to \infty}\frac{n!}{(n+1)!}n\\
&=& \lim_{n\to \infty} \frac{n}{n+1}\\
&=& 1,
}

de las dos propiedades anteriores se obtiene 

\formula{
\gam{z+1}&=&z(z-1)\cdots  2\cdot 1\\
&=&z!.
}

Pero la definición integral, conocida como integral definida de Euler, es a menudo más útil, aunque se restringe a valores de \form{z\in \mathds{C}} con \form{\re{z}>0}

\formul{
\gam{z}=\int_{0}^\infty \ex{-t}t^{z-1}dt, \ \re{z}>0
}{\label{eq:20}}

esta expresión de la función Gamma puede ser reescrita en dos formas que suelen ser útiles, tomando el cambio variable \form{t=u^2} \cite{arfken85} en la ecuación anterior

\formula{
\gam{z}&=&\int_0^\infty \ex{-t}t^{z-1}dt\\
&=& \int_0^\infty \ex{-u^2}(u^2)^{z-1}(2udu)\\
&=&2\int_0^\infty \ex{-u^2}u^{2z-2}udu\\
&=&2\int_0^\infty \ex{-t^2}t^{2z-1}dt,\ \re{z}>0,
}

si ahora se toma el cambio de variable \form{t=-\ln u} \cite{arfken85}

\formula{
\gam{z}&=&\int_0^\infty e^{-t}t^{z-1}dt\\
&=&\int_1^{0} u\left( -\ln u\right) ^{z-1}\left( -\frac{1}{u}du\right)\\
&=&-\int_1^0 \left(\ln \frac{1}{u}\right)^{z-1} du\\
&=& \int_0^1 \left(\ln \frac{1}{t}\right)^{z-1} dt, \ \re{z}>0 ,
}

estas formas equivalentes son de ayudan en ocasiones, por ejemplo tomando la primera es fácil obtener el valor de \form{\gam{\frac{1}{2}}}

\formula{
\gam{\frac{1}{2}} &=& 2\int_0^\infty \ex{-t^2}t^{2\cdot \frac{1}{2}-1}dt\\
&=&2\int_0^\infty \ex{-t^2}dt\\
&=& 2\left(\frac{\sqrt{\pi}}{2} \right)\\
&=& \sqrt{\pi},
}

para mostrar la equivalencia de las ecuaciones \eqref{eq:19} y \eqref{eq:20} se considera le ecuación de dos variables \cite{arfken85}

\formul{
F(z,n)=\int_0^n \left(1-\frac{t}{n} \right)^n t^{z-1}dt, \ \re{z}>0,
}{\label{eq:30}}

tomando el limite cuando \form{n\to \infty} en la ecuación \eqref{eq:30}

\formula{
\lim_{n\to \infty} F(z,n)&=&\lim_{n\to \infty} \int_0^n \left(1-\frac{t}{n} \right)^n t^{z-1}dt\\
&=& \int_0^n \lim_{n\to \infty} \left(1-\frac{t}{n} \right)^n t^{z-1}dt\\
&=&\int_0^\infty \ex{-t}t^{z-1}dt\\
&=& \gam{z},
}

por otra parte tomando el cambio de variable \form{t=nu}

\formula{
F(z,n)&=&\int_0^n \left(1-\frac{t}{n} \right)^n t^{z-1}dt\\
&=&\int_0^1 (1-u)^n (nu)^{z-1}ndu\\
&=&n^z\int_0^1(1-u)^nu^{z-1}du,
}

integrando por partes se obtiene 

\formula{
F(z,n)&=&n^z \left( \left. (1-u)^n \frac{u^z}{z}\right|_0^1+\frac{n}{z}\int_0^1 (1-u)^{n-1}u^zdu \right)\\
&=&n^z \left(\frac{n}{z}\int_0^1 (1-u)^{n-1}u^zdu\right),
}

entonces se puede deducir que

\formula{
F(z,n)&=& n^z\left(\frac{n(n-1)\cdots 2\cdot 1}{z(z+1)\cdots (z+n-1)}\int_0^1 u^{z+n-1}du\right)\\
&=&n^z \left(\frac{n(n-1)\cdots 2\cdot 1}{z(z+1)\cdots (z+n-1)}\left. \frac{u^{z+n}}{z+n}\right|_0^1\right)\\
&=& \left(\frac{n(n-1)\cdots 2\cdot 1}{z(z+1)\cdots (z+n-1)(z+n)}\right)n^z\\
&=&\frac{n!(z-1)!}{(z+n)!}n^z,
}

tomando el limite cuando \form{n\to \infty}

\formula{
\lim_{n\to \infty} F(z,n)&=&\lim_{n\to \infty}\frac{n!(z-1)!}{(z+n)!}n^z
\\
&=&\gam{z}.
}

Debido a la relación de recurrencia \form{\gam{z+1}=z\gam{z}} se puede mostrar que para un entero positivo \form{n}

\formula{
\gam{n+1}&=& n\gam{n}\\
&=& n(n-1)\gam{n-1}\\
&\vdots &\\
&=& n(n-1)\cdots 2\cdot 1 \gam{1}\\
&=& n!,
}

reescribiendo la relación de recurrencia como

\formula{
\gam{z-1}=\frac{\gam{z}}{(z-1)},
}

la relación de recurrencia también ayuda a extender la definición de la función Gamma a valores negativos para los cuales la definición \eqref{eq:20} no es valida. Esta relación muestra que \form{\gam{0}} es infinito, así como \form{\gam{-1}} y cualquier valor de la función Gamma en los enteros negativos. Sin embargo los cocientes entre funciones Gamma para enteros negativos son finitos \cite{oldham74}, entonces si \form{N} y \form{n} son enteros positivos

\formula{
\frac{\gam{-n}}{\gam{-N}}&=& \frac{\gam{1-n}}{\gam{1-N}}\frac{(-N)}{(-n)}\\
&=& \frac{\gam{2-n}}{\gam{2-N}}\frac{(1-N)(-N)}{(1-n)(-n)}\\
&=& \frac{\gam{3-n}}{\gam{3-N}}
\frac{(2-N)(1-N)(-N)}{(2-n)(1-n)(-n)}\\
&\vdots&\\
&=&\frac{\gam{-1}}{\gam{-1}}\frac{(-2)(-3)(-4)\cdots (2-N)(1-N)(-N)}{(-2)(-3)(-4)\cdots (2-n)(1-n)(-n)}\\
&=&(-1)^{N-n}\frac{N!}{n!}\\
&=& (-1)^{N-n}\frac{\gam{N+1}}{\gam{n+1}}.
}

El recíproco de la función Gamma es univaluda y finita para todo \form{z\in \mathds{C}}, es conocido como el producto infinito de Weierstrass \cite{arfken85}

\formula{
\frac{1}{\gam{z}}=z\ex{\gamma z}\prod_{n=1}^\infty \left(1+\frac{z}{n} \right) \ex{-\frac{z}{n}},
}

\begin{figure}[!ht]
\centering
\includegraphics[height=0.35\textheight, width=\textwidth]{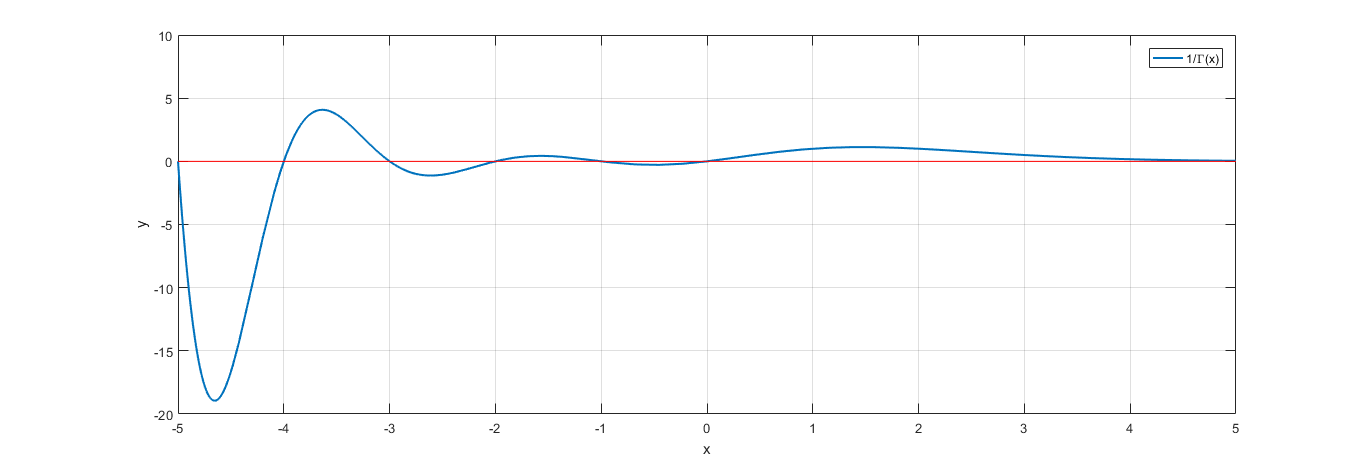}
\caption{Grafica de la función \form{\frac{1}{\gam{x}}} en el intervalo \form{[-5,5].}}\label{fig:1}
\end{figure}

donde $\gamma$ es la constante de Euler-Mascheroni, esta forma se puede deducir a partir del producto infinito de Euler, ya que

\formula{
\gam{z}&=&\lim_{n\to \infty}\frac{n!}{z(z+1)(z+2)\cdots (z+n)}n^z\\
&=& \lim_{n\to \infty } \cfrac{1}{z}\left[  \cfrac{(z+1)(z+2)\cdots (z+n)}{n!} \right]^{-1} n^z\\
&=& \lim_{n\to \infty } \cfrac{1}{z}\left[ \left(\frac{1+z}{1} \right) \left(\frac{2+z}{2} \right)\cdots \left(\frac{n+z}{n} \right)\right]^{-1} n^z\\
&=&\lim_{n \to \infty} \frac{1}{z}\prod_{m=1}^n \left( 1+\frac{z}{m}\right)^{-1}n^z,
}

invirtiendo la expresión anterior se obtiene

\formul{
\frac{1}{\gam{z}}=\lim_{n\to \infty} z\prod_{m=1}^n\left(1+\frac{z}{m} \right)n^{-z},
}{\label{eq:31}}

multiplicando  por 

\formula{
1= \prod_{m=1}^n \ex{-\frac{z}{m}}\ex{z\sum_{m=1}^n \frac{1}{m}},
}
 
la ecuación \eqref{eq:31}  se puede reescribir como

\formula{
\frac{1}{\gam{z}}&=&\lim_{n\to \infty} z\prod_{m=1}^n\left(1+\frac{z}{m} \right)n^{-z}\\
&=&\lim_{n\to \infty} z\prod_{m=1}^n\left(1+\frac{z}{m} \right)n^{-z}\prod_{m=1}^n \ex{-\frac{z}{m}}\ex{z\sum_{m=1}^n \frac{1}{m}}\\
&=&\lim_{n\to \infty} z \ex{z\sum_{m=1}^n \frac{1}{m}}  n^{-z} \prod_{m=1}^n\left(1+\frac{z}{m} \right)\prod_{m=1}^n \ex{-\frac{z}{m}}\\
&=&z \lim_{n\to \infty}   \ex{z\sum_{m=1}^n \frac{1}{m}} \ex{-z\ln\left( n\right)} \prod_{m=1}^n\left(1+\frac{z}{m} \right) \ex{-\frac{z}{m}}\\
&=&z \lim_{n\to \infty}   \ex{z\sum_{m=1}^n \frac{1}{m} -z\ln\left( n\right) } \lim_{n\to \infty} \prod_{m=1}^n\left(1+\frac{z}{m} \right) \ex{-\frac{z}{m}},
}

como la función exponencial es continua

\formula{
\lim_{n\to \infty}   \ex{z\sum_{m=1}^n \frac{1}{m} -z\ln\left( n\right) }
&=& \lim_{n\to \infty}   \ex{z\left(\sum_{m=1}^n \frac{1}{m} -\ln\left( n\right)\right) }\\
&=& \ex{z \lim_{n\to \infty} \left(\sum_{m=1}^n \frac{1}{m} -\ln\left( n\right)\right) },
}

de donde se obtiene la constante de Euler-Mascheroni

\formula{
\gamma &=&  \lim_{n\to \infty} \left(\sum_{m=1}^n \frac{1}{m} -\ln\left( n\right)\right)\\
&\simeq &0.577216\cdots,
}

entonces

\formula{
\frac{1}{\gam{z}}&=& z \lim_{n\to \infty}   \ex{z\sum_{m=1}^n \frac{1}{m} -z\ln\left( n\right) } \lim_{n\to \infty} \prod_{m=1}^n\left(1+\frac{z}{m} \right) \ex{-\frac{z}{m}}\\
&=& z \ex{z \gamma } \lim_{n\to \infty} \prod_{m=1}^n\left(1+\frac{z}{m} \right) \ex{-\frac{z}{m}}\\
&=& z \ex{z \gamma } \prod_{m=1}^\infty \left(1+\frac{z}{m} \right) \ex{-\frac{z}{m}}.
}

La figura \eqref{fig:1} muestra una gráfica de esta función para valores de \form{x\in \mathds{R}}. Se debe notar que se alternan los signos para  argumentos negativos de la función y que se presenta un acercamiento asintomático a cero conforme \form{x} se acerca infinito \cite{oldham74}, un comportamiento que puede ser descrito por 

\formula{
\frac{1}{\gam{x}}\sim \frac{x^{\frac{1}{2}-x}}{\sqrt{2\pi}}e^x ,\ x\rightarrow \infty.
}

Como se ha visto, la función Gamma de un entero positivo \form{n} es un entero positivo, mientras que para un entero negativo es invariablemente infinito. Las funciones \form{\gam{\frac{1}{2}+n}} y \form{\gam{\frac{1}{2}-n}} resultan ser múltiplos de \form{\sqrt{\pi}}

\formula{
\gam{\frac{1}{2}}&=& \sqrt{\pi},\\
\gam{\frac{1}{2}+n}&=& \frac{(2n)!\sqrt{\pi}}{2^{2n}n!},\\
\gam{\frac{1}{2}-n}&=& \frac{(-1)^n2^{2n}n!\sqrt{\pi}}{(2n)!}.
}

Dos propiedades más de la función Gamma que resultan ser útiles son la reflexión

\formula{
\gam{x}=\frac{\pi \csc(\pi x)}{\gam{1-x}}=\frac{\pi}{\gam{1-x}\si{\pi x}},
}

y la duplicación

\formula{
\gam{2x}=\frac{2^{2x-1}}{\sqrt{\pi}}
\gam{x}\gam{x+\frac{1}{2}}
}

siendo esta última una instancia de la fórmula de multiplicación de Gauss

\formula{
\gam{nx}=\sqrt{\frac{2\pi}{n}}\left[\frac{n^x}{\sqrt{\pi}} \right]^n \prod_{k=0}^{n-1}\gam{x+\frac{k}{n}}.
}

\subsubsection{Diferenciación e Integración de Orden Entero}

En los cursos convencionales de cálculo diferencial estamos acostumbrados al uso de la notación 

\formula{
\der{d}{x}{n}f,
}

para la derivada \form{n}-ésima de una función \form{f} con respecto a la variable \form{x} cuando \form{n} es un entero no negativo. Dado que la integración y la diferenciación son operaciones inversas, es natural asociar la notación

\formula{
\der{d}{x}{-1}f,
}

a la integral indefinida de \form{f} con respecto a \form{x}.
Sin embargo, es necesario proporcionar un límite inferior de integración para que una integral indefinida este completamente determinada. Se suele asociar la notación anterior con el límite inferior cero. Por lo tanto, se define

\formula{
\ifr{0}{D}{x}{-1}f= \der{d}{x}{-1}f\equiv \int_0^x f(y)dy.
}

La integración múltiple con el límite inferior cero puede ser simbolizada por

\formula{
\ifr{0}{D}{x}{-2}f&\equiv &\int_0^x dx_1 \int_0^{x_1}f(x_0)dx_0,\\
\ifr{0}{D}{x}{-3}f & \equiv &\int_0^x dx_2 \int_0^{x_2}dx_1 \int_0^{x_1}f(x_0)dx_0,\\
&\vdots&\\
\ifr{0}{D}{x}{-n}f&\equiv& \int_0^x dx_{n-1}\int_0^{x_{n-1}}dx_{n-2}\int_0^{x_{n-2}}dx_{n-3}\cdots \int_0^x dx_1 \int_0^{x_1}f(x_0)dx_0.
}

Teniendo en cuenta la identidad \cite{oldham74}

\formula{
\int_a^x f(y)dy=\int_0^{x-a}f(y+a)dy,
}

se puede extender el simbolismo anterior para casos en que el limite inferior de la integral sea menor a cero

\formula{
\ifr{a}{D}{x}{-1}f= \der{d}{(x-a)}{-1}f&\equiv & \int_a^x f(y)dy,\\
\ifr{a}{D}{x}{-2}f=\der{d}{(x-a)}{-2}f&\equiv & \int_a^x dx_1 \int_a^{x_1}f(x_0)dx_0,\\
&\vdots&\\
\ifr{a}{D}{x}{-n}f=\der{d}{(x-a)}{-n} f&\equiv& \int_a^x dx_{n-1}\int_a^{x_{n-1}}dx_{n-2}\int_a^{x_{n-2}}dx_{n-3}\cdots \int_a^x dx_1 \int_a^{x_1}f(x_0)dx_0.
}

Se debe tener cuidado con la equivalencia

\formula{
\der{d}{(x-a)}{n}=\der{d}{x}{n},
}

que es una característica de un operador local, ya que en general para órdenes negativos

\formula{
\der{d}{(x-a)}{-n}\neq \der{d}{x}{-n}.
}

El símbolo \form{f^{(n)}} tiene un uso frecuente en la literatura como abreviatura de \form{\der{d}{x}{n}f}. Asimismo, se utiliza ocasionalmente \form{f^{(-n)}} para simbolizar  una integral \form{n}-foleada de \form{f} con respecto a \form{x}, estando los límites inferiores sin especificar

\formula{
f^{(-n)}\equiv \int_{a_n}^xdx_{n-1}\int_{a_{n-1}}^{x_{n-1}}dx_{n-2}\int_{a_{n-2}}^{x_{n-2}}dx_{n-3}\cdots \int_{a_2}^{x_2}dx_1\int_{a_1}^{x_1}f(x_0)dx_0,
}

donde  \form{a_1,a_2,\cdots, a_n} son valores completamente arbitrarios. Sin embargo, cuando se tiene en cuenta una diferencia como \form{f^{(-n)}(x)-f^{(-n)}(a)} se asume que los límites inferiores \form{a_1,a_2,\cdots, a_n} de cada integral son los mismos. 

\subsubsection{Definiciones Generales Para Derivadas e Integrales}

Tomando la definición convencional de la primera derivada en términos de una diferencia hacia atrás

\formula{
\der{d}{x}{1}f= \frac{d}{dx}f(x)\equiv\lim_{\delta x \rightarrow 0}\cfrac{f(x)-f(x-\delta x)}{\delta x}.
}

De igual forma se pueden obtener

\formula{
\der{d}{x}{2}f&=& \frac{d^2}{dx^2}f(x)\\
&=& \lim_{\delta x \rightarrow 0}\cfrac{ \cfrac{d}{dx}f(x)-\cfrac{d}{dx}f(x-\delta x)}{\delta x}\\
&=& \lim_{\delta x \rightarrow 0}\cfrac{ \lim_{\delta x \rightarrow 0}\cfrac{f(x)-f(x-\delta x)}{\delta x}-\lim_{\delta x \rightarrow 0}\cfrac{f(x-\delta x)-f(x-2\delta x)}{\delta x}}{\delta x}\\
&\equiv &\lim_{\delta x \rightarrow 0}\cfrac{f(x)-2f(x-\delta x)+f(x-2\delta x)}{\delta x^2},
}

análogamente

\formula{
\der{d}{x}{3}f&=& \frac{d^3}{dx^3}f(x)\\
&=& \lim_{\delta x \rightarrow 0}\cfrac{ \cfrac{d^2}{dx^2}f(x)-\cfrac{d^2}{dx^2}f(x-\delta x)}{\delta x}\\
&=& \lim_{\delta x \rightarrow 0}\cfrac{ \lim_{\delta x \rightarrow 0}\cfrac{f(x)-2f(x-\delta x)+f(x-2\delta x)}{\delta x^2}-\lim_{\delta x \rightarrow 0}\cfrac{f(x-\delta x)-2f(x-2\delta x)+f(x-3\delta x)}{\delta x^2}}{\delta x}\\
&\equiv &\lim_{\delta x \rightarrow 0}\frac{f(x)-3f(x-\delta x)+3 f(x-2\delta x)-f(x-3\delta x)}{\delta x^3},
}

e.t.c., donde, se ha asumido que los límites indicados existen.

Nótese que cada derivada posee una evaluación más en la misma función que el orden de la derivada, los coeficientes se acumulan como coeficientes binomiales y alternan en signo. De lo anterior se asume que la fórmula general de la derivada para un entero positivo \form{n} es \cite{oldham74}

\formula{
\der{d}{x}{n}f\equiv \lim_{\delta x \rightarrow 0}\frac{1}{\delta x^n}\sum_{k=0}^n (-1)^k \binom{n}{k}f(x-k \delta x).
}

Si la derivada \form{n}-ésima de \form{f} existe, esta última ecuación define  \form{\der{d}{x}{n}f} como un límite sin restricciones, es decir, como un límite donde \form{\delta x} tiende a cero a través de valores sin restricciones. Para unificar esta fórmula con la que define una integral como el límite de una suma, es deseable definir las derivadas en términos de un límite restringido. Para ello, se elige \form{\delta_N x\equiv \frac{(x-a)}{N}, N=1,2,\cdots}, donde \form{a} es un número menor que \form{x} y desempeña un papel equivalente a un límite inferior. Entonces, si el límite no restringido existe también existe el límite restringido y son iguales, la \form{n}-ésima derivada puede definirse entonces como

\formula{
\der{d}{x}{n}f \equiv \lim_{\delta_N x \rightarrow 0} \frac{1}{(\delta_N x)^n}\sum_{k=0}^n(-1)^k \binom{n}{k}f(x-k\delta_N x).
}

Ahora ya que \form{\binom{n}{k}=0} si \form{k>n} cuando \form{n} es un entero, la ecuación anterior puede escribirse como 

\formul{
\der{d}{x}{n}f &\equiv & \lim_{\delta_N x \rightarrow 0} \frac{1}{(\delta_N x)^n} \sum_{k=0}^{N-1}(-1)^k \binom{n}{k}f(x-k\delta_N x) \nonumber \\ 
&\equiv &\lim_{N\to \infty} \frac{1}{\left(\frac{x-a}{N}\right)^n}\sum_{k=0}^{N-1}(-1)^k \binom{n}{k}f\left(x-k \left(\frac{x-a}{N}\right) \right).
}{\label{eq:17}}

La ecuación \eqref{eq:17} se tomara  como \form{\der{d}{x}{n}f} con el entendido de que el límite indicado existe en el sentido usual, sin restricciones.

Regresando nuestra atención a las derivadas e integrales, comenzamos con la definición usual de una integral como límite de una suma de Riemann \cite{oldham74}

\formula{
\ifr{a}{D}{x}{-1}f &= &\int_a^x f(y)dy,\\
&= & \lim_{\delta_N x \rightarrow 0}\delta_N x \left( f(x)+f(x-\delta_N x)+ f(x-2\delta_N x) +\cdots + f(a+\delta_N x) \right)\\
&\equiv & \lim_{\delta_N x \rightarrow 0}  \delta_N x \sum_{k=0}^{N-1}f(x-k \delta_N x),
}

donde \form{\delta_N x \equiv \frac{(x-a)}{N}}. La misma definición sobre una integral doble da

\formula{
\ifr{a}{D}{x}{-2}f &= & \int_a^x dx_1 \int_a^{x_1}f(x_0)dx_0\\
&= & \lim_{\delta x \rightarrow 0} (\delta_N x)^2 \left( f(x)+2f(x-\delta_N x)+3 f(x-2\delta_N x)+\cdots + N f(a+\delta_N x)\right)\\
&\equiv & (\delta_N x)^2\sum_{k=0}^{N-1}(k+1)f(x-k \delta_N x),
}

tomando una iteración más para obtener una imagen más clara de la fórmula general

\formula{
\ifr{a}{D}{x}{-3}f&=& \int_a^x dx_2 \int_a^{x_2}dx_1 \int_a^{x_1}f(x_0)dx_0\\
&=& \lim_{\delta x \rightarrow 0} (\delta_N x)^3 \left(  f(x)+ 3 f(x-\delta_N x)+ 6 f(x-2\delta_N x)+\cdots + \frac{N(N+1)}{2} f(a+\delta_N x)\right)\\
&\equiv & \lim_{\delta_N x \rightarrow 0}(\delta_N x)^3\sum_{k=0}^{N-1}\frac{(k+1)(k+2)}{2}f(x-k\delta_N x).
}

Se observa que los coeficientes son construidos de la forma \form{\binom{k+n-1}{k}}, donde \form{n} es el orden de la integral, y todos los signos son positivos. Por lo tanto,

\formul{
\ifr{a}{D}{x}{-n}f&\equiv &\lim_{\delta_N x \rightarrow 0}(\delta_N x)^n \sum_{k=0}^{N-1}\binom{k+n-1}{k}f(x-k\delta_N x)\nonumber \\
&\equiv & \lim_{N\to \infty}  \left(\frac{x-a}{N} \right)^n \sum_{k=0}^{N-1}\binom{k+n-1}{k}f\left(x-k\left( \frac{x-a}{N}\right) \right).
}{\label{eq:18}}

Comparando los coeficientes de las ecuaciones \eqref{eq:17} y \eqref{eq:18}

\formula{
(-1)^k \binom{n}{k}=\binom{k-n-1}{k}=\frac{\gam{k-n}}{\gam{k+1}\gam{-n}},
}

se puede construir una ecuación general que involucre tanto a la derivada como a la integral

\formula{
\ifr{a}{D}{x}{q}f\equiv \lim_{N\to \infty} \frac{1}{\gam{-q}}\frac{1}{\left( \frac{x-a}{N}\right)^q} \sum_{k=0}^{N-1}\frac{\gam{k-q}}{\gam{k+1}}f\left(x-k\left( \frac{x-a}{N}\right) \right),
}

donde \form{q} es un entero de cualquiera signo. 

\subsubsection{Función de Mittag-Leffler}

Una función que es de gran utilidad en la resolución de ecuaciones diferenciales fraccionarias es la función de Mittag-Leffler, esta función actúa de forma similar a como lo hace la función exponencial en la resolución de ecuaciones diferenciales ordinarias, esta función se define por

\formul{
E_{\alpha,\beta}(t)= \sum_{k=0}^\infty \frac{t^k}{\gam{\alpha k + \beta}}, \ \ \alpha,\beta>0,
}{\label{con6}}

\begin{figure}[!ht]
\includegraphics[height=0.35\textheight, width=\textwidth]{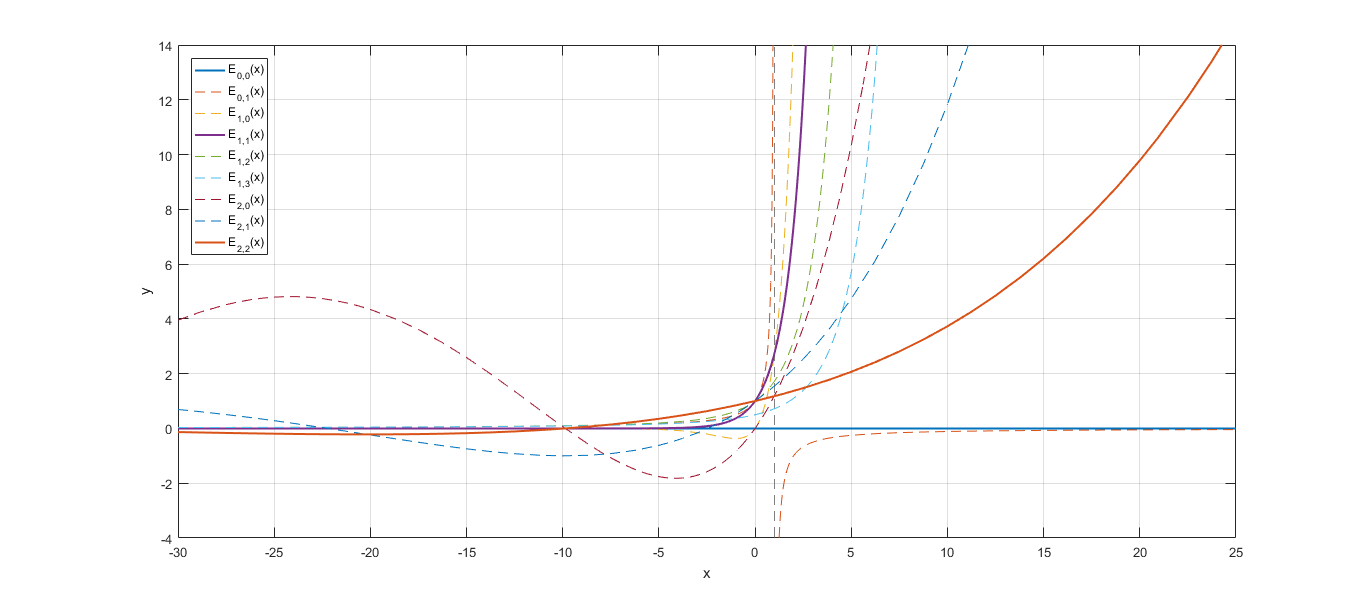}
\caption{Gráfica de la función de Mittag-Leffler $E_{\alpha,\beta}(x)$  para diferentes valores de $\alpha$ y $\beta$.}
\end{figure}

para obtener la transformada de Laplace de la ecuación de Mittag-Leffler tomemos la función para el parámetro $at^\alpha$ con $a=cte.$

\formula{
E_{\alpha,\beta}(at^\alpha)&=& \sum_{k=0}^\infty \frac{(at^\alpha)^k}{\gam{\alpha k + \beta}}\\
&=& \sum_{k=0}^\infty \frac{a^kt^{\alpha k}}{\gam{\alpha k + \beta}},
}

multiplicando por $t^{\beta-1}$ la expresión anterior se obtiene

\formula{
t^{\beta-1}E_{\alpha,\beta}(at^\alpha)
&=& \sum_{k=0}^\infty \frac{a^kt^{\beta-1}t^{\alpha k}}{\gam{\alpha k + \beta}}\\
&=& \sum_{k=0}^\infty \frac{a^kt^{\alpha k + \beta -1 }}{\gam{\alpha k + \beta}},
}

aplicando la transformada de Laplace y utilizando \eqref{con3}

\formula{
\mathcal{L}\left\{ t^{\beta-1}E_{\alpha,\beta}(at^\alpha)\right\}
&=& \sum_{k=0}^\infty \frac{a^k}{\gam{\alpha k + \beta} } \mathcal{L} \left\{ t^{\alpha k + \beta -1 } \right\}\\
&=& \sum_{k=0}^\infty \frac{a^k}{\gam{\alpha k + \beta} } \frac{\gam{\alpha k+\beta}}{s^{\alpha k+\beta}}\\
&=& \sum_{k=0}^\infty  \frac{a^k}{s^{\alpha k+\beta}}\\
\\
&=& s^{-\beta}\sum_{k=0}^\infty \left( \frac{a}{s^{\alpha }}\right)^k
\\
&=& s^{-\beta}\frac{1}{1-  \dfrac{a}{s^{\alpha }}},
}

por tanto

\formul{
\mathcal{L}\left\{ t^{\beta-1}E_{\alpha,\beta}(at^\alpha)\right\}=\frac{s^{\alpha-\beta}}{s^{\alpha}-a},
}{\label{con7}}

tomando $\beta=1$ obtenemos la transformada de Laplace para la función de Mittag-Leffler de un parámetro

\formul{
\mathcal{L}\left\{ E_{\alpha}(at^\alpha)\right\}= \mathcal{L}\left\{ \sum_{k=0}^\infty \frac{(at^\alpha)^k}{\gam{\alpha k+1}} \right\}=\frac{s^{\alpha-1}}{s^{\alpha}-a},
}{\label{con8}}

analizando la derivada de la función de Mittag-Leffler de un parámetro

\formula{
 \frac{d}{dt} E_{\alpha}(at^\alpha)
&=& \frac{d}{dt} \sum_{k=0}^\infty  \frac{a^kt^{\alpha k}}{\gam{\alpha k + 1}}\\
&=&  \sum_{k=1}^\infty  \frac{a^k}{\gam{\alpha k + 1}}\frac{d}{dt}t^{\alpha k}\\
&=&  \sum_{k=1}^\infty  \frac{a^k}{\gam{\alpha k + 1}}\frac{\gam{\alpha k+1}}{\gam{ \alpha k}}t^{\alpha k-1 }\\
&=&  \sum_{k=1}^\infty  \frac{a^k}{\gam{\alpha k }}t^{\alpha k-1 }\\
&=&  \sum_{k=0}^\infty  \frac{a^{k+1}}{\gam{\alpha k+\alpha }}t^{\alpha k+\alpha -1 }\\
&=& at^{\alpha-1} E_{\alpha,\alpha}(at^\alpha),
}

entonces

\formul{
 \frac{d}{dt} E_{\alpha}(\alpha t^\alpha)&=& at^{\alpha-1} E_{\alpha,\alpha}(at^\alpha),
}{\label{con9}}

analizando la derivada de la función de Mittag-Leffler de un parámetro

\formula{
 \frac{d}{dt} E_{\alpha}(at^\alpha)
&=& \frac{d}{dt} \sum_{k=0}^\infty  \frac{a^kt^{\alpha k}}{\gam{\alpha k + 1}}\\
&=&  \sum_{k=1}^\infty  \frac{a^k}{\gam{\alpha k + 1}}\frac{d}{dt}t^{\alpha k}\\
&=&  \sum_{k=1}^\infty  \frac{a^k}{\gam{\alpha k + 1}}\frac{\gam{\alpha k+1}}{\gam{ \alpha k}}t^{\alpha k-1 }\\
&=&  \sum_{k=1}^\infty  \frac{a^k}{\gam{\alpha k }}t^{\alpha k-1 }\\
&=&  \sum_{k=0}^\infty  \frac{a^{k+1}}{\gam{\alpha k+\alpha }}t^{\alpha k+\alpha -1 }\\
&=& at^{\alpha-1} E_{\alpha,\alpha}(at^\alpha),
}

entonces

\formul{
 \frac{d}{dt} E_{\alpha}(\alpha t^\alpha)&=& at^{\alpha-1} E_{\alpha,\alpha}(at^\alpha),
}{\label{con9}}

\subsubsection{Integral Iterada}

Sea $f(x)$ una función continua de $x\in \mathds{R}$, entonces se puede definir

\formula{
F(x)=\int_a^x f(t)dt,
}

integrando la función $F(x)$ se obtiene

\formul{
\int_a^x F(t)dt = \int_a^x \left(\int_a^tf(s)ds\right) dt,
}{\label{int1}}

realizando una integración por partes  tomando

\formula{
\begin{array}{lcl}
u=F(t) &\rightarrow & du=\dfrac{d}{dt}F(t)dt=f(t)dt,\\
dv=dt &\rightarrow & v=t,
\end{array}
}

entonces

\formula{
\int_a^x F(t)dt&=& vF(t)\Big{|}_a^x-\int_a^x tf(t)dt\\
&=& xF(x)-\int_a^x tf(t)dt\\
&=&x\int_a^x f(t)dt-\int_a^x tf(t)dt\\
&=&\int_a^x (x-t)f(t)dt
}

por lo tanto la integral de \eqref{int1} es

\formul{
\int_a^x F(t)dt =\int_a^x \left(\int_a^t f(s)ds \right)dt = \int_a^x (x-t)f(t)dt.
}{\label{int2}}

Definiendo ahora

\formula{
G(x)=\int_a^x F(t)dt,
}

integrando la función $G(x)$ se obtiene

\formula{
\int_a^x G(t)dt = \int_a^x \left(\int_a^tF(s)ds\right) dt=\int_a^x\left(\int_a^t \left(\int_a^s f(r)dr \right)ds \right)dt,
}

realizando una integración por partes  tomando

\formula{
\begin{array}{lcl}
u=G(t) &\rightarrow & du=\dfrac{d}{dt}G(t)dt=F(t)dt,\\
dv=dt &\rightarrow & v=t,
\end{array}
}

entonces

\formula{
\int_a^x G(t)dt&=& vG(t)\Big{|}_a^x-\int_a^x tF(t)dt\\
&=& xG(x)-\int_a^x tF(t)dt\\
&=&x\int_a^x F(t)dt-\int_a^x tF(t)dt\\
&=&\int_a^x (x-t)F(t)dt
}

realizando una integración por partes  tomando

\formula{
\begin{array}{lcl}
u=F(t) &\rightarrow & du=\dfrac{d}{dt}F(t)dt=f(t)dt,\\
dv=(x-t)dt &\rightarrow & v=-\dfrac{(x-t)^2}{2} ,
\end{array}
}

entonces

\formula{
\int_a^x G(t)dt&=& -\dfrac{(x-t)^2}{2}F(t)\Big{|}_a^x-\int_a^x -\dfrac{(x-t)^2}{2}f(t)dt\\
&=& \frac{1}{2}\int_a^x (x-t)^2f(t)dt,
}

por lo tanto la integral de \eqref{int2} es

\formul{
\int_a^x G(t)dt = \int_a^x \left(\int_a^tF(s)ds\right) dt=\int_a^x\left(\int_a^t \left(\int_a^s f(r)dr \right)ds \right)dt = \frac{1}{2}\int_a^x (x-t)^2f(t)dt
}{\label{int3}}

con el desarrollo anterior podemos deducir 

\formula{
\int_a^x\left(\int_a^t \left(\int_a^s\left( \int_a^x f(n)dn\right) dr \right)ds \right)dt = \frac{1}{3\cdot 2}\int_a^x (x-t)^3f(t)dt
}

tomando como $\ifr{a}{I}{x}{}$ el operador integral definido en el intervalo $(a,x)$ podemos escribir

\formula{
\ifr{a}{I}{x}{4}f(x)=\int_a^x\left(\int_a^t \left(\int_a^s\left( \int_a^x f(n)dn\right) dr \right)ds \right)dt = \frac{1}{3!}\int_a^x (x-t)^3f(t)dt,
}

lo que nos lleva a deducir una formula para la $n$-ésima integral de la función $f(x)$

\formul{
\ifr{a}{I}{x}{n}f(t) = \frac{1}{(n-1)!}\int_a^x (x-t)^{n-1}f(t)dt.
}{\label{int4}}

Procedemos a demostrar mediante el proceso de inducción que la ecuación anterior se cumple para todo $n\in\mathds{N}$. Como los casos para $n=1$ y $n=2$ se obtuvieron durante la construcción de \eqref{int4} procedemos a suponer que la ecuación es valida para $n=k$ con $k\in\mathds{N}$

\formula{
\ifr{a}{I}{x}{k}f(t) = \frac{1}{(k-1)!}\int_a^x (x-t)^{k-1}f(t)dt,
}

entonces aplicando nuevamente el operador integral y utilizando álgebra de operadores obtenemos

\formula{
\ifr{a}{I}{x}{}\left(\ifr{a}{I}{x}{k}f(t)\right) &=& \ifr{a}{I}{x}{}\left(\ifr{a}{I}{x}{k-1}\ifr{a}{I}{x}{}f(t)\right)= \ifr{a}{I}{x}{k}\left(\ifr{a}{I}{x}{}f(t)\right)\\
&=& \frac{1}{(k-1)!}\int_a^x(x-t)^{k-1}\left(\ifr{a}{I}{t}{}f(s) \right)dt,
}

realizando una integración por partes  tomando

\formula{
\begin{array}{lcl}
u=\ifr{a}{I}{t}{}f(s) &\rightarrow & du=f(t)dt,\\
dv=(x-t)^{k-1}dt &\rightarrow & v=-\dfrac{(x-t)^k}{k} ,
\end{array}
}

entonces

\formula{
\ifr{a}{I}{x}{}\left(\ifr{a}{I}{x}{k}f(t)\right) &=& \frac{1}{(k-1)!}\int_a^x(x-t)^{k-1}\left(\ifr{a}{I}{t}{}f(s) \right)dt\\&=& \frac{1}{(k-1)!} \left[ - \left.\dfrac{(x-t)^k}{k} \ifr{a}{I}{t}{}f(s) \right|_a^x -\int_a^x -\dfrac{(x-t)^k}{k} f(t)dt\right]\\
&=& \frac{1}{k!}\int_a^x (x-t)^{k}f(t)dt\\
&=& \ifr{a}{I}{x}{k+1}f(t),
}

lo cual demuestra que la ecuación \eqref{int4} es valida para todo $n\in \mathds{N}$. Tomando el hecho de que la función Gamma esta relacionada con el factorial por la igualdad $\gam{n}=(n-1)!$, podes escribir

\formula{
\ifr{a}{I}{x}{n}f(t) = \frac{1}{\gam{n}}\int_a^x (x-t)^{n-1}f(t)dt, \ \ n\in\mathds{R}\setminus (\mathds{Z}^{-}\cup \{0\}),
}

como el recíproco de la función Gamma, conocido como el producto infinito de Weierstrass, es univaludo y finito para todo \form{z\in \mathds{C} \setminus (\mathds{Z}^{-}\cup \{0\})} \cite{arfken85}, podemos reescribir la ecuación \eqref{int4} como

\formul{
\ifr{a}{I}{x}{\alpha}f(t) = \frac{1}{\gam{\alpha}}\int_a^x (x-t)^{\alpha-1}f(t)dt, \ \ \alpha\in\mathds{C}\setminus (\mathds{Z}^{-}\cup \{0\}),
}{\label{int5}}

por otro lado si se considera el intervalo $(x,b)$ y se lleva a cabo el mismo desarrollo se obtiene 

\formul{
\ifr{x}{I}{b}{\alpha}f(t) = \frac{1}{\gam{\alpha}}\int_x^b (t-x)^{\alpha-1}f(t)dt, \ \ \alpha\in\mathds{C}\setminus (\mathds{Z}^{-}\cup \{0\}).
}{\label{int6}}

\subsubsection{Operadores Diferointegrables}

Comencemos por analizar el comportamiento de las funciones trigonométricas seno y coseno para las derivadas e integrales de orden entero, para eso primero hay que notar como están relacionadas dichas funciones

\formula{
\si{x\pm \frac{n}{2}\pi}&=&\si{x}\co{ \frac{n}{2}\pi}\pm \si{ \frac{n}{2}\pi}\co{x},\\
\co{x \pm \frac{n}{2}\pi}&=& \co{x}\co{ \frac{n}{2}\pi}\mp \si{x}\si{ \frac{n}{2}\pi},
}

sin perdida de generalidad tomemos la función seno y obtengamos sus primeras derivadas

\formula{
\begin{array}{ccccc}
\der{d}{x}{}\si{x}&=&\co{x} &=&\si{x+\dfrac{1}{2}\pi},\\
\der{d}{x}{2}\si{x}&=& -\si{x}&=& \si{x+\pi},\\
\der{d}{x}{3}\si{x}&=& -\co{x}&=& \si{x+\dfrac{3}{2}\pi},\\
\der{d}{x}{4}\si{x}&=& \si{x}&=& \si{x+2\pi },
\end{array}
}

\begin{figure}[!ht]
\includegraphics[height=0.35\textheight, width=\textwidth]{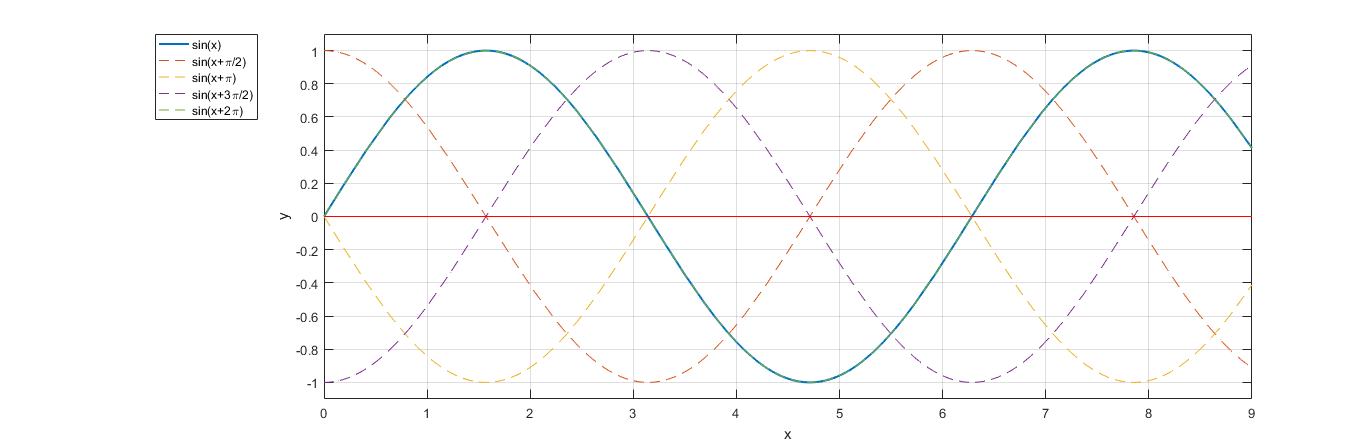}
\caption{Comportamiento de la función Seno para las derivadas enteras $\left(\frac{d}{dx} \right)^n\sin(x)$ con $n=0,1,2,3,4$}
\end{figure}

tomando ahora las integrales impropias se obtiene

\formula{
\begin{array}{ccccc}
\displaystyle{\int\si{x}dx}&=&\displaystyle{-\co{x}+a_0} &=&\displaystyle{ \si{x-\dfrac{1}{2}\pi}+a_0},\\
\displaystyle{ \left( \int \right)^2\si{x} d^2x}&=& \displaystyle{ -\si{x}+\sum_{k=0}^1 a_kx^k}&=& \displaystyle{\si{x-\pi}+\sum_{k=0}^1 a_kx^k},\\
\displaystyle{\left( \int\right)^3 \si{x} d^3x}&=& \displaystyle{\co{x}+\sum_{k=0}^2 a_kx^k}&=&\displaystyle{ \si{x-\dfrac{3}{2}\pi}+\sum_{k=0}^2 a_kx^k},\\
\displaystyle{ \left( \int \right)^4 \si{x} d^4x} &=&\displaystyle{ \si{x}+\sum_{k=0}^3 a_kx^k}&=&\displaystyle{ \si{x-2\pi }+\sum_{k=0}^3 a_kx^k},
\end{array}
}

se puede deducir entonces las formula para la $n$-ésima derivada y la $n$-ésima integral de la función seno

\formula{
\der{d}{x}{n}\si{x}&=& \si{x+\frac{n}{2}\pi},\\
\left(\int \right)^n \si{x}d^nx&=& \si{x-\frac{n}{2}\pi}+\sum_{k=0}^{n-1} a_kx^k,
}

de las ecuaciones anteriores obtenemos que el operador derivada induce un corrimiento hacia la derecha para la función seno mientas que el operador integral induce un corrimiento hacia la izquierda,  considerando que la derivada es el operador inverso por la izquierda de la integral podemos definir

\formula{
\der{d}{x}{n}\left(\int \right)^n f(x) d^nx &=& D^n I^n f(x)\\
&=& I^{-n}I^{n}f(x)\\
&=& D^n D^{-n}f(x)\\
&=& f(x),
}

de lo anterior se puede escribir

\formul{
\ifr{a}{I}{x}{\alpha}f(t)&=& \ifr{a}{D}{x}{-\alpha}f(t) \nonumber \\
&=&  D^n \ifr{a}{I}{x}{n} \ifr{a}{I}{x}{\alpha}f(t) =D^n \ifr{a}{I}{x}{n+\alpha}f(t) \nonumber \\
&=& \frac{1}{\gam{n+\alpha}} \der{d}{x}{n} \int_a^x(x-t)^{n+\alpha-1}f(t)dt
}{\label{int7}}

aunque en general la integral no es un operador por la izquierda de la derivada ( tómese por ejemplo un polinomio de grado $k<n$), para las funciones en las cuales se satisface $I^n D^n f(x)=f(x)$ podemos escribir

\formul{
\ifr{a}{I}{x}{\alpha}f(t)&=& \ifr{a}{D}{x}{-\alpha}f(t) \nonumber \\
&=&   \ifr{a}{I}{x}{\alpha} \ifr{a}{I}{x}{n} D^n f(t) = \ifr{a}{I}{x}{n+\alpha} D^n f(t)\nonumber \\
&=& \frac{1}{\gam{n+\alpha}}  \int_a^x(x-t)^{n+\alpha-1} \der{d}{t}{n} f(t)dt,
}{\label{int8}}

por ultimo se puede restringir el valor de $n$ utilizando la función piso, finalmente  tomando $n=\lfloor \alpha \rfloor +1$ y bajo el símbolo del operador derivada las ecuaciones \eqref{int7} y \eqref{int8} se pueden reescribir como

\formul{
\ifr{a}{D}{x}{\alpha}f(t) &=& \frac{1}{\gam{n-\alpha}} \der{d}{x}{n} \int_a^x  (x-t)^{n-\alpha-1}f(t)dt , \label{int9}\\
\ifrR{C}{a}{D}{x}{\alpha}f(t) 
&=& \frac{1}{\gam{n-\alpha}}  \int_a^x(x-t)^{n-\alpha-1} \der{d}{t}{n} f(t)dt \label{int10},
}{}

las ecuaciones anteriores se conocen como derivada fraccionarias de Riemann-Liouville y derivada fraccionaria de Caputo respectivamente.


Tomando $\alpha\in(0,1)$ y $f(x)=x$, se obtiene para la ecuación \eqref{int9}

\formula{
\ifr{a}{D}{x}{\alpha}f(t) &=& \frac{1}{\gam{1-\alpha}} \frac{d}{dx} \int_a^x  (x-t)^{-\alpha}t dt,
}

realizando una integración por partes  tomando

\formula{
\begin{array}{lcl}
u=t &\rightarrow & du=dt,\\
dv=(x-t)^{-\alpha}dt &\rightarrow & v=-\dfrac{1}{1-\alpha}(x-t)^{1-\alpha} ,
\end{array}
}

entonces

\formula{
\ifr{a}{D}{x}{\alpha}f(t)
&=& \frac{1}{\gam{1-\alpha}} \frac{d}{dx} \int_a^x  (x-t)^{-\alpha}t dt\\
&=& \frac{1}{\gam{1-\alpha}} \frac{d}{dx}  \left[ -\frac{1}{1-\alpha} t(x-t)^{1-\alpha} \Big{|}_a^x +\frac{1}{1-\alpha} \int_a^x  (x-t)^{1-\alpha}dt \right]\\
&=& \frac{1}{\gam{1-\alpha}} \frac{d}{dx}  \left[ \frac{1}{1-\alpha} a(x-a)^{1-\alpha}  -\frac{1}{1-\alpha} \frac{1}{2-\alpha}(x-t)^{2-\alpha}\Big{|}_{a}^x \right]\\
&=& \frac{1}{\gam{1-\alpha}} \frac{d}{dx}  \left[ \frac{1}{1-\alpha} a(x-a)^{1-\alpha}  +\frac{1}{(1-\alpha)(2-\alpha)}(x-a)^{2-\alpha} \right]\\
&=& \frac{1}{\gam{1-\alpha}}   \left[  a(x-a)^{-\alpha}  +\frac{1}{1-\alpha}(x-a)^{1-\alpha} \right]\\
&=& \frac{1}{\gam{1-\alpha}(1-\alpha)} (x-a)^{-\alpha}   \left[  a(1-\alpha)  +(x-a) \right]\\
&=& \frac{1}{\gam{2-\alpha}}(x-a)^{-\alpha} (x-\alpha a),
}

mientras que para la ecuación \eqref{int10} se obtiene

\formula{
\ifrR{C}{a}{D}{x}{\alpha}f(t) 
&=& \frac{1}{\gam{1-\alpha}}  \int_a^x(x-t)^{-\alpha} \frac{d}{dt}t dt\\
&=& \frac{1}{\gam{1-\alpha}}  \int_a^x(x-t)^{-\alpha}  dt\\
&=& \frac{-1}{\gam{1-\alpha}}\frac{1}{1-\alpha} (x-t)^{1-\alpha} \Big{|}_a^x\\
&=& \frac{1}{\gam{2-\alpha}}(x-a)^{1-\alpha},
}

tomando el intervalo $(0,x)$ se obtiene

\formula{
\ifr{0}{D}{x}{\alpha}f(t)=\ifrR{C}{0}{D}{x}{\alpha}f(t)=\frac{1}{\gam{2-\alpha}}x^{1-\alpha}.
}


\begin{figure}[!ht]
\includegraphics[height=0.35\textheight, width=\textwidth]{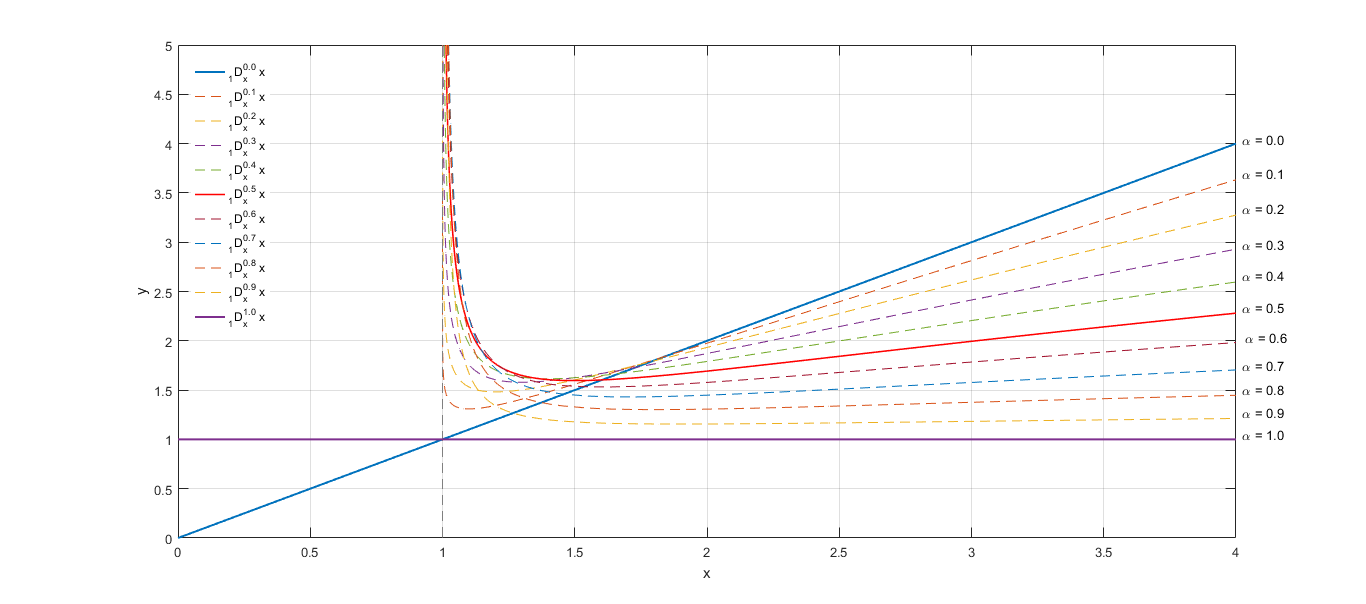}
\caption{Derivada fraccionaria de Riemann-Liouville de la función identidad en el intervalo $(1,x)$ para valores de $\alpha \in[0,1]$.}
\end{figure}

\begin{figure}[!ht]
\includegraphics[height=0.35\textheight, width=\textwidth]{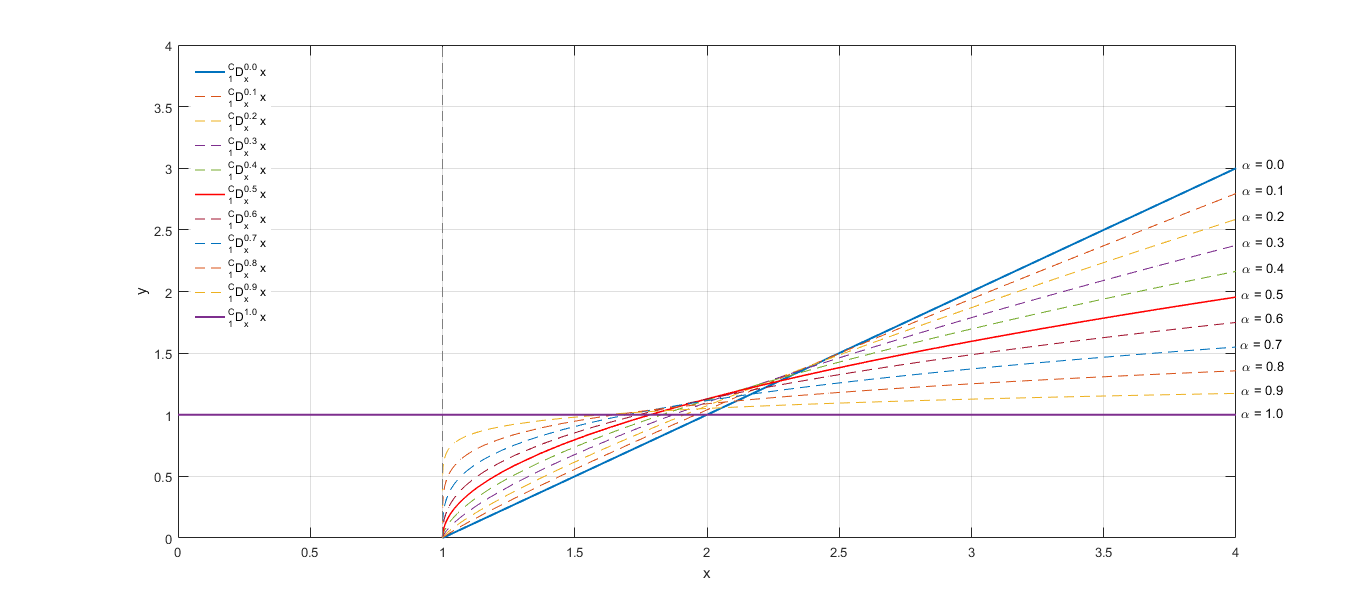}
\caption{Derivada fraccionaria de Caputo de la función identidad en el intervalo $(1,x)$ para valores de $\alpha \in[0,1]$.}
\end{figure}

\begin{figure}[!ht]
\includegraphics[height=0.35\textheight, width=\textwidth]{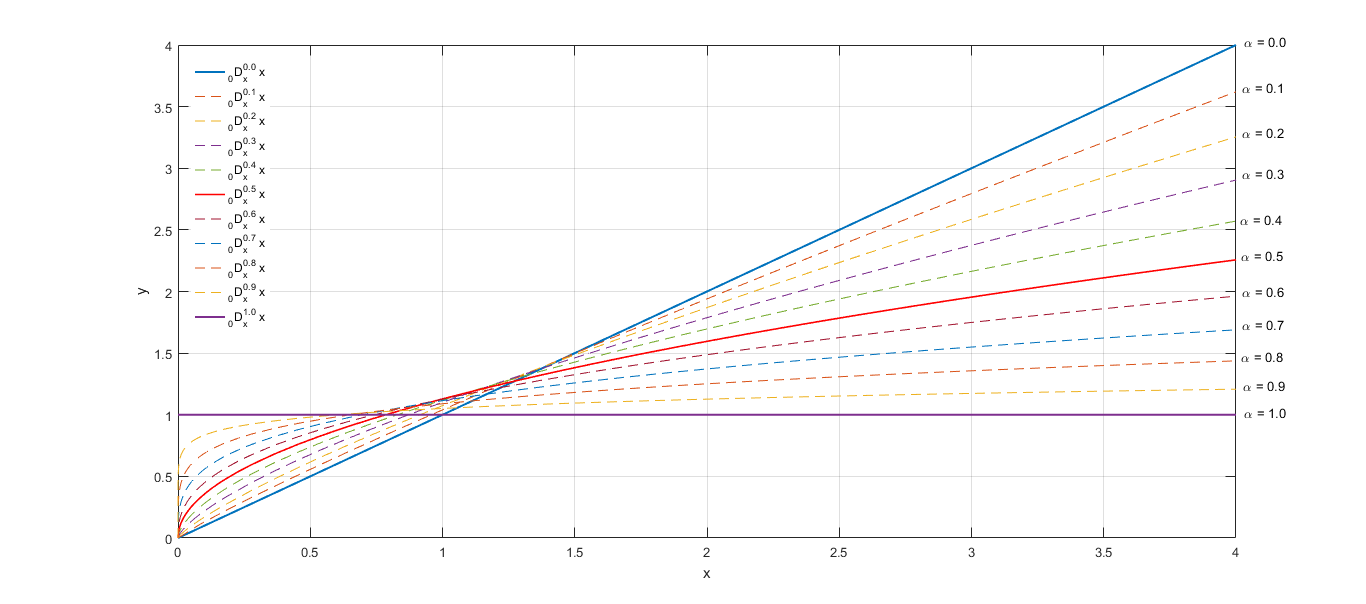}
\caption{Derivada fraccionaria de la función identidad en el intervalo $(0,x)$ para valores de $\alpha \in[0,1]$.}
\end{figure}


Ahora obtengamos un resultado un poco mas general, tomando $\alpha\in(0,1)$ y $f(x)=(x-c_0)^m$, se obtiene para la ecuación \eqref{int9}

\formula{
\ifr{a}{D}{x}{\alpha}f(t) &=& \frac{1}{\gam{1-\alpha}} \frac{d}{dx} \int_a^x  (x-t)^{-\alpha}(t-c_0)^m dt,
}

tomando el cambio de variable $t=c_0+(x-c_0)u$

\formula{
\ifr{a}{D}{x}{\alpha}f(t) 
&=& \frac{1}{\gam{1-\alpha}} \dfrac{d}{dx}\left\{ \int_\frac{a-c_0}{x-c_0}^1  \left( (x-c_0)-(x-c_0)u\right)^{-\alpha}\left( (x-c_0)u \right)^m (x-c_0)du  \right\}\\
&=& \frac{1}{\gam{1-\alpha}} \dfrac{d}{dx}\left\{ \left[ \int_0^1   (1-u)^{-\alpha} u^m du -\int_0^\frac{a-c_0}{x-c_0}   (1-u)^{-\alpha} u^m du\right](x-c_0)^{m-\alpha+1} \right\}\\
&=& \frac{1}{\gam{1-\alpha}} \dfrac{d}{dx}\left\{ \left[ B\left(1-\alpha,m+1 \right) -B\left(\dfrac{a-c_0}{x-c_0} ;1-\alpha,m+1 \right)\right](x-c_0)^{m-\alpha+1} \right\}\\
&=& \frac{1}{\gam{1-\alpha}} B\left(1-\alpha,m+1 \right)\dfrac{d}{dx}\left\{ \left[ 1 -\dfrac{ B\left(\dfrac{a-c_0}{x-c_0} ;1-\alpha,m+1 \right)}{B\left(1-\alpha,m+1 \right)}\right](x-c_0)^{m-\alpha+1} \right\}\\
&=& \frac{\gam{m+1}}{\gam{m-\alpha+2}} \dfrac{d}{dx}\left\{ \left[ 1 -\dfrac{ B\left(\dfrac{a-c_0}{x-c_0} ;1-\alpha,m+1 \right)}{B\left(1-\alpha,m+1 \right)}\right](x-c_0)^{m-\alpha+1} \right\},
}

entonces

\begin{eqnarray}
\ifr{a}{D}{x}{\alpha}(x-c_0)^m =\frac{\gam{m+1}}{(m-\alpha+1)\gam{m-\alpha+1}} \dfrac{d}{dx}\left\{ \left[ 1 -\dfrac{ B\left(\dfrac{a-c_0}{x-c_0} ;1-\alpha,m+1 \right)}{B\left(1-\alpha,m+1 \right)}\right](x-c_0)^{m-\alpha+1} \right\},
\end{eqnarray}

mientras que para la ecuación \eqref{int10} se obtiene

\formula{
\ifrR{C}{a}{D}{x}{\alpha}f(t) 
&=& \frac{1}{\gam{1-\alpha}}  \int_a^x(x-t)^{-\alpha} \frac{d}{dt}(t-c_0)^m dt\\
&=& \frac{1}{\gam{1-\alpha}}\frac{\gam{m+1}}{\gam{m}}  \int_a^x(x-t)^{-\alpha} (t-c_0)^{m-1}dt ,
}

tomando el cambio de variable $t=c_0+(x-c_0)u$

\formula{
\ifrR{C}{a}{D}{x}{\alpha}f(t) 
&=& \frac{\gam{m+1}}{\gam{1-\alpha}\gam{m}} \int_\frac{a-c_0}{x-c_0}^1((x-c_0)-(x-c_0)u)^{-\alpha} ((x-c_0)u)^{m-1}(x-c_0)du \\
&=& \frac{\gam{m+1}}{\gam{1-\alpha}\gam{m}}  \left[  \int_0^1 \left( 1-u\right)^{-\alpha} u^{m-1} du- \int_0^\frac{a-c_0}{x-c_0} \left( 1-u\right)^{-\alpha} u^{m-1} du\right](x-c_0)^{m-\alpha}  \\
&=& \frac{\gam{m+1}}{\gam{1-\alpha}\gam{m}} \left[ B\left( 1-\alpha,m\right)-  B\left(\dfrac{a-c_0}{x-c_0}; 1-\alpha,m \right)\right] (x-c_0)^{m-\alpha} \\
&=& \frac{\gam{m+1}}{\gam{1-\alpha}\gam{m}}B\left( 1-\alpha,m\right)  \left[ 1- \dfrac{ B\left(\dfrac{a-c_0}{x-c_0}; 1-\alpha,m \right)}{B\left( 1-\alpha,m\right)}\right](x-c_0)^{m-\alpha},
}

entonces

\begin{eqnarray}
\ifrR{C}{a}{D}{x}{\alpha}(x-c_0)^m = \frac{\gam{m+1}}{\gam{m-\alpha+1}} \left[ 1- \dfrac{ B\left(\dfrac{a-c_0}{x-c_0}; 1-\alpha,m \right)}{B\left( 1-\alpha,m\right)}\right](x-c_0)^{m-\alpha},
\end{eqnarray}

tomando $c_0=a$ se obtiene

\begin{eqnarray}
\ifr{a}{D}{x}{\alpha}(x-a)^m=\ifrR{C}{a}{D}{x}{\alpha}(x-a)^m=\frac{\gam{m+1}}{\gam{m-\alpha+1}}(x-a)^{m-\alpha}.
\end{eqnarray}


\begin{figure}[!ht]
\includegraphics[height=0.35\textheight, width=\textwidth]{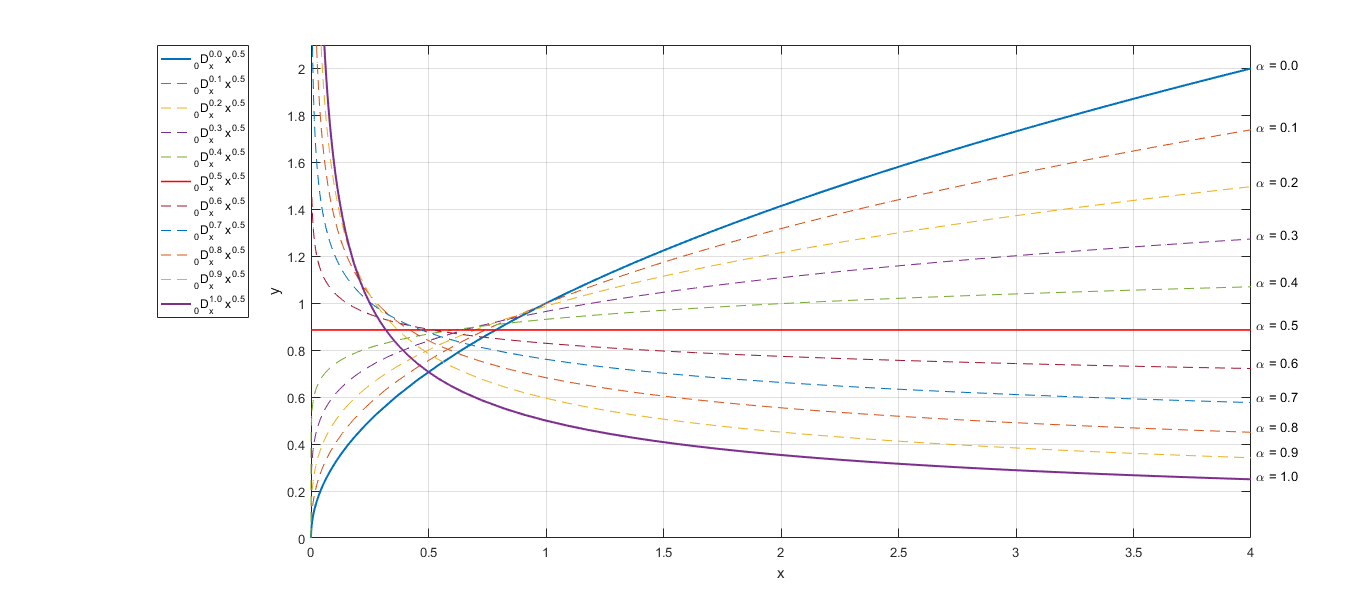}
\caption{Derivada fraccionaria de la función $x^\frac{1}{2}$ en el intervalo $(0,x)$ para valores de $\alpha \in[0,1]$.}
\end{figure}

\begin{figure}[!ht]
\includegraphics[height=0.35\textheight, width=\textwidth]{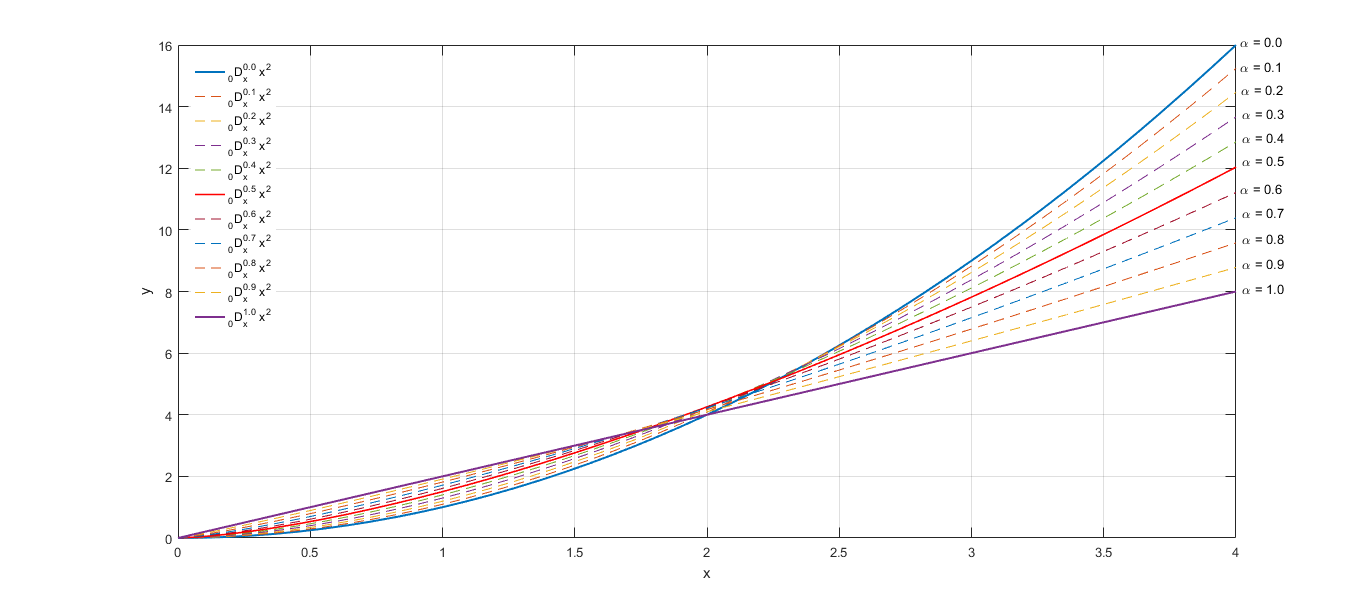}
\caption{Derivada fraccionaria de la función $x^2$ en el intervalo $(0,x)$ para valores de $\alpha \in[0,1]$.}
\end{figure}


Tomando $\alpha\in(0,1)$ y $f(x)=\si{c(x-c_0)}$, se obtiene para la ecuación \eqref{int9}

\formula{
\ifr{a}{D}{x}{\alpha}f(t) &=& \frac{1}{\gam{1-\alpha}} \frac{d}{dx} \int_a^x  (x-t)^{-\alpha}\si{c(t-c_0)} dt\\
&=& \frac{1}{\gam{1-\alpha}} \frac{d}{dx}\left[ \int_a^x  (x-t)^{-\alpha}\sum_{k=0}^\infty \dfrac{(-1)^k}{(2k+1)!}(c(t-c_0))^{2k+1}  dt \right]\\
&=& \frac{1}{\gam{1-\alpha}} \frac{d}{dx}\left\{ \sum_{k=0}^\infty \dfrac{(-1)^k}{\gam{2k+2}}c^{2k+1}   \int_a^x  (x-t)^{-\alpha}(t-c_0)^{2k+1}  dt  \right\},
}

tomando el cambio de variable $t=c_0+(x-c_0)u$

\formula{
\ifr{a}{D}{x}{\alpha}f(t) &=& \frac{1}{\gam{1-\alpha}} \frac{d}{dx}\left\{ \sum_{k=0}^\infty \dfrac{(-1)^k}{\gam{2k+2}}c^{2k+1}   \int_\frac{a-c_0}{x-c_0}^1  ((x-c_0)-(x-c_0)u)^{-\alpha}((x-c_0)u)^{2k+1}  (x-c_0)du  \right\} \\
&=& \frac{d}{dx}\left\{  \sum_{k=0}^\infty \dfrac{(-1)^k}{\gam{1-\alpha} \gam{2k+2}}  \left[ \int_0^1  \left( 1-u\right)^{-\alpha}u^{2k+1}  du  -\int_0^\frac{a-c_0}{x-c_0}  \left( 1-u\right)^{-\alpha}u^{2k+1}   du\right]c^{2k+1}(x-c_0)^{2k-\alpha+2} \right\} \\
&=& \frac{d}{dx}\left\{  \sum_{k=0}^\infty \dfrac{(-1)^k}{\gam{1-\alpha} \gam{2k+2}} \left[ B(1-\alpha,2k+2) - B\left(\dfrac{a-c_0}{x-c_0} , 1-\alpha,2k+2\right) \right] c^{2k+1}(x-c_0)^{2k+2-\alpha}  \right\} \\
&=& \frac{d}{dx}\left\{  \sum_{k=0}^\infty \dfrac{(-1)^k}{\gam{1-\alpha} \gam{2k+2}}  B(1-\alpha,2k+2)\left[ 1 - \dfrac{ B\left(\dfrac{a-c_0}{x-c_0} ; 1-\alpha,2k+2\right)}{B(1-\alpha,2k+2)} \right] c^{2k+1}(x-c_0)^{2k+2-\alpha}  \right\} \\
}

entonces

\begin{eqnarray}
\ifr{a}{D}{x}{\alpha}\si{c(x-c_0)} =
 \frac{d}{dx}\left\{  \sum_{k=0}^\infty \dfrac{(-1)^k}{ \gam{2k+3-\alpha}}\left[ 1 - \dfrac{ B\left(\dfrac{a-c_0}{x-c_0} ;1-\alpha,2k+2\right)}{B(1-\alpha,2k+2)} \right] c^{2k+1}(x-c_0)^{2k+2-\alpha}  \right\},
\end{eqnarray}

tomado  $c_0=a$ en la expresión anterior se obtiene

\begin{eqnarray}
\ifr{a}{D}{x}{\alpha}\si{c(x-a)} &=&
 \frac{d}{dx}\left\{  \sum_{k=0}^\infty \dfrac{(-1)^k}{ \gam{2k+3-\alpha}} c^{2k+1}(x-a)^{2k+2-\alpha}  \right\}\nonumber\\
&=& c \frac{d}{dx}\left\{ (x-a)^{2-\alpha} \sum_{k=0}^\infty \dfrac{(-1)^k}{ \gam{2k+3-\alpha}} c^{2k}(x-a)^{2k}  \right\}\nonumber\\
&=& c \frac{d}{dx}\left\{ (x-a)^{2-\alpha} E_{2,3-\alpha} \left(-c^2(x-a)^2 \right)  \right\},
\end{eqnarray}

\begin{figure}[!ht]
\includegraphics[height=0.35\textheight, width=\textwidth]{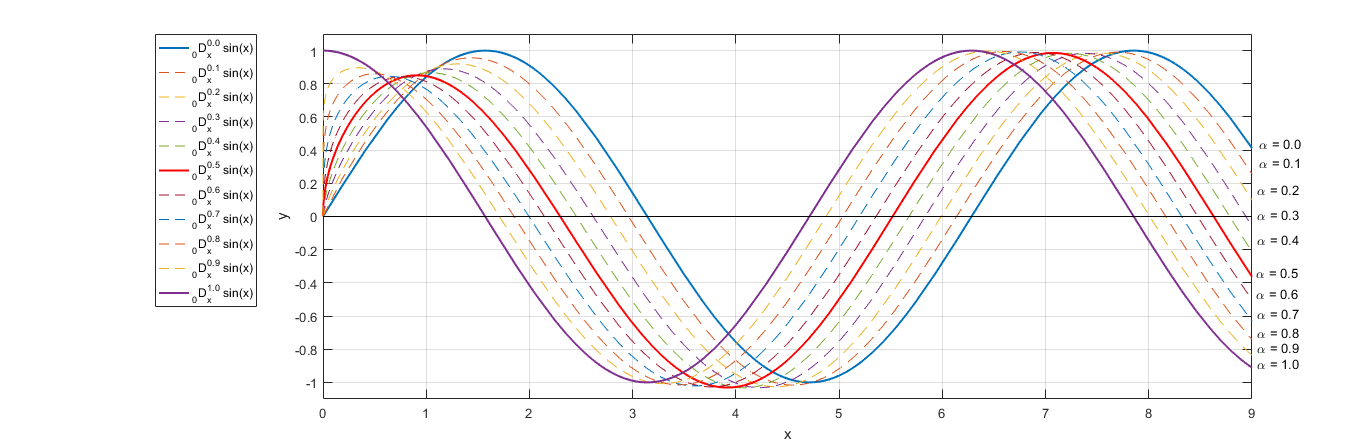}
\caption{Derivada fraccionaria de Riemann-Liouville de la función $\si{x}$ en el intervalo $(0,x)$ para valores de $\alpha \in[0,1]$.}
\end{figure}

mientras que para la ecuación \eqref{int10} se obtiene

\formula{
\ifrR{C}{a}{D}{x}{\alpha}f(t) 
&=& \frac{1}{\gam{1-\alpha}}  \int_a^x(x-t)^{-\alpha} \frac{d}{dt}\si{c(t-c_0)} dt\\
&=& \frac{c}{\gam{1-\alpha}}  \int_a^x(x-t)^{-\alpha} \co{c(t-c_0)} dt\\
&=& \frac{c}{\gam{1-\alpha}}  \int_a^x(x-t)^{-\alpha} \sum_{k=0}^\infty \frac{(-1)^k}{(2k)!}(c(t-c_0))^{2k} dt \\
&=& \frac{c}{\gam{1-\alpha}} \sum_{k=0}^\infty \frac{(-1)^k}{\gam{2k+1}}c^{2k} \int_a^x(x-t)^{-\alpha} (t-c_0)^{2k} dt ,
}

tomando el cambio de variable $t=c_0+(x-c_0)u$

\formula{
\ifrR{C}{a}{D}{x}{\alpha}f(t) 
&=& \frac{c}{\gam{1-\alpha}}  \sum_{k=0}^\infty \frac{(-1)^k}{\gam{2k+1}}c^{2k} \int_\frac{a-c_0}{x-c_0}^1((x-c_0)-(x-c_0)u)^{-\alpha} ((x-c_0)u)^{2k} (x-c_0)du\\
&=& c  \sum_{k=0}^\infty \frac{(-1)^k}{\gam{1-\alpha}\gam{2k+1}}\left[ \int_0^1 \left( 1-u\right)^{-\alpha} u^{2k} du -\int_0^\frac{a-c_0}{x-c_0} \left( 1-u\right)^{-\alpha} u^{2k} du \right]c^{2k}  (x-c_0)^{2k-\alpha+1}\\
&=& c  \sum_{k=0}^\infty \frac{(-1)^k}{\gam{1-\alpha}\gam{2k+1}}\left[ B(1-\alpha,2k+1)-B\left(\dfrac{a-c_0}{x-c_0};1-\alpha, 2k+1 \right) \right]c^{2k}  (x-c_0)^{2k-\alpha+1}\\
&=& c  \sum_{k=0}^\infty \frac{(-1)^k}{\gam{1-\alpha}\gam{2k+1}} B(1-\alpha,2k+1) \left[ 1-\dfrac{ B\left(\dfrac{a-c_0}{x-c_0};1-\alpha, 2k+1 \right)}{B(1-\alpha,2k+1)} \right]c^{2k}  (x-c_0)^{2k+1-\alpha},
}

entonces

\begin{eqnarray}
\ifrR{C}{a}{D}{x}{\alpha}\si{c(x-c_0)}  = c  \sum_{k=0}^\infty \frac{(-1)^k}{\gam{2k+2-\alpha}}  \left[ 1-\dfrac{ B\left(\dfrac{a-c_0}{x-c_0};1-\alpha, 2k+1 \right)}{B(1-\alpha,2k+1)} \right]c^{2k}  (x-c_0)^{2k+1-\alpha}.
\end{eqnarray}

tomado $c_0=a$ en la expresión anterior se obtiene

\begin{eqnarray}
\ifrR{C}{a}{D}{x}{\alpha}\si{c(x-a)}  &=& c  \sum_{k=0}^\infty \frac{(-1)^k}{\gam{2k+2-\alpha}}  c^{2k}  (x-a)^{2k+1-\alpha}\nonumber\\
&=& c (x-a)^{1-\alpha} \sum_{k=0}^\infty \frac{(-1)^k}{\gam{2k+2-\alpha}}  c^{2k}  (x-a)^{2k}\nonumber\\
&=& c (x-a)^{1-\alpha} E_{2,2-\alpha}  \left(- c^2  (x-a)^2\right).
\end{eqnarray}

\begin{figure}[!ht]
\includegraphics[height=0.35\textheight, width=\textwidth]{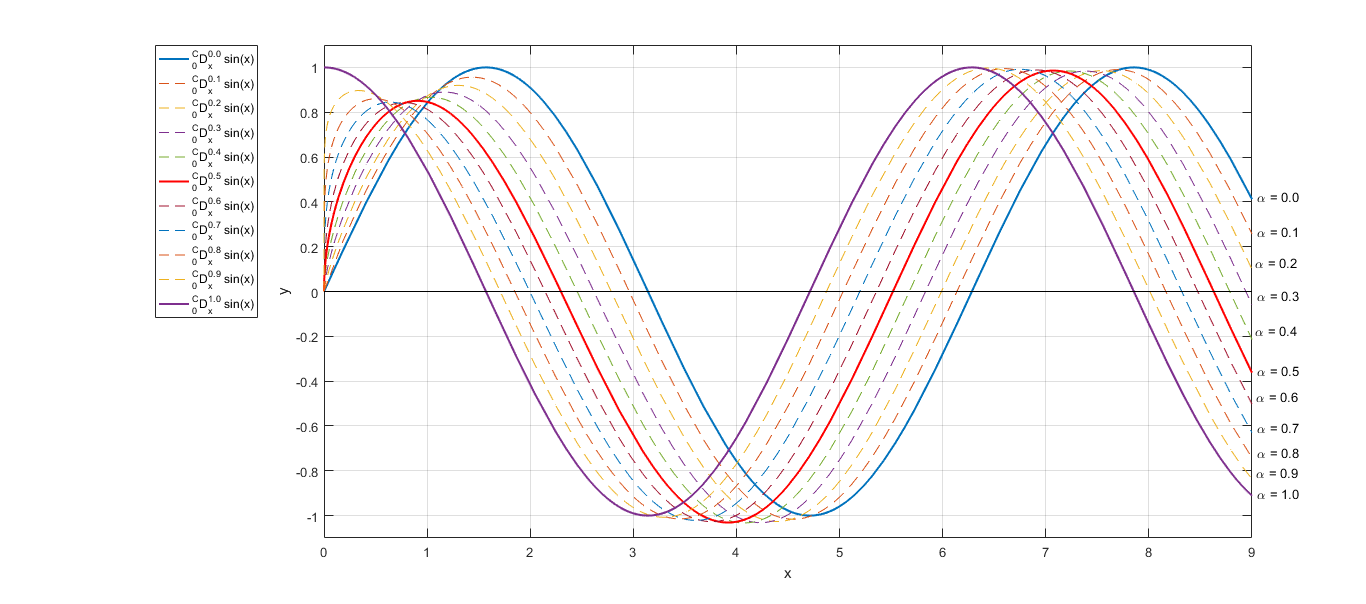}
\caption{Derivada fraccionaria de Caputo de la función $\si{x}$ en el intervalo $(0,x)$ para valores de $\alpha \in[0,1]$.}
\end{figure}

Tomando $\alpha\in(0,1)$ y $f(x)=\co{c(x-c_0)}$, se obtiene para la ecuación \eqref{int10}

\formula{
\ifrR{C}{a}{D}{x}{\alpha}f(t) 
&=& \frac{1}{\gam{1-\alpha}}  \int_a^x(x-t)^{-\alpha} \frac{d}{dt}\co{c(t-c_0)} dt\\
&=& \frac{-c}{\gam{1-\alpha}}  \int_a^x(x-t)^{-\alpha} \si{c(t-c_0)} dt\\
&=& \frac{-c}{\gam{1-\alpha}}  \int_a^x(x-t)^{-\alpha} \sum_{k=0}^\infty \frac{(-1)^k}{(2k+1)!}(c(t-c_0))^{2k+1} dt \\
&=& \frac{-c}{\gam{1-\alpha}} \sum_{k=0}^\infty \frac{(-1)^k}{\gam{2k+2}}c^{2k+1} \int_a^x(x-t)^{-\alpha} (t-c_0)^{2k+1} dt,
}

tomando el cambio de variable $t=c_0+(x-c_0)u$

\formula{
\ifrR{C}{a}{D}{x}{\alpha}f(t) 
&=& \frac{-c}{\gam{1-\alpha}}  \sum_{k=0}^\infty \frac{(-1)^k}{\gam{2k+2}}c^{2k+1} \left[ \int_\frac{a-c_0}{x-c_0}^1((x-c_0)-(x-c_0)u)^{-\alpha} ((x-c_0)u)^{2k+1} (x-c_0)du\right]\\
&=& -c  \sum_{k=0}^\infty \frac{(-1)^k}{\gam{1-\alpha}\gam{2k+2}}\left[ \int_0^1 \left( 1-u\right)^{-\alpha} u^{2k+1} du -\int_0^\frac{a-c_0}{x-c_0} \left( 1-u\right)^{-\alpha} u^{2k+1} du \right]c^{2k+1}  (x-c_0)^{2k-\alpha+2}\\
&=&- c  \sum_{k=0}^\infty \frac{(-1)^k}{\gam{1-\alpha}\gam{2k+2}}\left[ B(1-\alpha,2k+2)-B\left(\dfrac{a-c_0}{x-c_0};1-\alpha, 2k+2 \right) \right]c^{2k+1}  (x-c_0)^{2k-\alpha+2}\\
&=& -c  \sum_{k=0}^\infty \frac{(-1)^k}{\gam{1-\alpha}\gam{2k+2}} B(1-\alpha,2k+2) \left[ 1-\dfrac{ B\left(\dfrac{a-c_0}{x-c_0};1-\alpha, 2k+2 \right)}{B(1-\alpha,2k+2)} \right]c^{2k+1}  (x-c_0)^{2k+2-\alpha},
}

entonces

\begin{eqnarray}
\ifrR{C}{a}{D}{x}{\alpha}\co{c(x-c_0)}  = -c  \sum_{k=0}^\infty \frac{(-1)^k}{\gam{2k+3-\alpha}}  \left[ 1-\dfrac{ B\left(\dfrac{a-c_0}{x-c_0};1-\alpha, 2k+2 \right)}{B(1-\alpha,2k+2)} \right]c^{2k+1}  (x-c_0)^{2k+2-\alpha}.
\end{eqnarray}

tomado $c_0=a$ en la expresión anterior se obtiene

\begin{eqnarray}
\ifrR{C}{a}{D}{x}{\alpha}\co{c(x-a)}  &=& -c  \sum_{k=0}^\infty \frac{(-1)^k}{\gam{2k+3-\alpha}}  c^{2k+1}  (x-a)^{2k+2-\alpha}\nonumber\\
&=& -c^2 (x-a)^{2-\alpha} \sum_{k=0}^\infty \frac{(-1)^k}{\gam{2k+3-\alpha}}  c^{2k}  (x-a)^{2k}\nonumber\\
&=& -c^2 (x-a)^{2-\alpha} E_{2,3-\alpha}  \left(- c^2  (x-a)^2\right).
\end{eqnarray}

\subsubsection{Teorema de Taylor}

Para funciones suaves de una variable, el teorema de Taylor nos garantiza que

\begin{eqnarray}
f(x)=\sum_{k=0}^m \frac{f^{(k)}(c_0)}{k!}(x-c_0)^k + R_m(x,c_0),
\end{eqnarray}

donde $R_m(x,c_0)$ es el residuo

\begin{eqnarray}
R_m(x,c_0)&=&\int_{c}^x \frac{(x-t)^m}{m!}f^{(m+1)}(t)dt\nonumber \\
&=&\frac{1}{\gam{m+1}} \int_{c_0}^x (x-t)^m\der{d}{t}{m+1} f(t)dt\nonumber \\
&=& \ifr{c_0}{I}{x}{m+1}D^{m+1}f(t),
\end{eqnarray}

entonces

\formul{
f(x)=\sum_{k=0}^m \frac{f^{(k)}(c_0)}{k!}(x-c_0)^k + \ifr{c_0}{I}{x}{m+1}D^{m+1}f(t),
}{\label{taylor1}}

para $x$ próximo a $c_0$, el erro $R_n(x,c)$ es pequeño.


Tomando la derivada fraccionara de Riemann-Liouville de orden $\alpha$, con $\alpha\in(0,1)$, de la ecuación \eqref{taylor1} en el intervalo $(a,x)$ obtenemos

\formula{
\ifr{a}{D}{x}{\alpha} f(t)&=& \ifr{a}{D}{x}{\alpha} \sum_{k=0}^m \frac{f^{(k)}(c_0)}{k!}(t-c_0)^k +\ifr{a}{D}{x}{\alpha}  \ifr{c_0}{I}{x}{m+1}D^{m+1}f(t)\\
&=& \frac{1}{\gam{1-\alpha}}\frac{d}{dx}\left\{ \int_a^x (x-t)^{-\alpha}  \sum_{k=0}^m \frac{f^{(k)}(c_0)}{k!}(t-c_0)^k dt \right\} +D\ifr{a}{I}{x}{1-\alpha}  \ifr{c_0}{I}{x}{m+1}D^{m+1}f(t) \\
&=& \frac{d}{dx}\left\{ \sum_{k=0}^m \frac{f^{(k)}(c_0)}{\gam{1-\alpha}\gam{k+1}} \int_a^x (x-t)^{-\alpha}  (t-c_0)^k dt \right\} +D\ifr{a}{I}{x}{1-\alpha}  \ifr{c_0}{I}{x}{m+1}D^{m+1}f(t) ,
}

tomando el cambio de variable $t=c_0+(x-c_0)u$ en la integral de la suma 

\formula{
\ifr{a}{D}{x}{\alpha} f(t)&=& \frac{d}{dx}\left\{ \sum_{k=0}^m \frac{f^{(k)}(c_0)}{\gam{1-\alpha}\gam{k+1}} \int_\frac{a-c_0}{x-c_0}^1 ((x-c_0)-(x-c_0)u)^{-\alpha}  ((x-c_0)u)^k(x-c_0)du \right\} +D\ifr{a}{I}{x}{1-\alpha}  \ifr{c_0}{I}{x}{m+1}D^{m+1}f(t)\\
&=& \frac{d}{dx}\left\{ \sum_{k=0}^m \frac{f^{(k)}(c_0)}{\gam{1-\alpha}\gam{k+1}} \left[ \int_0^1 (1-u)^{-\alpha} u^k du -\int_0^\frac{a-c_0}{x-c_0} (1-u)^{-\alpha} u^k du\right](x-c_0)^{k-\alpha+1} \right\} +D\ifr{a}{I}{x}{1-\alpha}  \ifr{c_0}{I}{x}{m+1}D^{m+1}f(t)\\
&=& \frac{d}{dx}\left\{ \sum_{k=0}^m \frac{f^{(k)}(c_0)}{\gam{1-\alpha}\gam{k+1}} \left[ B\left(1-\alpha,k+1 \right) - B\left(\dfrac{a-c_0}{x-c_0};1-\alpha,k+1 \right)\right](x-c_0)^{k-\alpha+1} \right\} +D\ifr{a}{I}{x}{1-\alpha}  \ifr{c_0}{I}{x}{m+1}D^{m+1}f(t)\\
&=& \frac{d}{dx}\left\{ \sum_{k=0}^m \frac{f^{(k)}(c_0)}{\gam{1-\alpha}\gam{k+1}}B\left(1-\alpha,k+1 \right) \left[ 1 -\dfrac{ B\left(\dfrac{a-c_0}{x-c_0};1-\alpha,k+1 \right)}{B\left(1-\alpha,k+1 \right)}\right](x-c_0)^{k-\alpha+1} \right\} +D\ifr{a}{I}{x}{1-\alpha}  \ifr{c_0}{I}{x}{m+1}D^{m+1}f(t),
}

entonces

\begin{eqnarray}
\ifr{a}{D}{x}{\alpha} f(t)= \frac{d}{dx}\left\{ \sum_{k=0}^m \frac{f^{(k)}(c_0)}{\gam{k-\alpha+2}} \left[ 1 -\dfrac{ B\left(\dfrac{a-c_0}{x-c_0};1-\alpha,k+1 \right)}{B\left(1-\alpha,k+1 \right)}\right](x-c_0)^{k-\alpha+1} \right\} +D\ifr{a}{I}{x}{1-\alpha}  \ifr{c_0}{I}{x}{m+1}D^{m+1}f(t),
\end{eqnarray}

tomando $c_0=a$ en la expresión anterior y  la propiedad de semigrupo de la integral iterada 

\begin{eqnarray}
\ifr{a}{D}{x}{\alpha} f(t)= \frac{d}{dx}\left\{ \sum_{k=0}^m \frac{f^{(k)}(a)}{\gam{k-\alpha+2}} (x-a)^{k-\alpha+1} \right\} +D\ifr{a}{I}{x}{m-\alpha +2}  D^{m+1}f(t),
\end{eqnarray}

tomando el caso en que la función se expande en serie alrededor de $c_0=a$, tal que cuando $m\to \infty$ entonces $R_{m}(x,a)\to 0$ se obtiene

\begin{eqnarray}
\ifr{a}{D}{x}{\alpha} f(t)= \frac{d}{dx}\left\{ \sum_{k=0}^\infty \frac{f^{(k)}(a)}{\gam{k-\alpha+2}} (x-a)^{k-\alpha+1} \right\},
\end{eqnarray}


Tomando la derivada fraccionara de Caputo de orden $\alpha$, con $\alpha\in(0,1)$, de la ecuación \eqref{taylor1} en el intervalo $(a,x)$ obtenemos

\formula{
\ifrR{C}{a}{D}{x}{\alpha} f(t)&=& \ifrR{C}{a}{D}{x}{\alpha} \sum_{k=0}^m \frac{f^{(k)}(c_0)}{k!}(t-c_0)^k +\ifrR{C}{a}{D}{x}{\alpha}  \ifr{c_0}{I}{x}{m+1}D^{m+1}f(t)\\
&=& \frac{1}{\gam{1-\alpha}} \int_a^x (x-t)^{-\alpha}\frac{d}{dt} \left[ \sum_{k=0}^m \frac{f^{(k)}(c_0)}{k!}(t-c_0)^k\right]dt +\ifr{a}{I}{x}{1-\alpha}D  \ifr{c_0}{I}{x}{m+1}D^{m+1}f(t) \\
&=&   \sum_{k=1}^m \frac{f^{(k)}(c_0)}{\gam{1-\alpha}(k-1)!}  \left[ \int_a^x (x-t)^{-\alpha} (t-c_0)^{k-1}dt \right] +\ifr{a}{I}{x}{1-\alpha} \ifr{c_0}{I}{x}{m}D^{m+1}f(t),
}

tomando el cambio de variable $t=c_0+(x-c_0)u$ en la integral de la suma

\formula{
\ifrR{C}{a}{D}{x}{\alpha} f(t)&=& \sum_{k=0}^{m-1} \frac{f^{(k+1)}(c_0)}{\gam{1-\alpha}\gam{k+1}}  \left[ \int_\frac{a-c_0}{x-c_0}^1 ((x-c_0)-(x-c_0)u)^{-\alpha} ((x-c_0)u)^{k}(x-c_0)du \right] +\ifr{a}{I}{x}{1-\alpha} \ifr{c_0}{I}{x}{m}D^{m+1}f(t)\\
&=& \sum_{k=0}^{m-1} \frac{f^{(k+1)}(c_0)}{\gam{1-\alpha}\gam{k+1}}  \left[ \int_0^1 (1-u)^{-\alpha} u^{k}du - \int_0^\frac{a-c_0}{x-c_0} (1-u)^{-\alpha} u^{k}du\right](x-c_0)^{k-\alpha+1} +\ifr{a}{I}{x}{1-\alpha} \ifr{c_0}{I}{x}{m}D^{m+1}f(t)\\
&=& \sum_{k=0}^{m-1} \frac{f^{(k+1)}(c_0)}{\gam{1-\alpha}\gam{k+1}}  \left[ B\left(1-\alpha, k+1 \right) - B\left(\dfrac{a-c_0}{x-c_0};1-\alpha, k+1 \right)\right](x-c_0)^{k-\alpha+1} +\ifr{a}{I}{x}{1-\alpha} \ifr{c_0}{I}{x}{m}D^{m+1}f(t)\\
&=& \sum_{k=0}^{m-1} \frac{f^{(k+1)}(c_0)}{\gam{1-\alpha}\gam{k+1}} B\left(1-\alpha, k+1 \right) \left[ 1 - \dfrac{ B\left(\dfrac{a-c_0}{x-c_0};1-\alpha, k+1 \right)}{B\left(1-\alpha, k+1 \right)}\right](x-c_0)^{k-\alpha+1} +\ifr{a}{I}{x}{1-\alpha} \ifr{c_0}{I}{x}{m}D^{m+1}f(t),
}

entonces

\begin{eqnarray}
\ifrR{C}{a}{D}{x}{\alpha} f(t)=\sum_{k=0}^{m-1} \frac{f^{(k+1)}(c_0)}{\gam{k-\alpha+2}}  \left[ 1 - \dfrac{ B\left(\dfrac{a-c_0}{x-c_0};1-\alpha, k+1 \right)}{B\left(1-\alpha, k+1 \right)}\right](x-c_0)^{k-\alpha+1} +\ifr{a}{I}{x}{1-\alpha} \ifr{c_0}{I}{x}{m}D^{m+1}f(t),
\end{eqnarray}

tomando $c_0=a$ en la expresión anterior y  la propiedad de semigrupo de la integral iterada 

\begin{eqnarray}
\ifrR{C}{a}{D}{x}{\alpha} f(t)=\sum_{k=0}^{m-1} \frac{f^{(k+1)}(a)}{\gam{k-\alpha+2}}  (x-a)^{k-\alpha+1} + \ifr{a}{I}{x}{m-\alpha+1}D^{m+1}f(t),
\end{eqnarray}

tomando el caso en que la función se expande en serie alrededor de $c_0=a$, tal que cuando $m\to \infty$ entonces $R_{m}(x,a)\to 0$ se obtiene

\begin{eqnarray}
\ifrR{C}{a}{D}{x}{\alpha} f(t)=\sum_{k=0}^{\infty} \frac{f^{(k+1)}(a)}{\gam{k-\alpha+2}}  (x-a)^{k-\alpha+1}.
\end{eqnarray}

\subsubsection{Transformada de Laplace para los Operadores Diferointegrables}

Para una función $F(t)$ definida para $t>0$, la transformada de Laplace de $F(t)$ se define como

\formula{
\mathcal{L}\lbrace F(t) \rbrace=f(s) =\int_0^\infty \ex{-st}F(t)dt,
}

tomando la derivada de $F(t)$ y aplicando la transformada de Laplace se obtiene

\formula{
\mathcal{L}\left\lbrace \frac{d}{dt}F(t) \right\rbrace&=& \int_0^\infty \ex{-st}\frac{d}{dt}F(t)dt,
}

integrando por partes tomando

\formula{
\begin{array}{lcl}
u=\ex{-st} &\rightarrow & du=-s\ex{-st}dt,\\
dv=\dfrac{d}{dt}F(t)dt &\rightarrow & v=F(t),
\end{array}
}

entonces

\formula{
\mathcal{L}\left\lbrace \frac{d}{dt}F(t) \right\rbrace&=& \ex{-st}F(t)\Big{|}_0^\infty-\int_0^\infty -s\ex{-st}F(t)dt\\
&=& s\int_0^\infty \ex{-st}F(t)dt -F(0)\\
&=& s\mathcal{L}\left\{ F(t)\right\}-F(0),\\
}

por lo tanto

\formul{
\mathcal{L}\left\{ \dfrac{d}{dt}F(t) \right\} =s\mathcal{L}\left\{ F(t)\right\}-F(0),
}{\label{lap1}}

tomando la segunda derivada de $F(t)$, aplicando la transformada de Laplace y ocupando \eqref{lap1} se obtiene

\formula{
\mathcal{L}\left\{ \der{d}{t}{2}F(t) \right\} &=&s\mathcal{L}\left\{ \dfrac{d}{dt}F(t)\right\}-\dfrac{d}{dt}F(0)\\
&=& s\left[s\mathcal{L}\left\{ F(t)\right\}-F(0) \right]-\dfrac{d}{dt}F(0)\\
&=&s^2\mathcal{L}\left\{ F(t)\right\}-sF(0)-\dfrac{d}{dt}F(0)\\
&=& s^2\mathcal{L}\left\{F(t) \right\}-\sum_{k=0}^1 s^{1-k}\der{d}{t}{k}F(0)\\
&=& s^2\mathcal{L}\left\{F(t) \right\}-\sum_{k=0}^1 s^{k}\der{d}{t}{1-k}F(0),
}

tomando la tercera derivada de $F(t)$, aplicando la transformada de Laplace y ocupando \eqref{lap1} se obtiene

\formula{
\mathcal{L}\left\{ \der{d}{t}{3}F(t) \right\} &=&s\mathcal{L}\left\{ \der{d}{t}{2}F(t)\right\}-\der{d}{t}{2}F(0)\\
&=& s\left[s\mathcal{L}\left\{ \dfrac{d}{dt} F(t)\right\}-\dfrac{d}{dt}F(0) \right]-\der{d}{t}{2}F(0)\\
&=&s^2\mathcal{L}\left\{ \dfrac{d}{dt} F(t)\right\}-s\dfrac{d}{dt}F(0)-\der{d}{t}{2}F(0)\\
&=&s^2\left[s\mathcal{L}\left\{ F(t)\right\}-F(0)\right]-s\dfrac{d}{dt}F(0)-\der{d}{t}{2}F(0)\\
&=&s^3\mathcal{L}\left\{ F(t)\right\}-s^2F(0)-s\dfrac{d}{dt}F(0)-\der{d}{t}{2}F(0)\\
&=& s^3\mathcal{L}\left\{F(t) \right\}-\sum_{k=0}^2 s^{2-k}\der{d}{t}{k}F(0)\\
&=& s^3\mathcal{L}\left\{F(t) \right\}-\sum_{k=0}^2 s^{k}\der{d}{t}{2-k}F(0),
}

de los pasos anteriores podemos intuir una formula para la transformada de Laplace de la derivada $n$-ésima de la función $F(t)$

\formul{
\mathcal{L}\left\{ \der{d}{t}{n}F(t) \right\}&=& s^n\mathcal{L}\left\{F(t) \right\}-\sum_{k=0}^{n-1} s^{(n-1)-k}\der{d}{t}{k}F(0)\nonumber\\
&=&s^n\mathcal{L}\left\{F(t) \right\}-\sum_{k=0}^{n-1} s^{k}\der{d}{t}{(n-1)-k}F(0),
}{\label{lap2}}

procedemos a demostrar la ecuación \eqref{lap2} por inducción, los casos para $n=1$ y $n=2$ se obtuvieron durante la construcción, supongamos que la formula es valida para $n=k$ (con $k<n$), sea $n=k+1$ y ocupando la ecuación \eqref{lap1} se obtiene

\formula{
\mathcal{L}\left\{ \der{d}{t}{k+1}F(t) \right\} =s\mathcal{L}\left\{ \der{d}{t}{k}F(t)\right\}-\der{d}{t}{k}F(0),
}

ocupando la hipótesis de inducción se obtiene

\formula{
\mathcal{L}\left\{ \der{d}{t}{k+1}F(t) \right\} &=&s\left[ s^k\mathcal{L}\left\{F(t) \right\}-\sum_{m=0}^{k-1} s^{(k-1)-m}\der{d}{t}{m}F(0)\right]-\der{d}{t}{k}F(0)\\
&=& s^{k+1}\mathcal{L}\left\{F(t) \right\}-\sum_{m=0}^{k-1} s^{k-m}\der{d}{t}{m}F(0)-\der{d}{t}{k}F(0)\\
&=& s^{k+1}\mathcal{L}\left\{F(t) \right\}-\sum_{m=0}^{k} s^{k-m}\der{d}{t}{m}F(0)\\
&=& s^{k+1}\mathcal{L}\left\{F(t) \right\}-\sum_{m=0}^{k} s^{m}\der{d}{t}{k-m}F(0),
}

con lo cual se demuestra la validez de la ecuación \eqref{lap2} para todo $n\in \mathds{N}$.

Analicemos ahora el comportamiento de la transformada de Laplace para la integral de $F(t)$ en el intervalo $(a,t)$,  sea

\formula{
G(t)=\int_a^t F(u)du,
}

ocupando \eqref{lap1} se obtiene

\formula{
\mathcal{L}\left\{ \dfrac{d}{dt}G(t) \right\} = s\mathcal{L}\left\{ G(t)\right\}- G(0),
}

entonces

\formul{
\mathcal{L}\left\{ \int_a^t F(u)du\right\} = \dfrac{1}{s}\mathcal{L}\left\{\dfrac{d}{dt}\int_a^t F(u)du\right\} + \dfrac{1}{s} \left.\int_a^t F(u)du\right|_{t=0} ,
}{\label{lap3}}

aplicando la ecuación \eqref{lap3} a la ecuación \eqref{int4} para $n=2$ obtenemos la transformada de Laplace para $2$-ésima integral de una función $F(t)$

\formula{
\mathcal{L}\left\{ \ifr{a}{I}{t}{2}F(u)\right\} &=&\mathcal{L}\left\{ \dfrac{1}{(2-1)!} \int_a^t (t-u)^{2-1}F(u)du\right\}\\
&=&\mathcal{L}\left\{ \int_a^t (t-u)F(u)du\right\}\\
&=&\dfrac{1}{s}\mathcal{L}\left\{\dfrac{d}{dt}\int_a^t (t-u)F(u)du\right\} +\dfrac{1}{s} \left. \int_a^t (t-u)F(u)du\right|_{t=0} \\
&=&\dfrac{1}{s}\mathcal{L}\left\{\int_a^t \dfrac{\partial}{\partial t}(t-u)F(u)du + \left.\left[(t-u)F(u)\frac{d}{dt}u\right]\right|_a^t\right\} +\dfrac{1}{s} \left. \int_a^t (t-u)F(u)du\right|_{t=0} \\
&=&\dfrac{1}{s}\mathcal{L}\left\{\int_a^t F(u)du \right\} +\dfrac{1}{s} \left. \int_a^t (t-u)F(u)du\right|_{t=0} \\
&=&\dfrac{1}{s^2}\mathcal{L}\left\{\dfrac{d}{dt}\int_a^t F(u)du\right\} + \dfrac{1}{s^2}\left. \int_a^t F(u)du\right|_{t=0} +\dfrac{1}{s} \left. \int_a^t (t-u)F(u)du\right|_{t=0}\\
&=& \dfrac{1}{s^2}\mathcal{L}\left\{\dfrac{d}{dt}\int_a^t F(u)du\right\} + \sum_{k=0}^1 \frac{1}{s^{2-k}}\left. \int_a^t (t-u)^{k}F(u)du\right|_{t=0}\\
&=& \dfrac{1}{s^2}\mathcal{L}\left\{\dfrac{d}{dt}\int_a^t F(u)du\right\} + \sum_{k=1}^2 \frac{1}{s^{k}}\left. \int_a^t (t-u)^{2-k}F(u)du\right|_{t=0},
}

aplicando la ecuación \eqref{lap3} a la ecuación \eqref{int4} para $n=3$ obtenemos la transformada de Laplace para $3$-ésima integral de una función $F(t)$

\formula{
\mathcal{L}\left\{ \ifr{a}{I}{t}{3}F(u)\right\} &=&\mathcal{L}\left\{ \dfrac{1}{(3-1)!} \int_a^t (t-u)^{3-1}F(u)du\right\}\\
&=&\dfrac{1}{2!}\mathcal{L}\left\{ \int_a^t (t-u)^2F(u)du\right\}\\
&=&\dfrac{1}{2!} \left[\dfrac{1}{s}\mathcal{L}\left\{\dfrac{d}{dt}\int_a^t (t-u)^2F(u)du\right\} +\dfrac{1}{s}\left. \int_a^t (t-u)^2F(u)du\right|_{t=0} \right] \\
&=&\dfrac{1}{2!s}\mathcal{L}\left\{\int_a^t \dfrac{\partial}{\partial t}(t-u)^2F(u)du + \left.\left[(t-u)^2F(u)\frac{d}{dt}u\right]\right|_a^t\right\} +\dfrac{1}{2!s}\left. \int_a^t (t-u)^2F(u)du \right|_{t=0} \\
&=&\dfrac{1}{2!s}\mathcal{L}\left\{\int_a^t 2(t-u)F(u)du \right\} +\dfrac{1}{s} \ifr{a}{I}{t}{3}F(0) \\
&=&\dfrac{1}{s}\mathcal{L}\left\{\int_a^t (t-u)F(u)du \right\} +\dfrac{1}{s} \ifr{a}{I}{t}{3}F(0) \\
&=&\dfrac{1}{s}\mathcal{L}\left\{\ifr{a}{I}{t}{2}F(u) \right\} +\dfrac{1}{s} \ifr{a}{I}{t}{3}F(0) \\
&=&\dfrac{1}{s}\left[ \dfrac{1}{s^2}\mathcal{L}\left\{\dfrac{d}{dt}\int_a^t F(u)du\right\} + \sum_{k=0}^1 \frac{1}{s^{2-k}}\left.\int_a^0 (t-u)^{k}F(u)du\right|_{t=0} \right] +\dfrac{1}{s} \ifr{a}{I}{t}{3}F(0)\\
&=&  \dfrac{1}{s^3}\mathcal{L}\left\{\dfrac{d}{dt}\int_a^t F(u)du\right\} + \sum_{k=0}^1 \frac{1}{ s^{3-k}} \ifr{a}{I}{t}{k+1}F(0)+\dfrac{1}{s} \ifr{a}{I}{t}{3}F(0)\\
&=& \dfrac{1}{s^3}\mathcal{L}\left\{\dfrac{d}{dt}\int_a^t F(u)du\right\} + \sum_{k=0}^2 \frac{1}{ s^{3-k}} \ifr{a}{I}{t}{k+1}F(0)\\
&=& \dfrac{1}{s^3}\mathcal{L}\left\{\dfrac{d}{dt}\int_a^t F(u)du\right\} + \sum_{k=1}^3 \frac{1}{ s^{k}} \ifr{a}{I}{t}{4-k}F(0),
}

aplicando la ecuación \eqref{lap3} a la ecuación \eqref{int4} para $n=4$ obtenemos la transformada de Laplace para $4$-ésima integral de una función $F(t)$

\formula{
\mathcal{L}\left\{ \ifr{a}{I}{t}{4}F(u)\right\} &=&\mathcal{L}\left\{ \dfrac{1}{(4-1)!} \int_a^t (t-u)^{4-1}F(u)du\right\}\\
&=&\dfrac{1}{3!}\mathcal{L}\left\{ \int_a^t (t-u)^3F(u)du\right\}\\
&=&\dfrac{1}{3!} \left[\dfrac{1}{s}\mathcal{L}\left\{\dfrac{d}{dt}\int_a^t (t-u)^3F(u)du\right\} +\dfrac{1}{s}\left. \int_a^t(t-u)^3F(u)du\right|_{t=0} \right] \\
&=&\dfrac{1}{3!s}\mathcal{L}\left\{\int_a^t \dfrac{\partial}{\partial t}(t-u)^3F(u)du + \left.\left[(t-u)^3F(u)\frac{d}{dt}u\right]\right|_a^t\right\} +\dfrac{1}{3!s}\left. \int_a^t (t-u)^3F(u)du\right|_{t=0} \\
&=&\dfrac{1}{3!s}\mathcal{L}\left\{\int_a^t 3(t-u)^2F(u)du \right\} +\dfrac{1}{s} \ifr{a}{I}{t}{4}F(0) \\
&=&\dfrac{1}{2!s}\mathcal{L}\left\{\int_a^t (t-u)^2F(u)du \right\} +\dfrac{1}{s} \ifr{a}{I}{t}{4}F(0)  \\
&=& \dfrac{1}{s}\left[\dfrac{1}{2!}\mathcal{L}\left\{\int_a^t (t-u)^2F(u)du \right\}\right] +\dfrac{1}{s} \ifr{a}{I}{t}{4}F(0) \\
&=& \dfrac{1}{s}\mathcal{L}\left\{\ifr{a}{I}{t}{3}F(u) \right\} +\dfrac{1}{s} \ifr{a}{I}{t}{4}F(0)  \\
&=& \dfrac{1}{s}\left[ \dfrac{1}{s^3}\mathcal{L}\left\{\dfrac{d}{dt}\int_a^t F(u)du\right\} + \sum_{k=0}^2 \frac{1}{ s^{3-k}} \ifr{a}{I}{t}{k+1}F(0) \right] +\dfrac{1}{s} \ifr{a}{I}{t}{4}F(0)  \\
&=&  \dfrac{1}{s^4}\mathcal{L}\left\{\dfrac{d}{dt}\int_a^t F(u)du\right\} + \sum_{k=0}^2 \frac{1}{ s^{4-k}} \ifr{a}{I}{t}{k+1}F(0) +\dfrac{1}{s} \ifr{a}{I}{t}{4}F(0)  \\
&=&  \dfrac{1}{s^4}\mathcal{L}\left\{\dfrac{d}{dt}\int_a^t F(u)du\right\} + \sum_{k=0}^3 \frac{1}{ s^{4-k}} \ifr{a}{I}{t}{k+1}F(0)\\
&=&  \dfrac{1}{s^4}\mathcal{L}\left\{\dfrac{d}{dt}\int_a^t F(u)du\right\} + \sum_{k=1}^4 \frac{1}{ s^{k}} \ifr{a}{I}{t}{5-k}F(0),
}

del procedimiento anterior se puede deducir una formula para la transformada de Laplace de la $n$-ésima integral de $F(t)$

\formul{
\mathcal{L}\left\{ \ifr{a}{I}{t}{n}F(u) \right\}&=& \frac{1}{(n-1)!}\mathcal{L}\left\{\int_a^t (t-u)^{n-1}F(u)du \right\} \nonumber \\
&=& \dfrac{1}{s^n}\mathcal{L}\left\{\dfrac{d}{dt}\int_a^t F(u)du\right\} + \sum_{k=0}^{n-1} \frac{1}{ s^{n-k}}\ifr{a}{I}{t}{k+1} F(0) 
\nonumber \\
&=& \dfrac{1}{s^n}\mathcal{L}\left\{\dfrac{d}{dt}\int_a^t F(u)du\right\} + \sum_{k=1}^{n} \frac{1}{ s^{k}}\ifr{a}{I}{t}{(n+1)-k} F(0), 
}{\label{lap4}}

procedemos a demostrar la ecuación \eqref{lap4} por medio de el proceso de inducción, los casos para $n=1$ y $n=2$ se obtuvieron durante la construcción de la formula anterior, suponemos que la formula es cierta para $n=k$ (con $k<n$), sea $n=k+1$, entonces 

\formula{
\mathcal{L}\left\{ \ifr{a}{I}{t}{k+1}F(u)\right\} &=&\mathcal{L}\left\{ \dfrac{1}{((k+1)-1)!} \int_a^t (t-u)^{(k+1)-1}F(u)du\right\}\\
&=&\dfrac{1}{k!}\mathcal{L}\left\{ \int_a^t (t-u)^kF(u)du\right\}\\
&=&\dfrac{1}{k!} \left[\dfrac{1}{s}\mathcal{L}\left\{\dfrac{d}{dt}\int_a^t (t-u)^kF(u)du\right\} +\dfrac{1}{s}\left. \int_a^t (t-u)^kF(u)du\right|_{t=0} \right] \\
&=&\dfrac{1}{k!s}\mathcal{L}\left\{\int_a^t \dfrac{\partial}{\partial t}(t-u)^kF(u)du + \left.\left[(t-u)^kF(u)\frac{d}{dt}u\right]\right|_a^t\right\} +\dfrac{1}{k!s}\left. \int_a^t (t-u)^kF(u)du\right|_{t=0} \\
&=&\dfrac{1}{k!s}\mathcal{L}\left\{\int_a^t k(t-u)^{k-1}F(u)du \right\} +\dfrac{1}{s} \ifr{a}{I}{t}{k+1}F(0) \\
&=&\dfrac{1}{(k-1)!s}\mathcal{L}\left\{\int_a^t (t-u)^{k-1}F(u)du \right\} +\dfrac{1}{s} \ifr{a}{I}{t}{k+1}F(0) \\
&=& \dfrac{1}{s}\left[\dfrac{1}{(k-1)!}\mathcal{L}\left\{\int_a^t (t-u)^{k-1}F(u)du \right\}\right] +\dfrac{1}{s} \ifr{a}{I}{t}{k+1}F(0) \\
&=& \dfrac{1}{s}\mathcal{L}\left\{\ifr{a}{I}{t}{k}F(u)\right\} +\dfrac{1}{s} \ifr{a}{I}{t}{k+1}F(0) \\
&=& \dfrac{1}{s}\left[\dfrac{1}{s^k}\mathcal{L}\left\{\dfrac{d}{dt}\int_a^t F(u)du\right\} + \sum_{m=0}^{k-1} \frac{1}{ s^{k-m}}\ifr{a}{I}{t}{m+1} F(0)\right] +\dfrac{1}{s} \ifr{a}{I}{t}{k+1}F(0)\\
&=& \dfrac{1}{s^{k+1}}\mathcal{L}\left\{\dfrac{d}{dt}\int_a^t F(u)du\right\} + \sum_{m=0}^{k-1} \frac{1}{ s^{(k+1)-m}}\ifr{a}{I}{t}{m+1} F(0) +\dfrac{1}{s} \ifr{a}{I}{t}{k+1}F(0)\\
&=& \dfrac{1}{s^{k+1}}\mathcal{L}\left\{\dfrac{d}{dt}\int_a^t F(u)du\right\} + \sum_{m=0}^{k} \frac{1}{ s^{(k+1)-m}}\ifr{a}{I}{t}{m+1} F(0)\\
&=& \dfrac{1}{s^{k+1}}\mathcal{L}\left\{\dfrac{d}{dt}\int_a^t F(u)du\right\} + \sum_{m=1}^{k+1} \frac{1}{ s^{m}}\ifr{a}{I}{t}{((k+1)+1)-m} F(0),
}

con lo cual se demuestra la validez de la ecuación \eqref{lap4} para todo $n\in \mathds{N}$.


Una vez obtenidas las formulas de las transformadas de Laplace de las derivadas e integrales $n$-ésimas procedemos a calcular las transformadas de los operadores diferointegrables, comenzamos tomando la ecuación \eqref{int9}

\formula{
\ifr{a}{D}{t}{\alpha}F(u) 
&=& D^n\ifr{a}{I}{t}{n-\alpha}F(u)\\
&=&\frac{1}{\gam{n-\alpha}} \der{d}{t}{n}  \int_a^t(t-u)^{n-\alpha-1}  F(u)du, 
}

aplicando la transformada de Laplace y utilizando la ecuación \eqref{lap2} se obtiene

\formula{
\mathcal{L} \left\{  \ifr{a}{D}{t}{\alpha}F(u) \right\} &=& \mathcal{L}\left\{ D^n\ifr{a}{I}{t}{n-\alpha}F(u)\right\}\\
&=&s^n\mathcal{L}\left\{\ifr{a}{I}{t}{n-\alpha}F(u) \right\}-\sum_{k=0}^{n-1} s^{(n-1)-k}\der{d}{t}{k}\ifr{a}{I}{t}{n-\alpha}F(0)\\
&=&s^n\mathcal{L}\left\{\ifr{a}{I}{t}{n-\alpha}F(u) \right\}-\sum_{k=0}^{n-1} s^{k}\der{d}{t}{(n-1)-k}\ifr{a}{I}{t}{n-\alpha}F(0),
}

utilizando la ecuación \eqref{lap4} y tomando la función techo en el extremo superior de la suma se obtiene

\formul{
\mathcal{L} \left\{  \ifr{a}{D}{t}{\alpha}F(u) \right\} &=& s^n\mathcal{L}\left\{\ifr{a}{I}{t}{n-\alpha}F(u) \right\}-\sum_{k=0}^{n-1} s^{k}\der{d}{t}{(n-1)-k}\ifr{a}{I}{t}{n-\alpha}F(0) \nonumber\\
&=& s^n \left[\dfrac{1}{s^{n-\alpha}}\mathcal{L}\left\{\dfrac{d}{dt}\int_a^t F(u)du\right\} + \sum_{k=0}^{\lceil(n-\alpha)-1\rceil} \frac{1}{ s^{(n-\alpha)-k}}\ifr{a}{I}{t}{k+1} F(0) \right]-\sum_{k=0}^{n-1} s^{k}\der{d}{t}{(n-1)-k}\ifr{a}{D}{t}{\alpha-n}F(0) \nonumber\\
&=& s^n \left[\dfrac{1}{s^{n-\alpha}}\mathcal{L}\left\{\dfrac{d}{dt}\int_a^t F(u)du\right\} + \sum_{k=1}^{\lceil n-\alpha\rceil} \frac{1}{ s^{k}}\ifr{a}{I}{t}{(n-\alpha+1)-k} F(0) \right]-\sum_{k=0}^{n-1} s^{k}\ifr{a}{D}{t}{(\alpha-1)-k}F(0) \nonumber\\
&=& s^{\alpha}\mathcal{L}\left\{\dfrac{d}{dt}\int_a^t F(u)du\right\} + \sum_{k=1}^{\lceil n-\alpha\rceil}  s^{n-k}\ifr{a}{I}{t}{(n-\alpha+1)-k} F(0)-\sum_{k=0}^{n-1} s^{k}\ifr{a}{D}{t}{(\alpha-1)-k}F(0) \nonumber\\
&=& s^\alpha\mathcal{L}\left\{F(t) \right\}-\sum_{k=0}^{n-1} s^{k}\ifr{a}{D}{t}{(\alpha-1)-k}F(0)+ \sum_{k=1}^{\lceil n-\alpha\rceil}  s^{n-k}\ifr{a}{I}{t}{(n-\alpha+1)-k} F(0)
, \label{lap5}
}{}


tomando ahora la ecuación \eqref{int10}

\formula{
\ifrR{C}{a}{D}{t}{\alpha}F(u) 
&=& \ifr{a}{I}{t}{n-\alpha}D^nF(u)\\
&=&\frac{1}{\gam{n-\alpha}}  \int_a^t(t-u)^{n-\alpha-1} \der{d}{u}{n} f(u)du, 
}

aplicando la transformada de Laplace y utilizando la ecuación \eqref{lap4} y tomando la función techo en el extremo superior de la suma se obtiene

\formula{
\mathcal{L} \left\{  \ifrR{C}{a}{D}{t}{\alpha}F(u) \right\} &=& \mathcal{L}\left\{ \ifr{a}{I}{t}{n-\alpha}D^nF(u)\right\}\\
&=& \dfrac{1}{s^{n-\alpha}}\mathcal{L}\left\{\dfrac{d}{dt}\int_a^t D^nF(u)du\right\} + \sum_{k=0}^{\lceil (n-\alpha)-1 \rceil } \frac{1}{ s^{(n-\alpha)-k}}\ifr{a}{I}{t}{k+1} D^nF(0)\\
&=& \dfrac{1}{s^{n-\alpha}}\mathcal{L}\left\{\der{d}{t}{n} F(t)\right\} + \sum_{k=1}^{\lceil n-\alpha\rceil } \frac{1}{ s^{k}}\ifr{a}{I}{t}{(n-\alpha+1)-k}D^n F(0),
}

utilizando la ecuación \eqref{lap2} en la expresión anterior obtenemos

\formul{
\mathcal{L} \left\{  \ifrR{C}{a}{D}{t}{\alpha}F(u) \right\} &=& \dfrac{1}{s^{n-\alpha}}\mathcal{L}\left\{\der{d}{t}{n} F(t)\right\} + \sum_{k=1}^{\lceil n-\alpha\rceil } \frac{1}{ s^{k}}\ifr{a}{I}{t}{(n-\alpha+1)-k} D^nF(0) \nonumber\\
&=&\dfrac{1}{s^{n-\alpha}} \left[ s^n\mathcal{L}\left\{F(t) \right\}-\sum_{k=0}^{n-1} s^{(n-1)-k}\der{d}{t}{k}F(0)\right]+ \sum_{k=1}^{\lceil n-\alpha\rceil } \frac{1}{ s^{k}}\ifr{a}{I}{t}{(1-\alpha)-k} F(0) \nonumber \\
&=& s^\alpha\mathcal{L}\left\{F(t) \right\}-\sum_{k=0}^{n-1} s^{(\alpha-1)-k}\der{d}{t}{k}F(0)+ \sum_{k=1}^{\lceil n-\alpha\rceil} \frac{1}{ s^{k}}\ifr{a}{I}{t}{(1-\alpha)-k} F(0), \label{lap6}
}{}

\subsubsection{Transformada de Laplace a través de la Convolución}

A menudo al tratar de obtener la transformada de Laplace de una función $G(t)$ se obtiene que la transformada es el producto de las transformadas de dos funciones, entiéndase por

\formula{
\mathcal{L}\left\{ G(t)\right\}=\mathcal{L}\left\{ F_1(t)\right\}\mathcal{L}\left\{ F_2(t)\right\},
}

esto nos lleva a deducir que existe una relación entre la función $G(t)$ y el producto de las funciones $F_1(t)$ y $F_2(t)$, esta relación se da a través de un producto generalizado conocido como convolución la cual se define por

\formula{
F_1(t)*F_2(t)=\int_0^t F_1(t-u)F_2(u)du,
}{\label{con1}}

para comprobar que la ecuación \eqref{con1} satisface lo pedido comencemos analizando el producto de las transformadas de las funciones $F_1(t)$ y $F_2(t)$

\formula{
\mathcal{L}\left\{ F_1(t)\right\}\mathcal{L}\left\{ F_2(t)\right\}&=& \int_0^\infty e^{-st}F_1(t)dt \int_0^\infty e^{-st}F_2(t)dt\\
&=& \int_0^\infty e^{-sx}F_1(x)dx \int_0^\infty e^{-sy}F_2(y)dy\\
&=& \int_0^\infty \int_0^\infty e^{-s(x+y)}F_1(x) F_2(y)dxdy,
}

tomando el cambio de variable

\formula{
u&=&y,\\
t&=& x+y = x+u,
}

entonces

\formula{
\begin{matrix}
x>0&\rightarrow &t -u>0 &\rightarrow&t>u,\\
y>0& \rightarrow & u>0,
\end{matrix} 
}

además

\formula{
dxdy&=& \frac{\partial(x,y)}{\partial(t,u)}dtdu\\
&=& \left|\begin{matrix}
\dfrac{\partial x}{\partial t} & \dfrac{\partial y}{\partial t}\\
\dfrac{\partial x}{\partial u} & \dfrac{\partial y}{\partial u}
\end{matrix}\right| dtdu\\
&=& \left|\begin{matrix}
1 & 0\\
-1 & 1
\end{matrix}\right| dtdu,\\
&=& dtdu,
}

lo que implica

\formula{
\mathcal{L}\left\{ F_1(t)\right\}\mathcal{L}\left\{ F_2(t)\right\}
&=& \int_0^\infty \int_0^\infty e^{-s(x+y)}F_1(x) F_2(y)dxdy\\
&=& \int_0^\infty \int_0^t e^{-st}F_1(t-u) F_2(u)dtdu\\
&=& \int_0^\infty e^{-st} \left\{\int_0^t F_1(t-u) F_2(u) du\right\} dt\\
&=& \mathcal{L} \left\{\int_0^t F_1(t-u) F_2(u) du\right\}, 
}

entonces obtenemos que la transformada de Laplace de la convolución de dos funciones es el producto de las transformaciones de cada una de las funciones

\formul{
\mathcal{L}\left\{ F_1(t)*F_2(t) \right\}&=& \mathcal{L} \left\{\int_0^t F_1(t-u) F_2(u) du\right\}=\mathcal{L}\left\{ F_1(t)\right\}\mathcal{L}\left\{ F_2(t)\right\},
}{\label{con2}}

un resultado que nos sera de utilidad mas adelante es la transformada de Laplace de el monomio $t^n$ con $n>0$

\formul{
\mathcal{L}\left\{ t^n\right\}&=& \int_0^\infty e^{-st}t^n dt \nonumber\\
&=& \int_0^\infty (-1)^n\der{\partial}{s}{n} e^{-st} dt \nonumber \\
&=& (-1)^n \der{\partial}{s}{n} \int_0^\infty  e^{-st} dt \nonumber\\
&=& (-1)^n \der{\partial}{s}{n} s^{-1} \nonumber\\
&=& (-1)^n \frac{\gam{0}}{\gam{-n}} s^{-1-n}\nonumber\\
&=& (-1)^n (-1)^n\frac{\gam{n+1}}{\gam{1}} s^{-1-n}\nonumber\\
&=& \frac{\gam{n+1}}{s^{n+1}}.
}{\label{con3}}


Una vez obtenida la transformada de Laplace de la convolución  procedemos a calcular las transformadas de los operadores diferointegrables, comenzamos tomando la ecuación \eqref{int9} en el intervalo $(0,t)$

\formula{
\ifr{0}{D}{t}{\alpha}F(u) 
&=& D^n\ifr{0}{I}{t}{n-\alpha}F(u)\\
&=&\frac{1}{\gam{n-\alpha}} \der{d}{t}{n}  \int_0^t(t-u)^{n-\alpha-1}  F(u)du, 
}

aplicando la transformada de Laplace y utilizando la ecuación \eqref{lap2} se obtiene

\formula{
\mathcal{L} \left\{  \ifr{0}{D}{t}{\alpha}F(u) \right\} &=& \mathcal{L}\left\{ D^n\ifr{0}{I}{t}{n-\alpha}F(u)\right\}\\
&=&s^n\mathcal{L}\left\{\ifr{0}{I}{t}{n-\alpha}F(u) \right\}-\sum_{k=0}^{n-1} s^{(n-1)-k}\der{d}{t}{k}\ifr{0}{I}{t}{n-\alpha}F(0)\\
&=&s^n\mathcal{L}\left\{\ifr{0}{I}{t}{n-\alpha}F(u) \right\}-\sum_{k=0}^{n-1} s^{k}\der{d}{t}{(n-1)-k}\ifr{0}{I}{t}{n-\alpha}F(0),
}

utilizando la ecuaciones \eqref{con2} y \eqref{con3} se obtiene

\formul{
\mathcal{L} \left\{  \ifr{0}{D}{t}{\alpha}F(u) \right\} &=& s^n\mathcal{L}\left\{\ifr{a}{I}{t}{n-\alpha}F(u) \right\}-\sum_{k=0}^{n-1} s^{k}\der{d}{t}{(n-1)-k}\ifr{0}{I}{t}{n-\alpha}F(0) \nonumber\\
&=& \frac{s^n}{\gam{n-\alpha}} \mathcal{L}\left\{\int_0^t(t-u)^{n-\alpha-1} F(u)du\right\} -\sum_{k=0}^{n-1} s^{k}\der{d}{t}{(n-1)-k}\ifr{0}{D}{t}{\alpha-n}F(0) \nonumber\\
&=& \frac{s^n}{\gam{n-\alpha}} \mathcal{L}\left\{t^{n-\alpha-1}\right\}\mathcal{L}\left\{ F(t)\right\} -\sum_{k=0}^{n-1} s^{k}\ifr{0}{D}{t}{(\alpha-1)-k}F(0) \nonumber\\
&=& \frac{s^n}{\gam{n-\alpha}}\frac{\gam{n-\alpha}}{s^{n-\alpha}}\mathcal{L}\left\{ F(t)\right\} -\sum_{k=0}^{n-1} s^{k}\ifr{0}{D}{t}{(\alpha-1)-k}F(0) \nonumber\\
&=& s^\alpha\mathcal{L}\left\{ F(t)\right\} -\sum_{k=0}^{n-1} s^{k}\ifr{0}{D}{t}{(\alpha-1)-k}F(0)
, \label{con4}
}{}


tomando ahora la ecuación \eqref{int10} en el intervalo $(0,t)$

\formula{
\ifrR{C}{0}{D}{t}{\alpha}F(u) 
&=& \ifr{0}{I}{t}{n-\alpha}D^nF(u)\\
&=&\frac{1}{\gam{n-\alpha}}  \int_0^t(t-u)^{n-\alpha-1} \der{d}{u}{n} F(u)du, 
}

aplicando la transformada de Laplace y utilizando la ecuaciones \eqref{con2} y \eqref{con3} se obtiene

\formula{
\mathcal{L} \left\{  \ifrR{C}{0}{D}{t}{\alpha}F(u) \right\} &=& \mathcal{L}\left\{ \frac{1}{\gam{n-\alpha}}  \int_0^t(t-u)^{n-\alpha-1} \der{d}{u}{n} F(u)du\right\}\\
&=&\frac{1}{\gam{n-\alpha}} \mathcal{L}\left\{   \int_0^t(t-u)^{n-\alpha-1} \der{d}{u}{n} F(u)du\right\}\\
&=&\frac{1}{\gam{n-\alpha}} \mathcal{L}\left\{   t^{n-\alpha-1} \right\}\mathcal{L}\left\{    \der{d}{t}{n} F(t)\right\}\\
&=&\frac{1}{\gam{n-\alpha}} \frac{\gam{n-\alpha}}{s^{n-\alpha}} \mathcal{L}\left\{    \der{d}{t}{n} F(t)\right\}\\
&=& \frac{1}{s^{n-\alpha}} \mathcal{L}\left\{    \der{d}{t}{n} F(t)\right\},
}

utilizando la ecuación \eqref{lap2} en la expresión anterior obtenemos

\formul{
\mathcal{L} \left\{  \ifrR{C}{0}{D}{t}{\alpha}F(u) \right\} &=& \dfrac{1}{s^{n-\alpha}}\mathcal{L}\left\{\der{d}{t}{n} F(t)\right\} \nonumber\\
&=&\dfrac{1}{s^{n-\alpha}} \left[ s^n\mathcal{L}\left\{F(t) \right\}-\sum_{k=0}^{n-1} s^{(n-1)-k}\der{d}{t}{k}F(0)\right] \nonumber \\
&=& s^\alpha\mathcal{L}\left\{F(t) \right\}-\sum_{k=0}^{n-1} s^{(\alpha-1)-k}\der{d}{t}{k}F(0). \label{con5}
}{}


\subsubsection{Caída libre con resistencia del aire}

Tómese un cuerpo en caída libre en el cual  suponemos que las únicas fuerzas a las que esta sometido son la gravedad y la resistencia del aire, procedemos a calcular la velocidad bajo la condición $v(0)=v_0$, la aceleración a la que esta sometido el cuerpo es la gravedad y la resistencia del aire se expresa por $-bv$ con $b>0$, por l segunda ley de Newton se obtiene

\formula{
F=m\frac{d}{dt}v =mg-bv,
}

tomando la derivada fraccionara de Caputo con $a=0$ en la expresión anterior y con $0<\alpha<1$, colocaremos un subíndice del lado izquierdo del símbolo de la velocidad para dejar en claro la dependencia con el valor de $\alpha$ y dejar libre el lado derecho por si se consideran diferentes cuerpos

\formula{
m\ifrR{C}{0}{D}{t}{\alpha}{}_\alpha v=mg-b{}_\alpha v,
}

aplicando la transformada de Laplace y la ecuación \eqref{lap5}

\formula{
m\mathcal{L}\left\{ \ifrR{C}{0}{D}{t}{\alpha}{}_\alpha v \right\} &=& \mathcal{L}\left\{ mg-b{}_\alpha v \right\}\\
ms^\alpha\mathcal{L}\left\{{}_\alpha v \right\}-m s^{\alpha-1} v(0)&=& mg\frac{1}{s}-b\mathcal{L}\left\{ {}_\alpha v \right\} \\
(ms^\alpha+b) \mathcal{L}\left\{{}_\alpha v \right\}&=& mg\frac{1}{s}+ m v_0 s^{\alpha-1}\\
 \mathcal{L}\left\{{}_\alpha v \right\}&=&mg\frac{1}{s(ms^\alpha+b)}+ m v_0 \frac{s^{\alpha-1}}{(ms^\alpha+b)}\\
 &=&\frac{g}{\gam{\alpha}}\frac{\gam{\alpha}}{s^\alpha} \frac{s^{\alpha -1}}{\left( s^\alpha+\frac{b}{m}\right)}+  v_0 \frac{s^{\alpha-1}}{\left( s^\alpha+\frac{b}{m}\right) },
}

ocupando las ecuaciones \eqref{con3} y \eqref{con8} así como la convolución

\formula{
 \mathcal{L}\left\{{}_\alpha v \right\} &=&\frac{g}{\gam{\alpha}}\frac{\gam{\alpha}}{s^\alpha} \frac{s^{\alpha -1}}{\left( s^\alpha+\frac{b}{m}\right)}+ v_0 \frac{s^{\alpha-1}}{\left( s^\alpha+\frac{b}{m}\right) }\\
 &=&\frac{g}{\gam{\alpha}}\mathcal{L}\left\{ t^{\alpha-1}\right\} \mathcal{L}\left\{ E_\alpha\left(- \frac{b}{m} t^\alpha\right)\right\}+ v_0\mathcal{L}\left\{ E_\alpha\left( -\frac{b}{m} t^\alpha\right)\right\}\\
 &=&\frac{g}{\gam{\alpha}}\mathcal{L}\left\{ t^{\alpha-1}* E_\alpha\left(- \frac{b}{m} t^\alpha\right)\right\}+ v_0\mathcal{L}\left\{ E_\alpha\left(- \frac{b}{m} t^\alpha\right)\right\},
}

aplicando la transformada inversa de Laplace a la ecuación anterior

\formula{
 {}_\alpha v(t)  &=&\frac{g}{\gam{\alpha}}t^{\alpha-1}* E_\alpha\left(- \frac{b}{m} t^\alpha\right)+ v_0 E_\alpha\left( -\frac{b}{m} t^\alpha\right)\\
&=&\frac{g}{\gam{\alpha}}\int_0^t (t-u)^{\alpha-1} E_\alpha\left( -\frac{b}{m} u^\alpha\right)du+ v_0 E_\alpha\left(- \frac{b}{m} t^\alpha\right)\\
&=&\frac{g}{\gam{\alpha}}\int_0^t (t-u)^{\alpha-1} \sum_{k=0}^\infty \frac{1}{\gam{\alpha k+1}} \left( -\frac{b}{m} u^\alpha\right)^k du+ v_0 E_\alpha\left(- \frac{b}{m} t^\alpha\right)\\
&=&\frac{g}{\gam{\alpha}} \sum_{k=0}^\infty \frac{1}{\gam{\alpha k+1}} \left( -\frac{b}{m} \right)^k  t^{\alpha -1} \int_0^t \left( 1-\frac{u}{t} \right)^{\alpha-1}  u^{\alpha k} du+ v_0 E_\alpha\left(- \frac{b}{m} t^\alpha\right),
}

tomando el cambio de variable $s=u/t$

\formula{
 {}_\alpha v(t)  
&=&\frac{g}{\gam{\alpha}} \sum_{k=0}^\infty \frac{1}{\gam{\alpha k+1}} \left( -\frac{b}{m} \right)^k  t^{\alpha -1} \int_0^t \left( 1-\frac{u}{t} \right)^{\alpha-1}  u^{\alpha k} du+ v_0 E_\alpha \left(- \frac{b}{m} t^\alpha\right)\\
&=&g \sum_{k=0}^\infty \frac{1}{\gam{\alpha}\gam{\alpha k+1}} \left( -\frac{b}{m} \right)^k  t^{\alpha k + \alpha } \int_0^1 \left( 1-s \right)^{\alpha-1}  s^{\alpha k} ds+ v_0 E_\alpha\left(- \frac{b}{m} t^\alpha\right)\\
&=&g \sum_{k=0}^\infty \frac{1}{\gam{\alpha}\gam{\alpha k+1}} \left( -\frac{b}{m} \right)^k  B(\alpha,\alpha k +1)  t^{\alpha k + \alpha } + v_0 E_\alpha\left(- \frac{b}{m} t^\alpha\right)\\
&=&g t^\alpha \sum_{k=0}^\infty \frac{1}{\gam{\alpha k+\alpha +1}} \left( -\frac{b}{m} t^\alpha \right)^k    + v_0 E_\alpha\left(- \frac{b}{m} t^\alpha\right),
}

entonces

\formul{
 {}_\alpha v(t)  
&=&g t^\alpha E_{\alpha,\alpha+1} \left( -\frac{b}{m} t^\alpha \right)    + v_0 E_a\left(- \frac{b}{m} t^\alpha\right),
}{\label{cl1}}

tomando  $\alpha=1$ en la ecuación \eqref{cl1}

\formula{
 {}_1 v(t)  
&=&g t E_{1,2} \left( -\frac{b}{m} t \right)    + v_0 E_1\left(- \frac{b}{m} t\right)\\
&=&g t \left[ \dfrac{\ex{ -\dfrac{b}{m} t }-1}{-\dfrac{b}{m} t}\right]    + v_0 \ex{- \frac{b}{m} t}\\
&=&-\frac{mg}{b}  \left[ \ex{ -\dfrac{b}{m} t }-1\right]    + v_0 \ex{- \frac{b}{m} t}\\
&=& \frac{mg}{b} + \left[ v_0-\frac{mg}{b}\right] \ex{ -\frac{b}{m}t},
}

que corresponde a la solución para el caso clásico.

\begin{figure}[!ht]
\includegraphics[height=0.35\textheight,width=\textwidth]{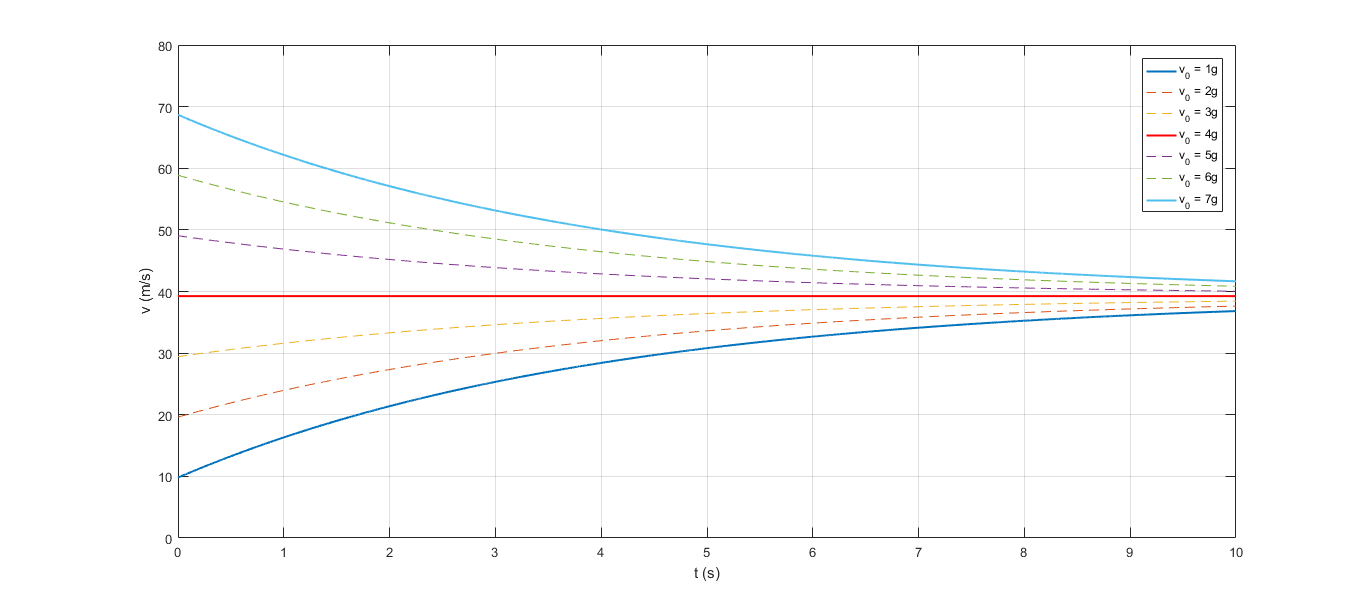}
\caption{Gráfica de ${}_1v(t)$, tomando $m/b= 4s, \ g=9.81 m/s^2$, para 7 velocidades iniciales distintas. De la gráfica podemos apreciar que cuando $v_0$ es mayor que la velocidad terminal el objeto reduce su velocidad, mientras que para el caso en que $v_0$ es menor que la velocidad terminal el objeto aumenta su velocidad.
}
\end{figure}

\subsection{Definiciones Básicas de las Derivadas Fraccionarias}

\subsubsection{Introducción a los Enfoques de Riemann-Liouville y Weyl}

A continuación se consideraran algunos puntos de partida para comenzar el estudio del cálculo fraccional. Se comienza con una generalización de la integral iterada. Si \form{f} es una función integrable localmente en el intervalo \form{(a,\infty)}, la integral \form{n}-foleada esta dada por \cite{hilfer00}

\formula{
\ifr{a}{I}{x}{n} f(x)&:=& \int_a^x du_1 \int_a^{u_1}du_2 \cdots \int_a^{u_{n-1}}f(u_n)du_n\\
&=&\frac{1}{(n-1)!}\int_a^x(x-u)^{n-1}f(u)du,
}

tomando en cuenta que \form{(n-1)!=\gam{n}}, de manera natural se puede obtener una generalización de la integral de \form{f} para un orden arbitrario \form{\alpha >0}

\begin{eqnarray}\label{eq:22}
\ifr{a}{I}{x}{\alpha} f(x) = \frac{1}{\gam{\alpha}}\int_a^x (x-u)^{\alpha -1}f(u)du, \ \mbox{(derecha)},
\end{eqnarray}

de manera similar en el caso en que \form{f} sea una función integrable localmente en el intervalo \form{(-\infty,b)} se tiene

\begin{eqnarray}\label{eq:23}
\ifr{x}{I}{b}{\alpha} f(x)=\frac{1}{\gam{\alpha}}\int_x^b(u-x)^{\alpha -1}f(u)du, \ \mbox{(izquierda)},
\end{eqnarray}

ambas formas están definidas para \form{f}. Cuando \form{a = -\infty} la ecuación \eqref{eq:22} es equivalente a la definición de Liouville, y cuando \form{a =0} se tiene la definición de Riemann (sin la función complementaria). Generalmente se habla de \form{\ifr{a}{I}{x}{\alpha} f} como la integral fraccionaria de Riemann-Liouville de orden \form{\alpha} de \form{f}. Por otra parte las expresiones

\formula{
\ifr{x}{W}{\infty}{\alpha} f(x)&=&\ifr{x}{I}{\infty}{\alpha} f(x)\\
&=& \frac{1}{\gam{\alpha}}\int_x^\infty(u-x)^{\alpha -1} f(u)du,
}

\formula{
\ifr{-\infty}{W}{x}{\alpha} f(x)&=& \ifr{-\infty}{I}{x}{\alpha} f(x)\\
&=& \frac{1}{\gam{\alpha}}\int_{-\infty}^x(x-u)^{\alpha -1}f(u)du,
}

son conocidas como  integrales fraccionarias de Weyl de orden \form{\alpha}.

Las integrales fraccionarias derecha e izquierda \form{\ifr{a}{I}{x}{\alpha} f(x) } y \form{\ifr{x}{I}{b}{\alpha} f(x)} están relacionadas por la igualdad de Parseval (integración fraccionaria por partes) \cite{hilfer00} que se da de forma conveniencia para los casos en que \form{a=0} y \form{b=\infty}

\formula{
\int_0^\infty f(x)\left({}_0I_x^\alpha g \right)(x) dx=\int_0^\infty \left({}_xW_\infty^\alpha f \right)(x)g(x)dx.
}

Las siguientes propiedades son validas para integrales fraccionarias derechas (para el caso de integrales fraccionarias izquierdas se presentan algunos cambios).

Con respecto a la existencia de integrales fraccionarias para \form{f\in L_{loc}^1(a,\infty)}. Si \form{a>-\infty}, la integral fraccionaria \form{\ifr{a}{I}{x}{\alpha} f(x)} es finita en casi  todas partes del intervalo \form{(a,\infty)} y pertenece a \form{L_{loc}^1(a,\infty)}. Si \form{a=-\infty}, se asume que \form{f} se comporta en \form{-\infty} de tal manera que la integral en la ecuación \eqref{eq:23} es convergente bajo las mismas suposiciones, las integrales fraccionarias satisfacen la propiedad de semigrupo

\formula{
\ifr{a}{I}{x}{\alpha} \ifr{a}{I}{x}{\beta} = \ifr{a}{I}{x}{\alpha + \beta},\ (\alpha,\beta>0)
}

la propiedad de semigrupo se puede probar haciendo uso de la fórmula de Dirichlet relacionada con el cambio en el orden de integración

\formula{
\ifr{a}{I}{x}{\alpha} \ifr{a}{I}{x}{\beta} f(x) &=&\frac{1}{\gam{\alpha}}\int_a^x(x-u)^{\alpha -1} \left[\frac{1}{\gam{\beta}}\int_a^u (u-t)^{\beta -1}f(t)dt \right]du\\
&=& \frac{1}{\gam{\alpha}\gam{\beta} }\int_a^x f(t)dt \int_t^u (x-u)^{\alpha -1}(u-t)^{\beta -1}du
}

tomando el cambio de variable \form{y=\frac{u-t}{x-t}} en la segunda integral de la derecha

\formula{
\ifr{a}{I}{x}{\alpha} \ifr{a}{I}{x}{\beta} f(x)&=& \frac{1}{\gam{\alpha}\gam{\beta} }\int_a^x f(t)dt \int_t^u (x-u)^{\alpha -1}(u-t)^{\beta -1}du\\
&=& \frac{1}{\gam{\alpha}\gam{\beta} }\int_a^x (x-t)^{\alpha + \beta -1} f(t)dt \int_0^1 (1-y)^{\alpha -1}y^{\beta -1}dy \\
&=& \frac{1}{\gam{\alpha}\gam{\beta} }\int_a^x (x-t)^{\alpha + \beta -1} f(t)dt B(\alpha,\beta)\\
&=& \frac{1}{\gam{\alpha}\gam{\beta} } \frac{\gam{\alpha}\gam{\beta}}{\gam{\alpha + \beta}} \int_a^x (x-t)^{\alpha + \beta -1} f(t)dt \\
&=&  \frac{1}{\gam{\alpha + \beta}} \int_a^x (x-t)^{\alpha + \beta -1} f(t)dt, 
}

en particular, se tiene

\formula{
\ifr{a}{I}{x}{n+\alpha}=\ifr{a}{I}{x}{n}\ifr{a}{I}{x}{\alpha}f, \ (n\in\mathds{N}, \alpha>0),
}

lo que implica una diferenciación \form{n}-foleada para casi cualquier \form{x}

\formula{
\der{d}{x}{n}\ifr{a}{I}{x}{n+\alpha}f(x)=\ifr{a}{I}{x}{\alpha}f(x), \ (n\in \mathds{N}, \alpha >0).
}

Los resultados anteriores también son válidos para valores complejos \form{\alpha} con \form{Re(\alpha)>0}. Entonces  \form{\ifr{a}{I}{x}{\alpha}}  sera considerada como una función holomorfa de \form{\alpha} con \form{Re(\alpha)>0} que puede tomarse en todo el plano complejo si es extendida de forma analítica para \form{f} lo suficientemente suave.

Para entender este hecho, se asume por conveniencia, que \form{f} es una función infinitamente diferenciable definida en \form{\mathds{R}} con soporte compacto contenido en \form{[a,\infty]}, si \form{ \alpha> -\infty},  \form{f^{(n)}(a)=0} para \form{n\in \mathds{N}\cup\{0\}}. Entonces, para cualquier \form{x>a} la integral \eqref{eq:22} es una función holomorfa en \form{\alpha} para \form{Re(\alpha)>0}. Ahora, la integración por partes \form{n}-veces da como resultado

\begin{eqnarray}\label{eq:24}
\ifr{a}{I}{x}{\alpha}f(x)=\ifr{a}{I}{x}{n+\alpha}f^{(n)}(x), \ (n\in \mathds{N},Re(\alpha) >0).
\end{eqnarray}

Aplicando la propiedad de semigrupo a la expresión de la derecha en la ecuación \eqref{eq:24} y diferenciando el resultado \form{n}-veces con respecto a \form{x}, se obtiene

\begin{eqnarray}\label{eq:25}
\der{d}{x}{n}\ifr{a}{I}{x}{\alpha}f(x)=\ifr{a}{I}{x}{\alpha}f^{(n)}(x), \ (n\in \mathds{N}, Re(\alpha)>0),
\end{eqnarray}

mostrando que bajo las hipótesis asumidas las operaciones de integración de orden fraccionario \form{\alpha} y diferenciación de orden entero \form{n} conmutan.

Regresando a la ecuación \eqref{eq:24} nos damos cuenta de que su lado derecho es una función holomorfa de \form{\alpha} en el dominio más amplio \form{\conj{\alpha \in \mathds{C}, Re(\alpha)>-n}}. Así podemos extender analíticamente \form{\ifr{a}{I}{x}{\alpha}f(x)} en el dominio \form{\conj{\alpha\in \mathds{C}, \re{\alpha}\leq 0}}, definiendo para \form{\alpha\in \mathds{C}} con \form{\re{\alpha}\leq 0}, 

\begin{eqnarray}\label{eq:26}
\ifr{a}{I}{x}{\alpha}f(x)&:=& \ifr{a}{I}{x}{n+\alpha}f^{(n)}(x)\nonumber \\
&=& \der{d}{x}{n}\ifr{a}{I}{x}{n+\alpha}f(x),
\end{eqnarray}

para cualquier entero \form{n>-\re{\alpha}}. En particular, se tiene que

\formula{
\ifr{a}{I}{x}{0}f(x)&=& f(x),\\
\ifr{a}{I}{x}{-n}f(x)&=& f^{(n)}(x), \ (n\in \mathds{N}).
}

El método elegante de la extension analítica desarrollado por Riesz \cite{hilfer00}  se limita a una clase de funciones bastante pequeña. Pero las expresiones que ocurren en la ecuación \eqref{eq:26} son significativamente para clases mucho más generales de funciones y esto da lugar a las siguientes definiciones de derivadas fraccionarias que se remontan a Liouville.

Sea \form{\alpha\in \mathds{C}} con \form{\re{\alpha}>0} y \form{\lfloor\re{\alpha}\rfloor+1}, donde \form{\lfloor\re{\alpha}\rfloor} representa la parte entera de \form{\re{\alpha}}. Entonces, la derivada fraccionaria de orden \form{\alpha} esta definida  por

\begin{eqnarray}\label{eq:27}
\ifr{a}{D}{x}{\alpha}f(x) = \der{d}{x}{n}\ifr{a}{I}{x}{n-\alpha}f(x), \ (n=\lfloor\re{\alpha}\rfloor+1),
\end{eqnarray}

para cualquier \form{f\in L_{loc}^1(a,\infty)}.

Se pueden unificar las definiciones sobre integrales y derivadas de orden arbitrario \form{\alpha, \re{\alpha}\neq 0}, para \form{n\in \mathds{N}},

\formula{
\ifr{a}{D}{x}{\alpha}f(x) = \left\{
\begin{array}{ll}
\frac{1}{\gam{-\alpha}}\int_a^x(x-u)^{-\alpha -1}f(u)du, &(\re{\alpha}<0)\\ \\
\der{d}{x}{n}\ifr{a}{I}{x}{n-\alpha}f(x), & (\re{\alpha}>0, \re{\alpha}\in (n-1,n)) 
\end{array}
\right.
}

esta expresión suele ser conocida como diferointegral de \form{f} de orden \form{\alpha} o también se conoce como integrodiferenciación fraccionaria de orden \form{\alpha}.

Nótese que la derivada fraccionaria izquierda de orden \form{\alpha} se define por

\formula{
\ifr{x}{D}{b}{\alpha}f(x)=(-1)^n \der{d}{x}{n}\ifr{x}{I}{b}{n-\alpha}f(x), 
}

con \form{ n=\lfloor\re{\alpha}\rfloor+1}.  La derivada fraccionaria de orden  imaginario \form{\alpha=i\theta, \theta\neq 0}, se define como

\formula{
\ifr{a}{D}{x}{i\theta}f(x)=\frac{1}{\gam{1-i\theta}}\frac{d}{dx}\int_a^x   (x-u)^{-i\theta}f(u)du,
}

debido que la integral fraccionaria de la ecuación \eqref{eq:22} diverge para \form{\alpha=i\theta} \cite{hilfer00}, la integral fraccionaria de orden  \form{\alpha=i\theta} esta definida por

\formula{
\ifr{a}{I}{x}{i\theta}f(x)&=&\frac{d}{dx}\ifr{a}{I}{x}{1+i\theta}f(x)\\
&=& \frac{1}{\gam{1+i\theta}}\frac{d}{dx}\int_a^x(x-u)^{i\theta}f(u)du,
}

la definición de integro-diferenciación fraccionaria para \form{\alpha \in \mathds{C}} se completa con la introducción del operador identidad 

\formula{
\ifr{a}{D}{x}{0}f:= \ifr{a}{I}{x}{0}=f,
}

los operadores fraccionarios tambien son lineales

\formula{
\ifr{a}{D}{x}{\alpha}[c_1f_1(x)+c_2f_2(x)]=c_1\ifr{a}{D}{x}{\alpha}f_1(x)+c_2\ifr{a}{D}{x}{\alpha} f_2(x),
}

con \form{c_1,c_2}  constantes.

Con respecto a condiciones suficientes para la existencia de las derivadas fraccionarias del tipo de la ecuación \eqref{eq:27} y su relación con la ecuación \eqref{eq:29}, considérese el caso \form{0<\re{\alpha}<1}, con \form{\alpha>-\infty}. Suponiendo que \form{f} es absolutamente continua en el intervalo finito  \form{[a,b]}, esto es \form{f\in AC[a,b]}, significa que \form{f} es diferenciable casi en todas partes en el intervalo \form{(a,b)} con \form{f^{(1)}\in L^1(a,b)} y tiene la representación en \form{[a,b]}

\formula{
f(x)=\int_a^x f^{(1)}(u)du + f(a) = \ifr{a}{I}{x}{1}f^{(1)}(x) + f(a).
}

Sustituyendo esta expresión en \form{\ifr{a}{I}{x}{1-\alpha}}f(x) y notando que por la propiedad de semigrupo los operadores \form{\ifr{}{I}{}{1-\alpha}} y \form{\ifr{}{I}{}{1}}  conmutan, se obtiene

\formula{
\ifr{a}{I}{x}{1-\alpha}f(x)=\ifr{a}{I}{x}{1}\ifr{a}{I}{x}{1-\alpha}
f^{(1)}(x)+\frac{1}{\gam{2-\alpha}}f(a)(x-a)^{1-\alpha}.
}

Al diferenciar con respecto a \form{x} se obtiene

\begin{eqnarray}\label{eq:28}
\ifr{a}{D}{x}{\alpha}f(x)&=& \frac{d}{dx}\ifr{a}{I}{x}{1-\alpha}f(x)\nonumber \\
&=& \ifr{a}{I}{x}{1-\alpha}f^{(1)}(x)+\frac{1}{\gam{1-\alpha}}f(a)(x-a)^{-\alpha},
\end{eqnarray}

Lo que demuestra que, en general, los operadores \form{\ifr{a}{I}{x}{1-\alpha}} y \form{\frac{d}{dx}} no conmutan.

Por la desigualdad de Hölder \cite{hilfer00} se deriva fácilmente de la ecuación \eqref{eq:28} que \form{\ifr{a}{D}{x}{\alpha}f\in L^r(a,b)} para \form{1<r<\frac{1}{\re{\alpha}}}.

La expresión en la ecuación \eqref{eq:28} puede extenderse a \form{\alpha} con \form{\re{\alpha}\geq 1}. Los resultados se resumen en la siguiente proposición. De antemano se introduce la siguiente notación: Para \form{n\in \mathds{N}, AC^{n-1}[a,b]} representa el conjunto de las funciones \form{f} que son \form{(n-1)}-veces diferenciables  en \form{[a,b]} tal que \form{f,f^{(1)},\cdots, f^{(n-1)}} son absolutamente continuas en \form{[a,b]}. Teniendo en cuenta que \form{AC^0[a,b]} es igual a \form{AC[a,b]} .

\begin{prop}
\form{{}}

\begin{itemize}
\item[a)] Si \form{f\in AC[a,b]} esta dada en el intervalo finito \form{[a,b]}, entonces \form{\ifr{a}{D}{x}{\alpha}f} y \form{\ifr{x}{D}{b}{\alpha}f} existen para \form{0<\re{\alpha}<1}. Además, \form{\ifr{a}{D}{x}{\alpha}f\in L^r} para \form{1\leq r < \frac{1}{\re{\alpha}}} con

\formula{
\ifr{a}{D}{x}{\alpha}f(x)=\frac{1}{\gam{1-\alpha}} \left[ \frac{f(a)}{(x-a)^\alpha}+\int_a^x (x-u)^{-\alpha}f^{(1)}(u)du\right].
}

\item[b)] Si \form{f\in AC^{n-1}[a,b], n=\lfloor \re{\alpha}\rfloor +1}, entonces \form{\ifr{a}{D}{x}{\alpha}f} existe para \form{\re{\alpha}\geq 0} y tiene la representación

\formula{
\ifr{a}{D}{x}{\alpha}f(x)=\sum_{k=0}^{n-1}\frac{f^{(k)}(a)}{\gam{1+k-\alpha}}(x-a)^{k-\alpha}+\frac{1}{\gam{n-\alpha}}\int_a^x (x-u)^{n-\alpha-1}f^{(n)}(u)du.
}

Una forma alternativa de definir una derivada fraccionaria de orden \form{\alpha}, también debido a Liouville, es

\formul{
\ifr{a}{\overline{D}}{x}{\alpha}f(x):= \ifr{a}{I}{x}{n-\alpha}f^{(n)}(x), 
}{\label{eq:29}}

con \form{n=\lfloor\re{\alpha}\rfloor+1}, además \form{f} tiene que ser \form{n}-veces diferenciable para que el lado derecho de la ecuación \eqref{eq:29} exista.

\end{itemize}

\end{prop}


\subsubsection{Introducción a la Derivada Fraccionaria de Caputo }

Michele Caputo \form{(1969)} publicó un libro  en el que introdujo una nueva derivada fraccionaria, que había sido descubierta de forma independiente por Gerasimov \form{(1948)}. Esta derivada fraccionaria es de suma importancia, ya que permite dar una interpretación física a problemas de condiciones iniciales, además se utiliza para modelar el tiempo fraccionario. En algunos textos se conoce como la derivada fraccionaria de Gerasimov-Caputo.

Sea \form{[a, b]} un intervalo finito de la recta real \form{\mathds{R}}, para \form{ \alpha \in \mathds{C} (\re{\alpha} \geq 0)}. Las derivadas derecha e izquierda  fraccionarias de Caputo  se definen como:

\formula{
\left(\ifrR{C}{a}{D}{x}{\alpha}y\right)(x)&=& \left(\ifrR{RL}{a}{I}{x}{n-\alpha}D^n y \right)(x) \\
&=&\dfrac{1}{\gam{n-\alpha}}\int_{a}^{x} (x-t)^{n-\alpha -1} y^{(n)}(t)dt ,
}
\formula{
\left(\ifrR{C}{x}{D}{b}{\alpha}y\right)(x)&=&(-1)^n\left(\ifrR{RL}{x}{I}{b}{n-\alpha}D^n y \right)(x)\\
&=&\dfrac{(-1)^n}{\gam{n-\alpha}}\int_{x}^{b} (t-x)^{n-\alpha -1} y^{(n)}(t)dt,
}

donde \form{D=\frac{d}{dx}} y \form{n=-\lfloor -\re{\alpha}\rfloor}, i.e., \form{n=\re{\alpha}+1} para \form{\alpha \in \mathds{N}_0} y \form{n=\alpha} para \form{\alpha \notin\mathds{N}_0}. Y si \form{0<\re{\alpha}<1}

\formula{
\left(\ifrR{C}{a}{D}{x}{\alpha} y\right)(x)&=&\left(\ifrR{RL}{a}{I}{x}{1- \alpha}D y\right)(x)\\
&=&\dfrac{1}{\gam{1-\alpha}}\int_{a}^{x} (x-t)^{-\alpha} y^{(1)}(t)dt,
}

\formula{
\left(\ifrR{C}{x}{D}{b}{\alpha} y\right)(x)&=&-\left(\ifrR{RL}{x}{I}{b}{1- \alpha}D y\right)(x)\\
&=&\dfrac{-1}{\gam{1-\alpha}}\int_{x}^{b} (t-x)^{-\alpha} y^{(1)}(t)dt.
}

La conexión entre las derivadas fraccionarias de Caputo y  de Riemann está dada por las relaciones:

\formul{
\left(\ifrR{C}{a}{D}{x}{\alpha}y\right)(x)&=&\left(\ifrR{RL}{a}{D}{x}{\alpha}{\left[y(t)-\sum_{k=0}^{n-1}\dfrac{y^{(k)}(a)}{k!}(t-a)^k\right]}\right)(x),
}{\label{eq:44}}

\formul{
\left(\ifrR{C}{x}{D}{b}{\alpha}y\right)(x)&=&\left(\ifrR{RL}{x}{D}{b}{\alpha}{\left[y(t)-\sum_{k=0}^{n-1}\dfrac{y^{(k)}(b)}{k!}(b-t)^k\right]}\right)(x),
}{\label{eq:45}}

En particular, si \form{ 0 <\re{\alpha} <1 }, las relaciones \eqref{eq:44} y \eqref{eq:45} toman las siguientes formas:

\formula{
\left(\ifrR{C}{a}{D}{x}{\alpha}y\right)(x)&=&\left(\ifrR{RL}{a}{D}{x}{\alpha}{\left[y(t)-y(a)\right]}\right)(x),\\
\left(\ifrR{C}{x}{D}{b}{\alpha}y\right)(x)&=&\left(\ifrR{RL}{x}{D}{b}{\alpha}{\left[y(t)-y(b)\right]}\right)(x).
}

Para \form{\alpha = n}, las derivadas fraccionarias de Caputo corresponden a las derivadas clásicas, excepto por el signo de la derivada derecha.

Sin embargo, para \form{ k = 0,1, \dots, n-1 }, se tiene:

\formula{
\left(\ifrR{C}{a}{D}{x}{\alpha}(t-a)^k\right)(x)=0,\\\left(\ifrR{C}{x}{D}{b}{\alpha}(b-t)^k\right)(x)=0,
}

en particular

\formula{
\left(\ifrR{C}{a}{D}{x}{\alpha}1\right)(x)=0,\\
\left(\ifrR{C}{x}{D}{b}{\alpha}1\right)(x)=0.
}

Por otro lado, si \form{ \re{\alpha}> 0}  y  \form{\lambda> 0}, entonces

\formula{
\left(\ifrR{C}{a}{D}{x}{\alpha}{\ex{\lambda t}}\right)(x)\neq\lambda^\alpha\ex{\lambda t}, \ (\alpha\in \mathds{R}).
}

Las derivados fraccionarias de Caputo se comportan como operadores inversos por la izquierda para las integrales fraccionarias de Riemann-Liouville \form{ \ifrR{RL}{a}{I}{x}{\alpha} } y \form{\ifrR{RL}{x}{I}{b}{\alpha}}, si \form{ \re{\alpha}> 0} y \form{ y (x) \in C [a, b] }

\formula{
\left(\ifrR{C}{a}{D}{x}{\alpha}{\left(\ifrR{RL}{a}{I}{x}{\alpha} y\right)}\right)(x)=y(x),\\
\left(\ifrR{C}{x}{D}{b}{\alpha}{\left(\ifrR{RL}{x}{I}{b}{\alpha} y\right)}\right)(x)=y(x).
}

Por otro lado, si \form{ \re{\alpha}> 0 } y \form{ n = - \lfloor - \re{\alpha} \rfloor }, entonces para condiciones apropiadas para \form{ y (x)}

\formula{
\left(\ifrR{RL}{a}{I}{x}{\alpha}{\left(\ifrR{C}{a}{D}{x}{\alpha}y\right)}\right)(x)&=&y(x)-\sum_{k=0}^{n-1}\dfrac{y^{(k)}(a)}{k!}(x-a)^k,\\
\left(\ifrR{RL}{x}{I}{b}{\alpha}{\left(\ifrR{C}{x}{D}{b}{\alpha}y\right)}\right)(x)&=&y(x)-\sum_{k=0}^{n-1}(-1)^k\dfrac{y^{(k)}(b)}{k!}(b-x)^k.
}

En particular si, \form{ 0 <\re{\alpha} \leq 1 }, entonces

\formula{
\left(\ifrR{RL}{a}{I}{x}{\alpha}{\left(\ifrR{C}{a}{D}{x}{\alpha}y\right)}\right)(x)&=&y(x)-y(a),\\
\left(\ifrR{RL}{x}{I}{b}{\alpha}{\left(\ifrR{C}{x}{D}{b}{\alpha}y\right)}\right)(x)&=&y(x)-y(b).
}

En sus primeros artículos y varios después de eso, Caputo utilizó una transformada de Laplace de la derivada fraccionaria de Caputo, que esta dada por:

\formula{
\left(\mathcal{L}{\left\lbrace \ifrR{C}{0}{D}{x}{\alpha}y \right\rbrace }\right)(s)&=& s^\alpha(\mathcal{L}{y})(s)-\sum_{k=0}^{n-1}s^{\alpha-k-1}(D^ky)(0).
}

Cuando \form{ 0 < \alpha \leq 1}, entonces

\formula{
\left(\mathcal{L}{\left\lbrace \ifrR{C}{0}{D}{x}{\alpha}y \right\rbrace }\right)(s)&=& s^\alpha(\mathcal{L}{y})(s)-s^{\alpha-1}y(0).
}

Estas derivadas fraccionarias se pueden definir sobre todo el eje real dando como resultado las expresiones:

\formula{
\left(\ifrR{C}{-\infty}{D}{x}{\alpha}y\right)(x)&=&\dfrac{1}{\gam{n-\alpha}}\int_{-\infty}^{x} (x-t)^{n-\alpha -1} y^{(n)}(t)dt ,\\
\left(\ifrR{C}{x}{D}{\infty}{\alpha}y\right)(x)&=&\dfrac{(-1)^n}{\gam{n-\alpha}}\int_{x}^{\infty} (t-x)^{n-\alpha -1} y^{(n)}(t)dt,
}

con \form{x\in\mathds{R}}.

\subsubsection{Regla de Leibniz }

La regla clásica de Leibniz para la derivada $n$-ésima de un producto de dos funciones $f, g$ como suma de productos de operaciones realizadas en cada función viene dada, siempre que $f$ y $g$ sean $n$-veces diferenciables en $z$, por

\formula{
D_z^n \left(f(z)g(z) \right)=\sum_{k=0}^n \binom{n}{k}D_z^{n-k}f(z)D_z^kg(z),
}

donde \form{D_z^n =\left( \frac{d}{dz}\right)^n}.Esta regla puede extenderse a valores fraccionarios de $ \alpha $: reemplazando $n$ por $\alpha $ tenemos para funciones holomorfas $f,g$

\formul{
\ifr{a}{D}{z}{\alpha} \left(f(z)g(z) \right)=\sum_{k=0}^\infty \binom{\alpha}{k}\ifr{a}{D}{z}{\alpha-k}f(z)\ifr{a}{D}{z}{k}g(z), \ \  \alpha\in \mathds{C},
}{\label{r1}}

Un resultado que básicamente se remonta a Liouville $(1832)$, donde $\ifr{a}{D}{z}{\alpha}$ se entiende ahora en el sentido de Osler. Esta fórmula sufre el aparente inconveniente de que el intercambio de $f (z)$ y $g (z)$ en el lado derechos de \eqref{r1} no es obvio. Una generalización interesante de esta regla sin el inconveniente, debido a Y.Watanabe \form{(1931)} y Osler \form{(1970)}, es

\formul{
\ifr{a}{D}{z}{\alpha} \left(f(z)g(z) \right)=\sum_{k=-\infty}^\infty \binom{\alpha}{k+\mu}\ifr{a}{D}{z}{\alpha-k-\mu}f(z)\ifr{a}{D}{z}{k+\mu}g(z), \ \  \alpha\in \mathds{C},
}{\label{r2}}

donde $\mu$ es un número arbitrario, racional, irracional o complejo. El caso especial $\mu = 0 $ se reduce a \eqref{r1}.

Observe que si $\re{\alpha} <0 $, entonces la fórmula \eqref{r1} es realmente la contraparte de la regla de Leibniz para las integrales fraccionales. Existe también una regla Leibniz fraccional simétrica en la forma

\formul{
\ifr{a}{D}{z}{\alpha} \left(f(z)g(z) \right)=\sum_{k=-\infty}^\infty  c\binom{\alpha}{ck+\mu}\ifr{a}{D}{z}{\alpha-ck-\mu}f(z)\ifr{a}{D}{z}{ck+\mu}g(z), \ \  \alpha\in \mathds{C},
}{\label{r4}}

donde $ \alpha, \mu \in \mathds{C} $ para el cual $ \binom{\alpha} {ck + \mu} $ está bien definido, y $ 0 <c \leq 1 $. Se obtiene \eqref{r2} para $c = 1$, y se reduce a \eqref{r1} para $c = 1$ y $ \mu = 0 $.

Como casos especiales de \eqref{r2}, se tiene para $f(z) = 1$, y $g (z)$ renombrado como $f(z)$,

\formula{
\ifr{a}{D}{z}{\alpha}f(z)=\frac{\gam{\alpha +1}\si{(\alpha-\mu)\pi}}{\pi}\sum_{k=-\infty}^{\infty}(-1)^k \frac{(z-a)^{k+\mu -\alpha}}{(\alpha-\mu -k)\gam{\mu+k+1}}\ifr{\alpha}{D}{z}{\mu + k}f(z),
}

para $\re{\alpha}>-1$, $\alpha-\mu\notin \mathds{Z}$, notando que $\gam{z}\gam{1-z}=\pi \csc(\pi z)$ para $z\notin \mathds{Z}$.

Un resultado adicional, también debido a Osler, en el caso $ \alpha = 0 $ de \eqref{r2}, es 

\formula{
f(z)g(z)=\frac{\si{\mu \pi}}{\pi}\sum_{k=-\infty}^\infty (-1)^k\frac{\ifr{a}{D}{z}{-(\mu+k)}f(z)\ifr{a}{D}{z}{\mu+k}g(z)}{\mu + k}, \ \ \mu \notin \mathds{Z}.
}

Obsérvese que \eqref {r2} tiene bajo condiciones adecuadas el análogo integral interesante

\formul{
\ifr{a}{D}{z}{\alpha}\left(f(z)g(z) \right)=\int_{-\infty}^\infty \binom{\alpha}{\tau +\mu}\ifr{a}{D}{z}{\alpha-\mu-\tau}f(z)\ifr{a}{D}{z}{\tau + \mu}g(z) d\tau,
}{\label{r3}}

donde $ \alpha, \mu \in\mathds {C} \setminus \mathds{Z}^{-} $. Asume una forma elegante para $\mu = 0 $. Al establecer $\alpha = 0 $ la formula \eqref{r3} se reduce fácilmente a

\formula{
f(z)g(z)=\frac{1}{\pi}\int_{-\infty}^\infty \frac{\si{\pi(\mu + \tau)}}{\mu + \tau}\ifr{a}{D}{z}{-(\mu + \tau)}f(z) \ifr{a}{D}{z}{\mu + \tau}g(z)d\tau, \ \ \mu \notin \mathds{Z}^{-}.
}

Recientemente Kalia y su compañero de trabajo deducen de la regla de Leibniz \eqref{r4} un número de fórmulas de expansión interesantes asociadas con la función gamma $ \gam{z} $, la función Psi $ \psi(z) = \frac{\Gamma'(z)}{\gam{z}} $, con la función Gamma incompleta $ \gamma(a, z) $, definida por

\formula{
\gamma(a,z)=\int_0^a t^{z-t}e^{-t}dt, \ \ \re{z}>0,
}

así como con la función Gamma incompleta entera $ \gamma^* (a, z) $, definida por

\formula{
a^z\gamma^*(a,z)=\frac{\gamma(a,z)}{\gam{z}}, \ \ (|arg(z)|\leq \pi - \epsilon, \ 0<\epsilon <\pi).
}

Uno es el resultado bien conocido con respecto a la función psi

\formula{
\psi(\beta - \alpha +1)-\psi(\beta +1)=\sum_{k=1}^\infty \frac{(-\alpha)_k}{k(\beta -\alpha +1)_k},
}

valido para $\re{\beta}>-1$ y $\beta-\alpha+1\notin \mathds{Z}^{-}$. aquí 

\formula{
(\alpha)_k=\alpha(\alpha+1)\cdots (\alpha +k -1),
}
con $k\in \mathds{N}$. Una expansión para  la función Gamma incompleta entera, esta dada para $\re{\beta}>-1$; $\alpha \notin \mathds{Z}^{-}, \  \alpha-\nu \notin \mathds{Z}$ por

\formul{
\gamma^*(\beta - \alpha,z)=\frac{\si{\pi(\alpha-\nu)}}{\pi}\gam{\alpha +1}\sum_{k=-\infty}^\infty (-1)^k \frac{\gamma^*(\beta-\nu-k,z)}{(\alpha-\nu-k)\gam{\nu+k+1}}.
}{\label{r5}}

El caso especial de \eqref{r5} cuando $ \beta = 0 $, que se reduce fácilmente a

\formula{
\gamma(\alpha,z)=-\frac{\si{\pi(\alpha+\nu)}\si{\pi\nu}}{\pi \si{\pi \alpha}}\sum_{-\infty}^\infty \frac{z^{\alpha+\nu +k}}{\alpha +\nu + k}\gamma(-\nu-k,z),
}

con $ \alpha, \alpha + \nu \notin \mathds{Z} $. El caso particular $ \nu = 0$ de la ecuación \eqref{r5} da

\formula{
\gamma^*(\beta-\alpha,z)=\frac{\si{\pi \alpha}}{\pi}\gam{\alpha + 1}\sum_{k=0}^\infty (-1)^k \frac{\gamma^*(\beta-k,z)}{(\alpha-k)k!},
}

valido para $\re{\beta}>1, \alpha \notin \mathds{Z}$.

Una interesante aplicación de la regla de Leibniz para las funciones hipergeométricas viene dada por la identidad

\formul{
\ifr{2}{F}{1}{}(a,b,c;1)=\frac{\gam{c}\gam{c-a-b}}{\gam{c-a}\gam{c-b}},
}{\label{r6}}

válida para $\re{c}> \re{a+b},\ c \notin \mathds{Z}_0^{-} $, que se deriva directamente de la representación integral que puede ser establecida por la regla de Leibniz

\formula{
\ifr{2}{F}{1}{}(a,b,c;z)=\frac{\gam{c}}{\gam{a}\gam{c-a}}\int_0^1 u^{a-1}(1-u)^{c-a-1}(1-uz)^{-b}du,
}

con $0<\re{a}<\re{c}, \ |z|<1$, o aplicando la regla de Leibniz para las integrales fraccionales al producto de $f(x)=x^\mu$ y $g(x)=x^\lambda, \ \lambda,\mu \geq 0$, se obtiene

\formula{
\frac{\gam{\lambda +\mu +1}}{\gam{\lambda +\mu +\nu +1}}=\frac{\gam{\mu +1}}{\gam{\mu +\nu +1}}\ifr{2}{F}{1}{}(-\lambda,\nu,\mu+\nu+1;1).
}

La notación más convencional con $a=-\lambda,\ b=\nu,\ c=\mu+\nu+1$ da la ecuación \eqref{r6}. Obsérvese que la fórmula \eqref{r6} es en realidad un caso particular del teorema de muestreo de Shannon para el análisis de señales. Otro caso es dado por

\formula{
\ifr{2}{F}{1}{}(a,b,c;z)=(1-z)^{-b} \ifr{2}{F}{1}{}\left(c-a,b,c; \frac{z}{z-1} \right).
}

Añadamos finalmente un análogo integral inusual de la versión fraccionaria del teorema de Taylor

\formula{
f(z)=c\sum_{k=-\infty}^\infty \frac{\ifr{a}{D}{w}{ck+\mu}f(w)}{\gam{ck+\mu +1}}(z-w)^{ck+\mu}, \ \ 0<c\leq 1, \ \mu\in\mathds{C},
}


\subsection{Ejemplos Básicos de la Derivada Fraccionaria}

Comencemos analizando la derivada clásica de un monomio de la forma  \form{(ax+b)^m} con \form{m\in\mathds{Z}} 

\formula{
\der{d}{x}{}(ax+b)^m&=& am(ax+b)^{m-1},\\ \\
\der{d}{x}{2}(ax+b)^m&=&a^2 m(m-1)(ax+b)^{m-2},\\ \\
\der{d}{x}{3}(ax+b)^m&=& a^3m(m-1)(m-2)(ax+b)^{m-3},\\
&\vdots &\\
\der{d}{x}{n}(ax+b)^m&=&a^nm(m-1)(m-2)\cdots[m-(n-1)] (ax+b)^{m-n},
}

multiplicando tanto el numerador como el denominador de la ultima expresión por \form{(m-n)!} se obtiene

\formula{
\der{d}{x}{n}(ax+b)^m&=&a^n m(m-1)(m-2)\cdots[m-(n-1)]\frac{(m-n)!}{(m-n)!} (ax+b)^{m-n}\\
&=& a^n\frac{m!}{(m-n)!}(ax+b)^{m-n},
}

reemplazando los enteros \form{m} y \form{n} por números arbitrarios \form{\mu} y \form{\nu }, así como utilizando la función Gamma se obtiene

\formul{
\der{d}{x}{\nu}(ax+b)^\mu = a^\nu\frac{\Gamma(\mu+1)}{\Gamma(\mu-\nu+1)}(ax+b)^{\mu-\nu}.
}{\label{eq:34}}

Para el siguiente ejemplo se evalua la derivada fraccionaria \form{\ifr{0}{D}{x}{\alpha}f(x)} tomando la función \form{f(x)=x^\mu} con \form{\mu>-1}. Haciendo uso de la función Beta, válida para \form{p,q\in \mathds{C}} con \form{\re{p},\re{q}>0} \cite{arfken85}

\formula{
B(p,q)&=& \int_0^1 u^{p-1}(1-u)^{q-1}du\\
&=& \frac{\gam{p}\gam{q}}{\gam{p+q}},
}

se tiene para \form{m-1\leq \alpha <m} con \form{ m\in \mathds{N}},

\formula{
\ifr{0}{D}{x}{\alpha}x^\mu&=& \der{d}{x}{m}\left[\frac{1}{\gam{m-\alpha}}\int_0^x(x-u)^{m-\alpha-1}u^\mu du \right]
}

tomando el cambio de variable \form{u=xv}

\formula{
\ifr{0}{D}{x}{\alpha}x^\mu &=& \frac{1}{\gam{m-\alpha}}\der{d}{x}{m}\left[x^{m-\alpha + \mu}\int_0^1 (1-v)^{m-\alpha-1}v^\mu dv \right] 
\\
&=& \frac{1}{\gam{m-\alpha}}B(m-\alpha, \mu+1)\der{d}{x}{m}x^{m-\alpha+\mu}\\
&=& \frac{1}{\gam{m-\alpha}}\frac{\gam{m-\alpha}\gam{\mu+1}}{\gam{m-\alpha+\mu+1}}\der{d}{x}{m}x^{m-\alpha+\mu},
}

utilizando la ecuación \eqref{eq:34} se obtiene

\formul{
\ifr{0}{D}{x}{\alpha}x^\mu &=& \frac{\gam{\mu+1}}{\gam{m-\alpha+\mu+1}}\frac{\gam{(m-\alpha+\mu)+1}}{\gam{(m-\alpha + \mu)-m+1}}x^{(m-\alpha+\mu)-m}\nonumber \\
&=& \frac{\gam{\mu+1}}{\gam{\mu-\alpha+1}}x^{\mu-\alpha}.
}{\label{eq:35}}

Utilizando la ecuación \eqref{eq:35} para el caso en que \form{f(x)=K}, con \form{K} una constante, se obtiene

\formula{
\ifr{0}{D}{x}{\alpha}K&=& K \ifr{0}{D}{x}{\alpha}x^0\\
&=& K\frac{\gam{0+1}}{\gam{0-\alpha+1}}x^{0-\alpha}\\
&=& K\frac{1}{\gam{1-\alpha}}x^{-\alpha},
}

para cualquier \form{\alpha>0}. Así, la derivada fraccionaria de una constante \form{K} es cero sólo para valores enteros positivos de \form{\alpha=n\in\mathds{N}} ya que \form{\gam{1-n}\to \infty}. Por otra parte, para cualquier \form{\alpha>0}, se tiene que \form{\ifr{0}{D}{x}{\alpha}f(x)\equiv 0} para \form{f(x)=x^{\alpha-k}} con \form{k\in\conj{1,2,\cdots (\lfloor\alpha \rfloor+1)}}.

Para el siguiente ejemplo considérese \form{\ifr{0}{I}{x}{\alpha}f} para \form{f(x)=\log x}

\formula{
\ifr{0}{I}{x}{\alpha}f(x)&=& \frac{1}{\gam{\alpha}}\int_0^x (x-u)^{\alpha-1}\log u du,
}

tomando el cambio de variable \form{u=x(1-v)} e integrando por partes se obtiene

\formula{
\ifr{0}{I}{x}{\alpha}f(x)&=& \frac{x^\alpha \log x}{\gam{\alpha}}\int_0^1 v^{\alpha-1}dv + \frac{x^\alpha}{\gam{\alpha}}\int_0^1 v^{\alpha-1}\log(1-v)dv\\
&=& \frac{x^\alpha}{\gam{\alpha+1}}\log x -\frac{x^\alpha}{\alpha \gam{\alpha}}\int_0^1 \log(1-v)d(1-v^\alpha)\\
&=& \frac{x^\alpha}{\gam{\alpha+1}}\left[\log x -(1-v^\alpha)\log(1-v)\Big{|}_0^1-\int_0^1 \frac{1-v^\alpha}{1-v}dv \right],
}

teniendo en cuenta la relación \cite{hilfer00}

\formula{
\int_0^1 \frac{v^x-v^y}{1-v}dv=\psi(y+1)-\psi(x+1),\ (\re{x},\re{y}>-1)
}

donde la función \form{\psi}, esta definida por

\formula{
\psi(x)=\frac{1}{[\gam{x}]}\frac{d}{dx}\gam{x},
}

además, satisface la  relación de recurrencia

\formula{
\psi(x+1)-\psi(x)=\frac{1}{x},
}

con \form{
-\psi(1)=\gamma= 0.5772157\cdots
}, con lo cual se obtiene

\formula{
\ifr{0}{I}{x}{\alpha}f(x)&=& \frac{x^\alpha}{\gam{\alpha+1}}\left[\log x -(1-v^\alpha)\log(1-v)\Big{|}_0^1-\int_0^1 \frac{1-v^\alpha}{1-v}dv \right]\\
&=& \frac{x^\alpha}{\gam{\alpha +1}}\left( \log x -\psi(\alpha+1)+\psi(1)\right).
}

Por lo tanto, para \form{m-1\leq \re{\alpha}<m},

\formula{
\ifr{0}{D}{x}{\alpha}\log x &=& \der{d}{x}{m}\ifr{0}{I}{x}{m-\alpha}\log x\\
&=& \der{d}{x}{m}\frac{x^{m-\alpha}}{\gam{m-\alpha+1}}\left(\log x -\psi(1-\alpha)-\gamma \right),
}

en el caso \form{\alpha=n\in \mathds{N}} la expresión a la derecha debe interpretarse como el caso límite en el que \form{\alpha \rightarrow n}. De hecho 

\formula{
\lim_{\alpha \to  n}\frac{\psi(1-\alpha)}{\gam{1-\alpha}}=(-1)^{-n}\gam{n},
}

la formula

\formula{\der{d}{x}{n}\log x = -\frac{\gam{n}}{(-x)^n},} 

se sigue fácilmente. Sin embargo, para \form{\alpha=-n\in \mathds{N}} se obtiene el resultado clásico \cite{hilfer00}

\formula{
\ifr{0}{D}{x}{-n}\log x &=& \ifr{0}{I}{x}{n}\log x \\
&=& \frac{x^n}{n!}\left[\log x -\sum_{k=1}^n \frac{1}{k}\right].
}

Para el ultimo ejemplo se toma la definición de Weyl, \form{m-1\leq \alpha < m, m\in \mathds{N},}

\formula{
\ifr{x}{D}{\infty}{\alpha}f(x)=(-1)^m\der{d}{x}{m}\ifr{x}{W}{\infty}{m-\alpha}f(x),
}

para \form{f(x)=\ex{-px},p>0}, tomando el cambio de variable \form{u-x=\frac{y}{p}},

\formula{
\ifr{x}{W}{\infty}{\alpha}\ex{-px}&=& \frac{1}{\gam{\alpha}}\int_x^\infty (u-x)^{\alpha -1}\ex{-pu}du\\
&=& \frac{\ex{-px}}{p^\alpha \gam{\alpha}}\int_0^\infty y^{\alpha -1}\ex{-y}dy\\
&=& \frac{\ex{-px}}{p^\alpha \gam{\alpha}}\gam{\alpha}\\
&=&\frac{\ex{-px}}{p^\alpha}, \ (\alpha>0),
}

obteniendo finalmente  para el caso \form{p>0},

\formula{
\ifr{x}{D}{\infty}{\alpha}\ex{-px}&=& (-1)^m \der{d}{x}{m}p^{-(m-\alpha)}\ex{-px}\\
&=& p^\alpha\ex{-px}, \ (m-1\leq \alpha < m).
}

\subsubsection{Problema de la Tautócrona}

Terminamos esta sección estudiando la primera aplicación del cálculo fraccional hecha por Abel en \form{1823} para resolver el problema de la tautócrona, i.e., el problema de determinar la forma de una curva de manera que el tiempo de descenso de una masa puntual sin fricción que se deslice por la curva bajo la acción de la gravedad sea independiente del punto de partida.

Sea una curva suave y una partícula de masa \form{m} la cual  parte del reposo en el punto \form{P_0=(x_0,y_0)}, asumiendo que la partícula se desliza hacia el origen sin fricción sobre la curva  bajo la acción de su propio peso en el punto \form{P=(x,y)} \cite{aguilar}

\begin{figure}[!ht]
\centering
\includegraphics[height=6cm,width=7cm]{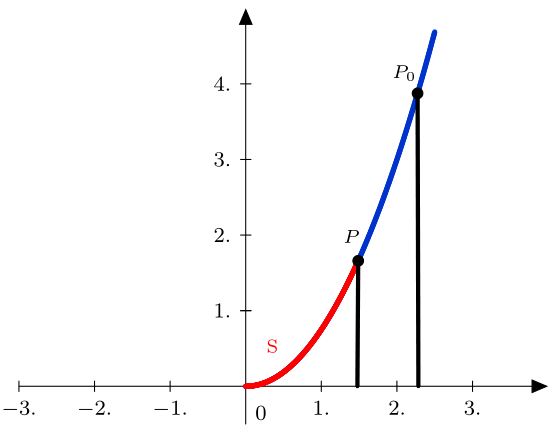}
\caption{Esquema para simplificar el  problema de la tautócrona}
\end{figure}

Partiendo del principio de conservación de la energía mecánica

\formula{
H(x,y)&=&K(x,y)+U(x,y),
}

donde \form{K(x,y)} y \form{U(x,y)} representan la energía cinética y potencial respectivamente,  debido a que el Hamiltoniano es una constante

\formula{
H(x,y)&=&H(x_0,y_0),
}

se puede obtener 

\formula{
K(x,y)+U(x,y)&=& K(x_0,y_0)+U(x_0,y_0)\\
 \frac{1}{2}m\left( \frac{d}{dt}x\right)^2+mgy &=&\frac{1}{2}m\left( \frac{d}{dt}x_0\right)^2 + mgy_0,\\
}

ya que la partícula parte del reposo 

\formula{
\frac{d}{dt}x_0=0,
}

tomando \form{S} como el segmento de  curva que recorre la partícula desde la posición inicial hasta el origen

\formula{
\frac{d}{dt}x\simeq\frac{d}{dt}S,
}

de lo anterior se tiene que

\formula{
\frac{1}{2}m\left( \frac{d}{dt}x\right)^2+mgy &=&\frac{1}{2}m\left( \frac{d}{dt}x_0\right)^2 + mgy_0\\
\frac{1}{2}m\left( \frac{d}{dt}S\right)^2+mgy &=& mgy_0\\
\frac{1}{2}\left(\frac{d}{dt}S \right)^2  &=& gy_0-gy,
}

despejando \form{\frac{d}{dt}S} de la relación anterior

\formula{
\left( \frac{d}{dt}S\right)^2 &=& 2g(y_0-y)\\
\frac{d}{dt}S&=&\pm \left( g(y_0-y)\right)^\frac{1}{2} ,
}

como consecuencia de que la partícula se acerca al origen conforme el tiempo \form{t} aumenta se tiene que la distancia \form{S} disminuye, lo que implica que \form{\frac{d}{dt}S<0}, debido a esto se toma la parte negativa de la raíz

\formula{
\frac{d}{dt}S &=& -\left( 2g(y_0-y)\right)^\frac{1}{2},
}

invirtiendo la expresión anterior y despejando \form{dt} se obtiene

\formula{
dt &=& -\left( 2g\right)^{-\frac{1}{2}} \left( y_0-y\right)^{-\frac{1}{2}} dS,
}

sea \form{T} el tiempo total que la partícula necesita para llegar desde el punto \form{P_0} hasta el origen \cite{aguilar}

\formula{
T &=& \int_{y_0}^0 dt,\\
&=& \int_{y_0}^0-\left( 2g\right)^{-\frac{1}{2}} \left( y_0-y\right)^{-\frac{1}{2}} dS\\
&=& \left( 2g\right)^{-\frac{1}{2}} \int_0^{y_0} \left( y_0-y\right)^{-\frac{1}{2}} dS.
}

Sea \form{\phi(y_0)=\left( 2g\right)^{\frac{1}{2}}T}, entonces 

\formula{
\phi(y_0)=\int_0^{y_0} \left( y_0-y\right)^{-\frac{1}{2}} dS,
}

haciendo tender \form{\frac{1}{2}\to n} (con \form{n\in[0,1]}) para tomar la expresión anterior de forma más general

\formula{
\phi(y_0)=\int_0^{y_0} \left( y_0-y\right)^{-n} dS,
}

Se la función Beta para los parámetros \form{\alpha} y  \form{\beta= 1-n} 

\formula{
B(\alpha,1-n)&=& \int_0^1 t^{\alpha -1}(1-t)^{(1-n)-1}dt\\
&=& \int_0^1 t^{\alpha-1} (1-t)^{-n} dt,
}

tomando el cambio de variable \form{z=y_0t} en la ecuación anterior

\formula{
B(\alpha,1-n) &=& \int_0^{y_0} (y_0^{-1}z)^{\alpha -1} (1-y_0^{-1}z)^{-n} (y_0^{-1}dz)\\
&=& y_0^{1-\alpha}y_0^{-1}y_0^{n} \int_0^{y_0} z^{\alpha-1} (y_0-z)^{-n} dz\\
&=& y_0^{n-\alpha}\int_0^{y_0} z^{\alpha -1} (y_0-z)^{-n}dz,
}

de lo anterior se obtiene 

\formula{
B(\alpha,1-n) y_0^{\alpha-n } &=&   \int_0^{y_0} z^{\alpha -1} (y_0-z)^{-n} dz ,
}

ahora multiplicando por \form{(y-y_0)^{n-1}} e integrando con respecto a \form{y_0}  desde \form{0} a \form{y}

\formul{
B(\alpha,1-n) \int_{0}^{y} y_0^{\alpha-n } (y-y_0)^{n-1} dy_0 &=&  \int_{0}^{y} (y-y_0)^{n-1} dy_0 \int_0^{y_0} z^{\alpha -1} (y_0-z)^{-n} dz,
}{\label{eq:36}}

tomando el cambio de variable \form{y_0=yt} en el lado izquierdo de la ecuación \eqref{eq:36}

\formula{
B(\alpha, 1-n)\int_0^y y_0^{\alpha-n} (y-y_0)^{n-1}  dy_0&=& B(\alpha ,1-n)\int_0^1 (yt)^{\alpha-n} (y-ty)^{n-1} ydt\\
&=& B(\alpha,1-n) \int_0^1 y^{n-1} y^{\alpha-n}y  t^{\alpha-n}(1-n)^{n-1} dt\\
&=& y^\alpha B(\alpha ,1-n)\int_0^1 t^{(\alpha+1-n)-1}(1-t)^{n-1}dt,
}

el cambio de variable anterior permite obtener una función Beta para los parámetros \form{(\alpha+1-n)} y \form{n}

\formula{
B(\alpha, 1-n)\int_0^y \frac{y_0^{\alpha-n}}{(y-y_0)^{1-n}} dy_0 &=& y^\alpha B(\alpha ,1-n)\int_0^1 t^{(\alpha+1-n)-1}(1-t)^{n-1}dt\\
&=& y^\alpha B(\alpha,1-n)B(\alpha+1-n,n)\\
&=& y^\alpha \left[ \frac{\gam{\alpha}\gam{1-n}}{\gam{\alpha +1 -n}}\right]\left[ \frac{\gam{\alpha+1-n} \gam{n}}{\gam{\alpha+1}}\right]\\
&=&y^\alpha \frac{\gam{\alpha}\gam{n}\gam{1-n}}{\gam{\alpha + 1}}\\
&=&y^\alpha \frac{\gam{\alpha}\gam{n}\gam{1-n}}{\alpha\gam{\alpha}}\\
&=& \frac{y^\alpha}{\alpha}\gam{n}\gam{1-n},
}

por otro lado se tiene que la función Gamma satisface la relación \cite{arfken85}

\formula{
\gam{n}\gam{1-n}=\frac{\pi}{\sin(n\pi)},
}

entonces el lado izquierdo de la ecuación \eqref{eq:36} toma la forma

\formula{
B(\alpha, 1-n)\int_0^y y_0^{\alpha-n} (y-y_0)^{n-1} dy_0&=&\gam{n}\gam{1-n} \frac{y^\alpha}{\alpha}\\
&=&\frac{\pi}{\sin(n\pi)} \frac{y^\alpha}{\alpha},
}

en consecuencia  la ecuación \eqref{eq:36} se reescribe como

\formula{
\frac{\pi}{\sin(n\pi)} \frac{y^\alpha}{\alpha} &=& \int_0^y (y-y_0)^{n-1}dy_0\int_0^{y_0}z^{\alpha-1} (y_0-z)^{-n} dz,
}

despejando \form{\frac{y^\alpha}{\alpha}} en la expresión anterior

\formula{
\frac{y^\alpha}{\alpha} &=& \frac{\sin(n\pi)}{\pi} \int_0^y (y-y_0)^{n-1}dy_0\int_0^{y_0}z^{\alpha-1} (y_0-z)^{-n} dz,
}

multiplicando por \form{\alpha \phi(\alpha)d\alpha} e integrando de forma indefinida se obtiene

\formul{
 \int y^\alpha \phi(\alpha)d\alpha
&=& \frac{\sin(n\pi)}{\pi} \int_0^y  (y-y_0)^{n-1}dy_0\int_0^{y_0} \left[\int z^{\alpha-1} \alpha \phi(\alpha)d\alpha \right] (y_0-z)^{-n} dz,
}{\label{eq:37}}

tomando \cite{aguilar}

\formula{
f(y)&=&\int  y^\alpha \phi(\alpha)d\alpha,\\
\frac{d}{dy}f(y)&=& \int  y^{\alpha-1}\alpha \phi(\alpha)d\alpha,
}

y sustituyendo en la ecuación \eqref{eq:37} se obtiene

\formul{
f(y)&=& \frac{\sin(n\pi)}{\pi} \int_0^y  (y-y_0)^{n-1}dy_0\int_0^{y_0} \left[\frac{d}{dz}f(z)\right] (y_0-z)^{-n} dz.
}{\label{eq:32}}

Recordando que

\formula{
\phi(y_0)=\int_0^{y_0} \left( y_0-y\right)^{-n} dS,
}

y multiplicando por \form{ \frac{\sin(n\pi)}{\pi}(y-y_0)^{n-1}
} e integrando con respecto a \form{y_0} desde \form{0} a \form{y}

\formul{
\frac{\sin(n\pi)}{\pi} \int_0^y   (y-y_0)^{n-1} \phi(y_0) dy_0 = \frac{\sin(n\pi)}{\pi} \int_0^y (y-y_0)^{n-1} dy_0 \int_0^{y_0} (y_0-y)^{-n} dS,
}{\label{eq:33}}

comparando las ecuaciones \eqref{eq:32} y \eqref{eq:33} se obtiene que

\formula{
S=\frac{\sin(n\pi)}{\pi} \int_0^y   (y-y_0)^{n-1} \phi(y_0) dy_0.
}

Para el caso particular en que \form{n=\frac{1}{2}} se obtiene

\formul{
\phi(y_0)&=&\int_0^{y_0} \left( y_0-y\right)^{-\frac{1}{2}} dS,\nonumber \\
S&=&\frac{1}{\pi} \int_0^y   (y-y_0)^{-\frac{1}{2}} \phi(y_0) dy_0,
}{\label{eq:38}}

cuando el tiempo de deslizamiento es una constante conocida, la ecuación integral de Abel se obtiene multiplicando la ecuación \eqref{eq:38} por \form{\pi}

\formula{
\pi S= \int_0^y   (y-y_0)^{-\frac{1}{2}} \phi(y_0) dy_0,
}

la integral en la ecuación anterior es, excepto por el factor multiplicativo \form{\gam{\frac{1}{2}}^{-1}}, un caso particular de una integral fraccionaria de orden \form{\frac{1}{2}} \cite{miller93}.


\newpage

\bibliographystyle{unsrt}
\bibliography{biblio}

\end{document}